\documentclass{siamltex}
\usepackage{amsmath}
\usepackage{amssymb}
\usepackage{dsfont}
\usepackage{color}
\usepackage{graphicx}
\usepackage{subfig}
\usepackage{multirow}
\usepackage{morefloats}
\usepackage{pgfplots}
\usetikzlibrary{plotmarks}
\usepackage{tikz}
\usepackage{epstopdf}

\begin{document}
\title{A meshfree particle method for a vision-based macroscopic pedestrian model}

\author{N.K. Mahato \footnotemark[1]\and A. Klar\footnotemark[1] \footnotemark[2]
       \and S. Tiwari \footnotemark[1] }
\footnotetext[1]{Technische Universit\"at Kaiserslautern, Department of Mathematics, Erwin-Schr\"odinger-Stra{\ss}e, 67663 Kaiserslautern, Germany 
  (\{mahato,klar,tiwari\}@mathematik.uni-kl.de)}
\footnotetext[2]{Fraunhofer ITWM, Fraunhoferplatz 1, 67663 Kaiserslautern, Germany}

\maketitle
\begin{abstract}
\noindent In this paper we present numerical simulations of a macroscopic vision-based model \cite{degond} derived from microscopic situation rules described in \cite{ondrej}. This model describes an approach to collision  avoidance between pedestrians by taking decisions of turning or slowing down based on  basic interaction rules, where the dangerousness level of an interaction with another pedestrian is measured in terms of the derivative of the bearing angle  and of the time-to-interaction. A meshfree particle method is used to solve the equations of the model. Several numerical cases are considered to compare this model with models established in the field, for example, social force model coupled to an Eikonal equation \cite{etikyala}. Particular emphasis is put on the comparison of evacuation and computation times.
\end{abstract}



\section{Introduction}
To achieve a better comprehension of crowd behavior and increase the reliability of predictions, numerical modeling and simulation of human crowd motion have become a major subject of research with a wide field of applications. These include planning the architecture of buildings, supporting transport planners or managers to design timetables, improving efficiency and safety in public places such as airport terminals, train stations, theaters and shopping malls.

Several models of pedestrian behavior have been proposed on different levels of description in recent years. The Individual-Based Models (or IBM) such as models based on heuristic rules \cite{reynolds}, Newton type equations \cite{helbing,helbing2,maury,piccoli}, traffic following models \cite{lemer}, cellular-automata \cite{B01,luo}, agent-based models \cite{karb,omer}, optimal control theory models \cite{hoog}, mechanical models \cite{helbing4,helbing5} and vision-based models \cite{ondrej,Rhughes,guy,huang,paris,pettre,berg} have been developed. Macroscopic pedestrian flow model involving equation for density and mean velocity of the flow are derived in Refs. \cite{degond,etikyala,mahato,degond2,bello,helbing3,treuille}.

Pedestrians can estimate the positions and velocities of moving obstacles such as other pedestrians  with fairly good accuracy, compare   \cite{cutting} and are able to process this information in order to determine the dangerousness level of an encounter \cite{warren}. With these considerations, the heuristic-based model of \cite{mous} proposes that pedestrians follow two phases: a perception phase and a decision-making phase.

In the perception phase, pedestrians make an assessment of the dangerousness of an encounter with the other pedestrians. In the decision-making phase, they take decisions of turning or slowing down based on this assessment. To make this assessment, each pedestrian senses the indicators of the dangerousness of the collision such as the derivative of the Bearing Angle (DBA), the Time-To-Interaction (TTI) and the minimal distance (MD) with each of his interaction partners. The Bearing angle is the angle between the walking direction and the line connecting the two interaction partners and its time derivative provides a sensor of the likeliness of a collision, a constant Bearing angle in time being associated with a risk of future collision. The TTI is the time needed by the subject to reach the interaction  point from his current position. The MD is the distance where the distance between the collision subjects is minimal. Finally, the risk of future collisions is determined by combining these indicators DBA, TTI and MD. For detailed explanations and geometrical interpretations we refer \cite{degond}. 

Using these considerations,  a macroscopic model has been derived in \cite{degond} from the vision-based Individual-Based Models (IBM) of \cite{ondrej}. Moreover, there, a procedure has been presented to approximate the macroscopic  mean-field  equations with a local approximation. These models will be investigated  numerically in the present paper. 

One main objective of the present paper is to develop a numerical scheme to simulate the vision-based pedestrian model on a macroscopic scale  and to compare it to other models established in the field, for example, classical social force models \cite{helbing,etikyala}. We use, as in Ref. \cite{mahato}, a particle method which is based on a Lagrangian formulation of these equations, where particles are used as moving grid points. Preliminary numerical comparisons of the macroscopic models are presented. Moreover, we note that the method presented here is easily extended to more complicated ``real" life situations.

The reminder of the paper is organized as follows. Section \ref{sec:2} contains the macroscopic vision-based model for pedestrian flow from \cite{degond}. The local approximation of the macroscopic model is presented in section \ref{sec:3}. Particle method used in the simulations is described in section \ref{sec:4}. Section \ref{sec:5} contains the numerical results. We consider a bidirectional evacuation problem. A comparison between the solutions of the macroscopic equations from the vision-based model and social force model is presented together with a comparison of the associated evacuation as well as computational times. Finally, section \ref{sec:6} concludes the work.

\section{Vision-based macroscopic pedestrian flow model from \cite{degond}}\label{sec:2}
For completeness of the presentation we  review  the  macroscopic vision-based model for pedestrian flow, see  Refs.\cite{degond,ondrej}.
For further details we refer to these papers.  There, denoting  the density of pedestrians by $\rho(x, \xi, t)$ and by $U(x, \xi, t)$  the mean velocity at position $x$ with target point $\xi$ at time $t$, the model is given by the following mass and momentum balance equations
\begin{align}
	\partial_t \rho + \nabla_x \cdot (c\rho U) &= 0, \label{eq:1} \\
	\partial_t U + cU \cdot \nabla_x U & =  \omega_{\rho \delta_U}(x,U(x,\xi,t),\xi,t) U^\perp(x,\xi,t), \label{eq:2}
\end{align}
where $c$ is the constant speed of the pedestrians and $\omega_{\rho\delta_U}(x,U(x,\xi,t),\xi,t)$  is the angular velocity, defined in  the following way.
First one has to define the indicators of the dangerousness of the collision measured between a particle located at position $x\in \mathbb{R}^2$ with velocity direction $u \in \mathbb{S}^1$  to the relative target $\xi \in \mathbb{R}^2$ and its collision particle located at position $y \in \mathbb{R}^2$ with velocity direction $v\in\mathbb{S}^1$ to the relative target $\eta \in \mathbb{R}^2$.
They are  respectively denoted by $\tilde{\dot{\alpha}}(x, u, y, v)$, $\tilde{\tau}(x, u, y, v)$, $\tilde{D}(x, u, y, v)$ and $\tilde{\dot{\alpha}}_g(x, u, \xi)$ and defined as
\begin{align}
	\tilde{\dot{\alpha}}(x,\xi,y,\eta,t) &= c\dfrac{(y-x)\times (U(y,\eta,t)-U(x,\xi,t))}{\arrowvert y - x \arrowvert^2}, \label{eq:3}\\
	\tilde{\tau}(x,\xi,y,\eta,t) &= - \frac{1}{c}\dfrac{(y-x)\cdot (U(y,\eta,t)-U(x,\xi,t))}{\arrowvert U(y,\eta,t)-U(x,\xi,t) \arrowvert^2},\label{eq:4}\\
	\tilde{D}(x,\xi,y,\eta,t) &= \left(\arrowvert y - x \arrowvert^2 - \left( (y - x) \cdot \dfrac{U(y,\eta,t)-U(x,\xi,t)}{\arrowvert U(y,\eta,t)-U(x,\xi,t) \arrowvert}\right)^2\right)^{1/2}, \label{eq:5}\\ 
	\tilde{\dot{\alpha}}_g(x,u,\xi) &= - c\dfrac{(\xi-x)\times U(x,\xi,t)}{\arrowvert \xi - x \arrowvert^2}. \label{eq:6}
\end{align}
Moreover, one defines 
\begin{equation}\label{eq:7}
	\begin{aligned}
		\tilde{S}_{\pm}(x,\xi,t) = &\{(y,\eta)\in \mathbb{R}^2 \times \mathbb{R}^2 |\pm\tilde{\dot{\alpha}}(x,\xi,y,\eta,t)>0,  \tilde{\tau}(x,\xi,y,\eta,t)>0,\\
		&\hspace{9mm} \tilde{D}^2(x,\xi,y,\eta,t)<R^2, |\tilde{\dot{\alpha}}(x,\xi,y,\eta,t)|< \sigma(|\tilde{\tau}(x,\xi,y,\eta,t)|) \}
	\end{aligned}
\end{equation}
and
\begin{equation}\label{eq:8}
	\begin{aligned}
		\tilde{\Phi}_{\pm}&(x,\xi,t) \\
		&= \dfrac{\int_{(y,\eta)\in \tilde{S}_{\pm}(x,\xi,t)} \Phi(|\tilde{\dot{\alpha}}(x,\xi,y,\eta,t)|,|\tilde{\tau}(x,\xi,y,\eta,t)|)\rho(y,\eta,t)dyd\eta}{\int_{(y,\eta)\in \tilde{S}_{\pm}(x,\xi,t)}\rho(y,\eta,t)dyd\eta},
	\end{aligned}
\end{equation}
where the function $\Phi$ is given by
\begin{align}\label{eq:9}
	\Phi(|\tilde{\dot{\alpha}}(x,\xi,y,\eta,t)|,&|\tilde{\tau}(x,\xi,y,\eta,t)|) \nonumber\\
	&= \Phi_0 \max \{\sigma(|\tilde{\tau}(x,\xi,y,\eta,t)|)-|\tilde{\dot{\alpha}}(x,\xi,y,\eta,t)|, 0\}
\end{align}
and
\begin{equation}\label{eq:10}
	\sigma(|\tilde{\tau}(x,\xi,y,\eta,t)|) = a + \dfrac{b}{(|\tilde{\tau}(x,\xi,y,\eta,t)| + \tau_0)^c},
\end{equation}
where $\tau_0, a, b$ are positive constants and $\Phi_0$ is the positive proportionality constant. 

\noindent Then $\omega_{\rho\delta_U}(x,U(x,\xi,t),\xi,t)$  is defined in the following:
\begin{enumerate}
	\item If the current deviation to the goal is small, this means that $|\tilde{\dot{\alpha}}_{g}|$ is smaller than the reaction induced by collision avoidance, i.e.
	\begin{equation}\label{eq:11}
		- \tilde{\Phi}_{-}(x,\xi,t) \leq \tilde{\dot{\alpha}}_{g}(x,\xi,t) \leq \tilde{\Phi}_{+}(x,\xi,t).
	\end{equation}
	Then
	\begin{equation}\label{eq:12}
		\begin{aligned}
			&\omega_{\rho\delta_U}(x,U(x,\xi,t),\xi,t) =\\ 
			& -\tilde{\Phi}_{+}(x,\xi,t)H(|\tilde{\Phi}_{-}(x,\xi,t) - |\tilde{\dot{\alpha}}_{g}(x,\xi,t)|| - |\tilde{\Phi}_{+}(x,\xi,t) - |\tilde{\dot{\alpha}}_{g}(x,\xi,t)||)\\
			& + \tilde{\Phi}_{-}(x,\xi,t)H(|\tilde{\Phi}_{+}(x,\xi,t) - |\tilde{\dot{\alpha}}_{g}(x,\xi,t)|| - |\tilde{\Phi}_{-}(x,\xi,t) - |\tilde{\dot{\alpha}}_{g}(x,\xi,t)||),
		\end{aligned}
	\end{equation}
	
	where H is the Heaviside function. This formula states that the pedestrian determines the worst case in each direction (expression of $\tilde{\Phi}_+$ and $\tilde{\Phi}_-$) and then chooses the turning direction as the one which produces the smallest deviation to the goal.  
	\item If the deviation to the goal is large, i.e.
	\begin{equation}\label{eq:13}
		\tilde{\dot{\alpha}}_{g}(x,\xi,t) < - \tilde{\Phi}_{-}(x,\xi,t) \quad\text{or}\quad \tilde{\dot{\alpha}}_{g}(x,\xi,t) > \tilde{\Phi}_{+}(x,\xi,t).
	\end{equation}
	
	Then
	\begin{equation}\label{eq:14}
		\omega_{\rho\delta_U}(x,U(x,\xi,t),\xi,t) = \tilde{\dot{\alpha}}_{g}(x,\xi,t).
	\end{equation}
	In this case, the deviation to the goal is larger than the reaction to collisions and the decision is to restore a direction of motion more compatible to the goal.	
\end{enumerate}

\noindent The intensity of the force is non-local in space as the pedestrian anticipates the behavior of neighboring pedestrians to take a decision. In the next section  the spatially local approximation is discussed .


\section{Local approximations to the macroscopic model}\label{sec:3}
In this section we state  the spatially local approximation to the macroscopic model using  the procedure explained in  \cite{degond} for the local approximation to the mean-field kinetic model. One assumes that there exists a  small dimensionless quantity $\lambda \ll 1 $ such that
\begin{equation}\label{eq:15}
	R = \lambda \hat{R}, \hspace{2mm} \Phi_0 = \lambda \hat{\Phi_0},\hspace{2mm} a = \frac{1}{\lambda} \hat{a},\hspace{2mm} \tau_0 = \lambda \hat{\tau}_0,\hspace{2mm} b = \lambda^{c-1} \hat{b},\hspace{2mm} \sigma = \frac{1}{\lambda} \hat{\sigma},
\end{equation}
where all `hat' quantities are assumed to be $\mathcal{O}(1)$. One introduces a  change of variables $y = x + \lambda \zeta$, with $\zeta \in \mathbb{R}^2$ in all expressions involving $y$. This yields the expression for equations (\ref{eq:3}) - (\ref{eq:5}) as 
\begin{align*}
	\tilde{\dot{\alpha}}(x,\xi,y,\eta,t) &= \frac{1}{\lambda} \hat{\dot{\alpha}}(\zeta,U(y,\eta,t)-U(x,\xi,t)), \hspace{2mm}\mbox{where}\\ 
	&\hat{\dot{\alpha}}(\zeta,U(y,\eta,t)-U(x,\xi,t)) = c\dfrac{\zeta \times (U(y,\eta,t)-U(x,\xi,t))}{\arrowvert \zeta \arrowvert^2}, \\
	\tilde{\tau}(x,\xi,y,\eta,t) &= \lambda \hat{\dot{\tau}}(\zeta,U(y,\eta,t)-U(x,\xi,t)), \hspace{2mm}\mbox{where}\\
	&\hat{\dot{\tau}}(\zeta,U(y,\eta,t)-U(x,\xi,t)) = - \frac{1}{c}\dfrac{\zeta\cdot (U(y,\eta,t)-U(x,\xi,t))}{\arrowvert U(y,\eta,t)-U(x,\xi,t) \arrowvert^2},\\
	\tilde{D}(x,\xi,y,\eta,t) &= \lambda \hat{D}(\zeta,U(y,\eta,t)-U(x,\xi,t)), \hspace{2mm}\mbox{where}\\
	\hat{D}(\zeta, & U(y,\eta,t)-U(x,\xi,t)) = \left(\arrowvert \zeta \arrowvert^2 - \left( \zeta \cdot \dfrac{U(y,\eta,t)-U(x,\xi,t)}{\arrowvert U(y,\eta,t)-U(x,\xi,t) \arrowvert}\right)^2\right)^{1/2}. 
\end{align*}
The function $\Phi(|\tilde{\dot{\alpha}}|,|\tilde{\tau}|)$ is changed into $\hat{\Phi}(|\hat{\dot{\alpha}}|,|\hat{\tau}|)$, such that
\begin{align*}
	\Phi(|\tilde{\dot{\alpha}}(x,\xi,y,\eta,t&)|,|\tilde{\tau}(x,\xi,y,\eta,t)|)\\
	= & \hat{\Phi}(|\hat{\dot{\alpha}}(\zeta,U(y,\eta,t)-U(x,\xi,t))|,|\hat{\tau}(\zeta,U(y,\eta,t)-U(x,\xi,t))|),\\
	\hat{\Phi}(|\hat{\dot{\alpha}}|,|\hat{\tau}|)  = \hat{\Phi}_0 & \max \{\hat{\sigma}(|\hat{\tau}|) -  |\hat{\dot{\alpha}}|, 0\},\\
	\hat{\sigma}(|\hat{\tau}|)  = \hat{a} + & \dfrac{\hat{b}}{(|\hat{\tau}|+\hat{\tau}_0)^c} .
\end{align*}
Thus, formula (\ref{eq:8})  can be written as:
\begin{equation}\label{eq:16}
	\tilde{\Phi}_{\pm}(x,\xi,t)
	= \dfrac{\int_{(\zeta,\eta)\in \tilde{S}_{\pm}(x,\xi,t)} \hat{\Phi}(|\hat{\dot{\alpha}}|,|\hat{\tau}|) \rho(x+\lambda \zeta,\eta,t)d\zeta d\eta}{\int_{(\zeta,\eta)\in \tilde{S}_{\pm}(x,\xi,t)}\rho(x+\lambda\zeta,\eta,t)d\zeta d\eta},
\end{equation}
where
\begin{equation*}
	\begin{aligned}
		\tilde{S}_{\pm}(x,\xi,t) & = \{(\zeta,\eta)\in \mathbb{R}^2 \times \mathbb{R}^2 |\pm\hat{\dot{\alpha}}(\zeta,U(y,\eta,t)-U(x,\xi,t))>0, \\ 
		&\hat{\tau}(\zeta,U(y,\eta,t)-U(x,\xi,t))>0, \hat{D}^2(\zeta,U(y,\eta,t)-U(x,\xi,t))<\hat{R}^2,\\
		& \hspace{8mm} |\hat{\dot{\alpha}}(\zeta,U(y,\eta,t)-U(x,\xi,t))|< \hat{\sigma}(|\hat{\tau}(\zeta,U(y,\eta,t)-U(x,\xi,t))|) \}.
	\end{aligned}
\end{equation*}
For $\lambda \longrightarrow 0$  the dependence of $\rho$ upon $\zeta$ disappears and $\hat{\Phi}$ can be integrated out with respect to $\zeta$. Therefore, (\ref{eq:16}) leads to:
\begin{equation}\label{eq:17}
	\tilde{\Phi}_{\pm}(x,\xi,t)
	= \dfrac{\int_{\eta\in \mathbb{R}^2} \Psi_{\pm}(|U(x,\eta,t)-U(x,\xi,t)|) \rho(x,\eta,t) d\eta}{\int_{\eta\in \mathbb{R}^2}\rho(x,\eta,t) d\eta},
\end{equation}
where
\begin{equation*}
	\Psi_{\pm}(|U(x,\eta,t)-U(x,\xi,t)|) = \dfrac{\int_{\zeta\in \hat{S}_{\pm}(U(x,\eta,t)-U(x,\xi,t))} \hat{\Phi}(|\hat{\dot{\alpha}}|,|\hat{\tau}|) d\zeta }{\text{Area}(\hat{S}_{\pm}(U(x,\eta,t)-U(x,\xi,t)))}
\end{equation*}
and the set
\begin{equation*}
	\begin{aligned}
		\hat{S}_{\pm}(U(x,\eta,t)& -U(x,\xi,t)))  = \{\zeta\in\mathbb{R}^2|\pm\hat{\dot{\alpha}}(\zeta,U(y,\eta,t)-U(x,\xi,t))>0, \\ 
		&\hat{\tau}(\zeta,U(y,\eta,t)-U(x,\xi,t))>0, \hat{D}^2(\zeta,U(y,\eta,t)-U(x,\xi,t))<\hat{R}^2,\\
		& \hspace{8mm} |\hat{\dot{\alpha}}(\zeta,U(y,\eta,t)-U(x,\xi,t))|< \hat{\sigma}(|\hat{\tau}(\zeta,U(y,\eta,t)-U(x,\xi,t))|) \}.
	\end{aligned}
\end{equation*}
After the computation of the function $\tilde{\Phi}_{\pm}$, one  determines $\omega_{\rho \delta_U}$ as in  section \ref{sec:2}.

The approximation (\ref{eq:17}) is called ``local" since only the values of $\rho$ at position $x$ are needed to evaluate this integral.

\section{Numerical method}\label{sec:4}
In this section, we discuss the numerical method for the macroscopic vision-based pedestrian models (\ref{eq:1}) and (\ref{eq:2}). To solve these equations numerically we use a meshfree particle method, see, for example \cite{fpm}. Meshfree particle methods are an appropriate scheme to solve pedestrian flow problems due to the appearance of situations with complicated geometries, free and moving boundaries and potentially large deformations of the domain of computation, i.e. the region where the density of pedestrians is non-zero. The particle method is based on a Lagrangian formulation of equations (\ref{eq:1}) and (\ref{eq:2}), which is given as
\begin{eqnarray}
	\frac{dx}{dt} &=& cU \nonumber \\
	\frac{d\rho}{dt} &=& -c \rho \nabla_x \cdot U \nonumber \\
	\frac{dU}{dt} &=& \omega_{\rho \delta_U}(x,U(x,\xi,t),\xi,t) U^\perp(x,\xi,t),\nonumber
\end{eqnarray}
where $\frac{d}{dt} = \partial_t + cU\cdot\nabla_x$.

A meshfree Lagrangian method is used to approximate the spatial differential and integral operators appearing on the right hand side of the above  equations, where a difference approximation at the particle locations from the surrounding neighboring particles a weighted least square approximation has been used. For the present computation we use a weight function $w$ with compact support of radius h, restricting in this way the number of neighboring particles. 
The Gaussian weight functions are of the form
\begin{eqnarray} 
	w = w(  x; h) =
	\left\{ 
	\begin{array}{l}  
		\exp (- \alpha \frac{ \vert x  \vert^2 }{h^2} ), 
		\quad \mbox{if    }  \frac{|  x  |}{h} \le 1 
		\\
		0,  \qquad \qquad \quad \quad \quad \mbox{else},
	\end{array}
	\right.
	\label{weight}
\end{eqnarray}
where $\alpha > 0$, a user defined constant is taken in the range of $2$ to $6$.   
The radius h is chosen to include initially enough particles for a stable approximation of the spatial derivatives, which is approximately three times the initial spacing of the particles. For details of the implementation, we refer to \cite{tiwari}.

The simplest way to evaluate the integral over the interaction potential is to use  a straightforward first order integration rule using an approximation of the local area around a particle determined by nearest neighbor search. This works fine for a well resolved situation with a sufficiently large number of grid points. After the approximations of the spatial derivatives and the integral, we obtain a system of time dependent ODEs. The resulting system of ODEs is then solved  by a suitable  time discretization method. 


\section{Numerical Results}\label{sec:5}

\begin{figure}
	\includegraphics[scale=0.25]{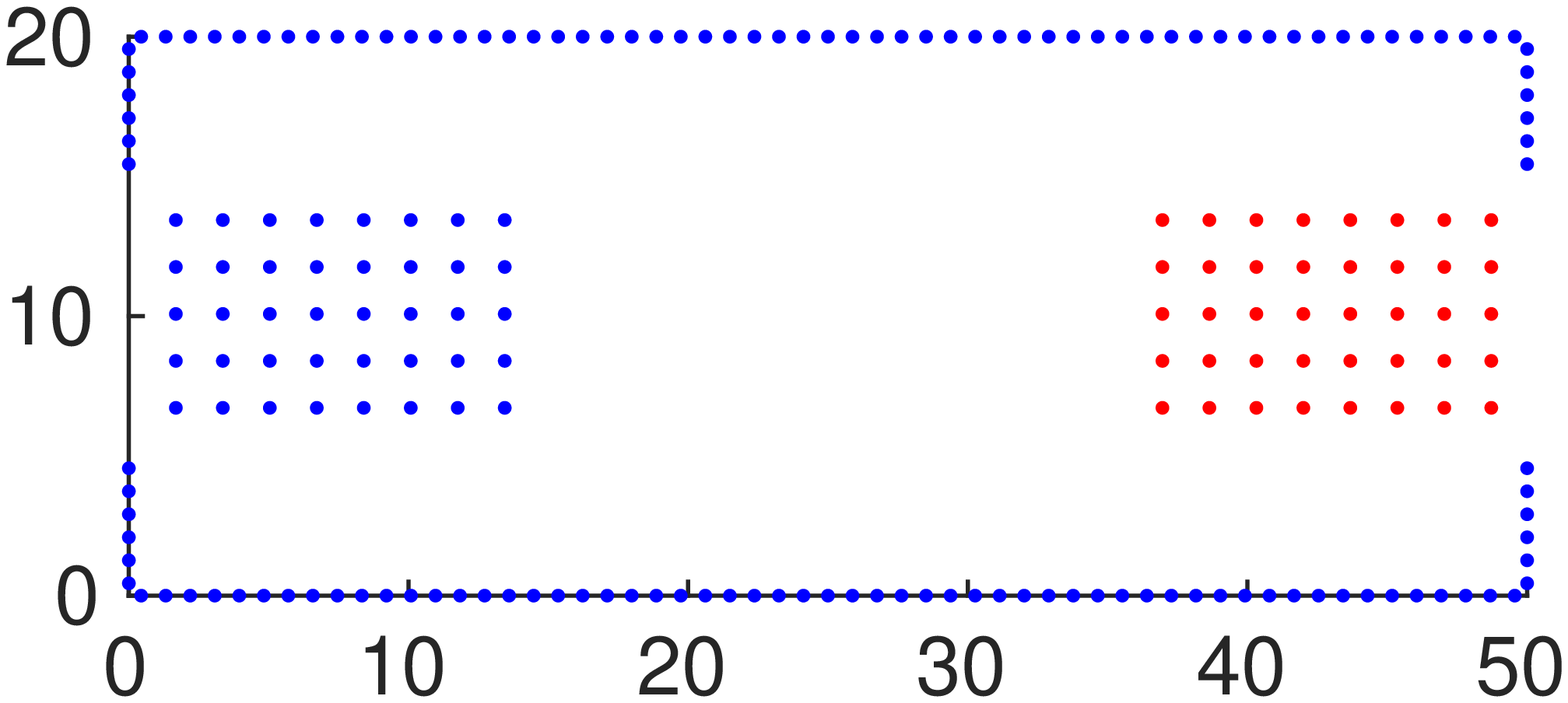}
	\includegraphics[scale=0.25]{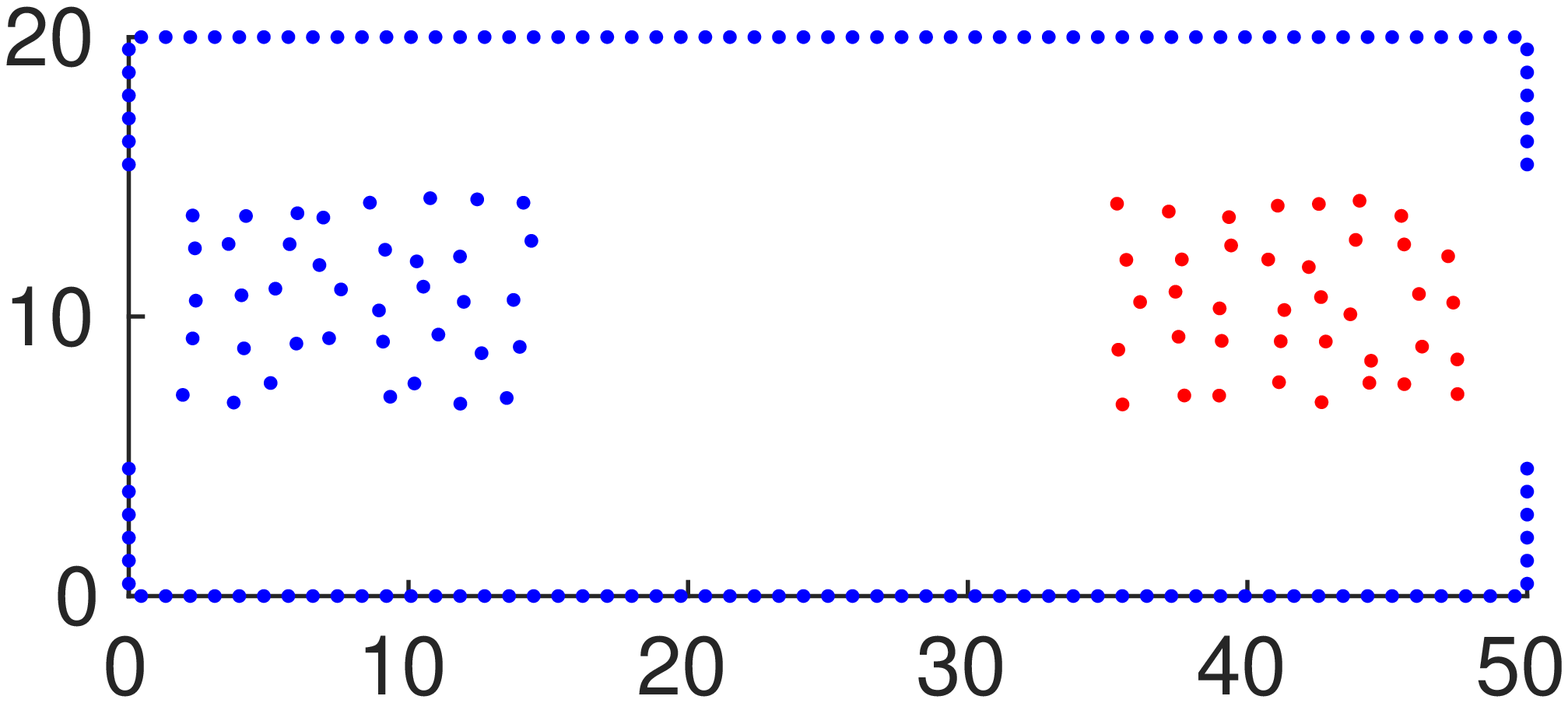}
	\caption{Initial distribution of grid particles for the systematic distribution case (left) and the random distribution case (right).}
	\label{fig:1} 
\end{figure}

In this section, we present a series of numerical experiments for the macroscopic vision-based pedestrian model  equations (\ref{eq:1}), (\ref{eq:2}) applied to bidirectional evacuation problem. We investigate the model numerically for a corridor of  length $50$ and  width $20$ with two entrances and two exits of width $10$ on both ends as in figure \ref{fig:1}. All distance are measured in meter $m$. Densities are measured in $1/m^2$. Time is measured in seconds $s$ and velocity in $m/s$. Initially, pedestrians are concentrated at the left boundary (blue particles) and right boundary (red particles). Pedestrians on the left boundary can leave through right exit and on right boundary through left exit. We choose parameters as $a = 0$, $b = 0.6$, $c = 1.5$, $\tau_0 = 1.0$ and $\Phi_0 = 1.0$. We are choosing different value for $R$ depending on the initial spacing of the particles. 

\subsection{Understanding decision-making phase}\label{sec:5.1}  
In this subsection we consider systematic and random initial distribution of 80 grid points (40 grid points  on each end) as in figure \ref{fig:1}. Pedestrians at left entrance have the destination point  (50,10). The pedestrians entering at right entrance have (0,10) on the opposite end. Here we try to understand collision avoidance between particles by the comparison of without control on direction and with control on direction. In the case of without control on direction, pedestrians have only the information to satisfy goal (i.e. the angular velocity is given by (\ref{eq:14}) if the deviation to the goal is large, otherwise, zero). We use the constant time step $\Delta t = 0.00042$, initial spacing of particles $\Delta x = 1.68$ and $R = 1.68$.
\begin{figure}
	\includegraphics[width=.5\linewidth]{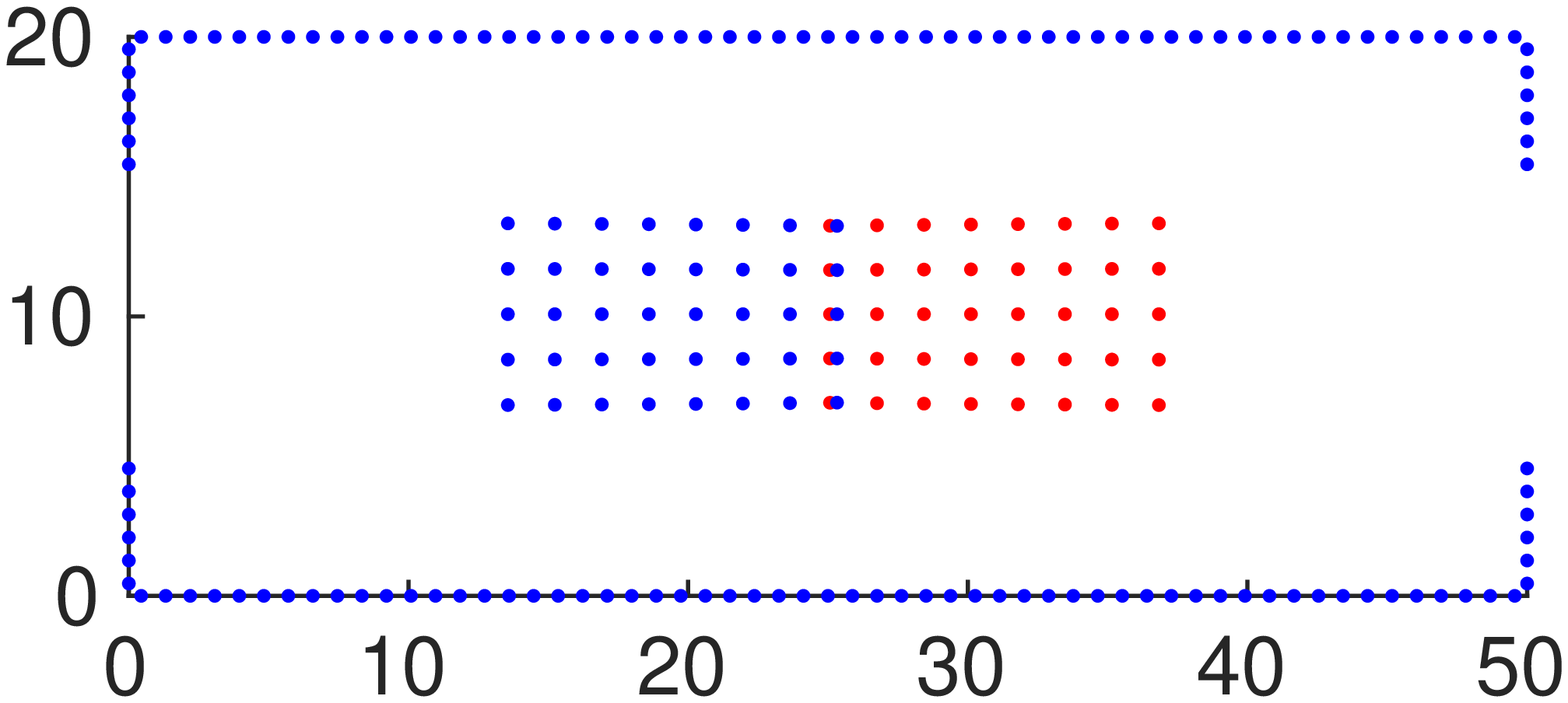}\hfill
	\includegraphics[width=.5\linewidth]{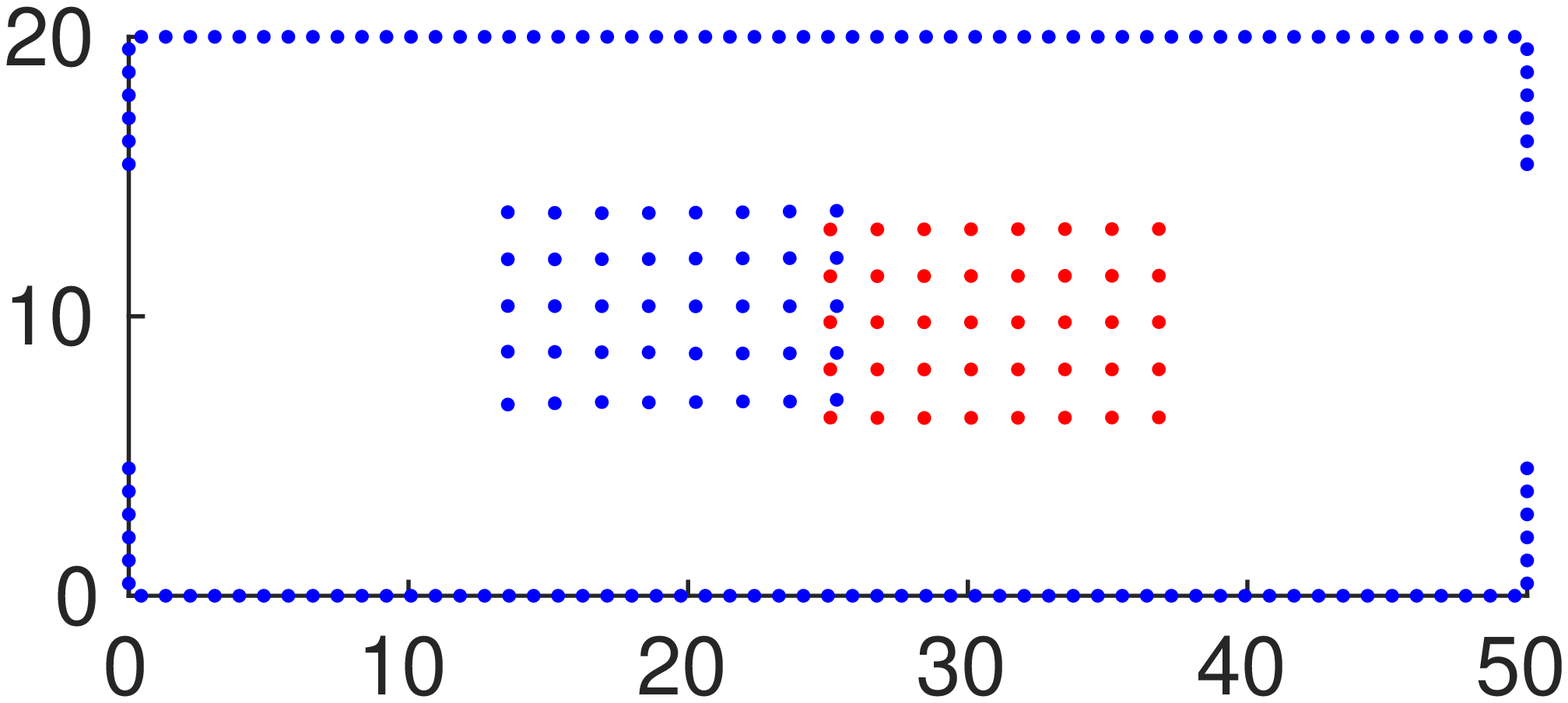}\\
	\includegraphics[width=.5\linewidth]{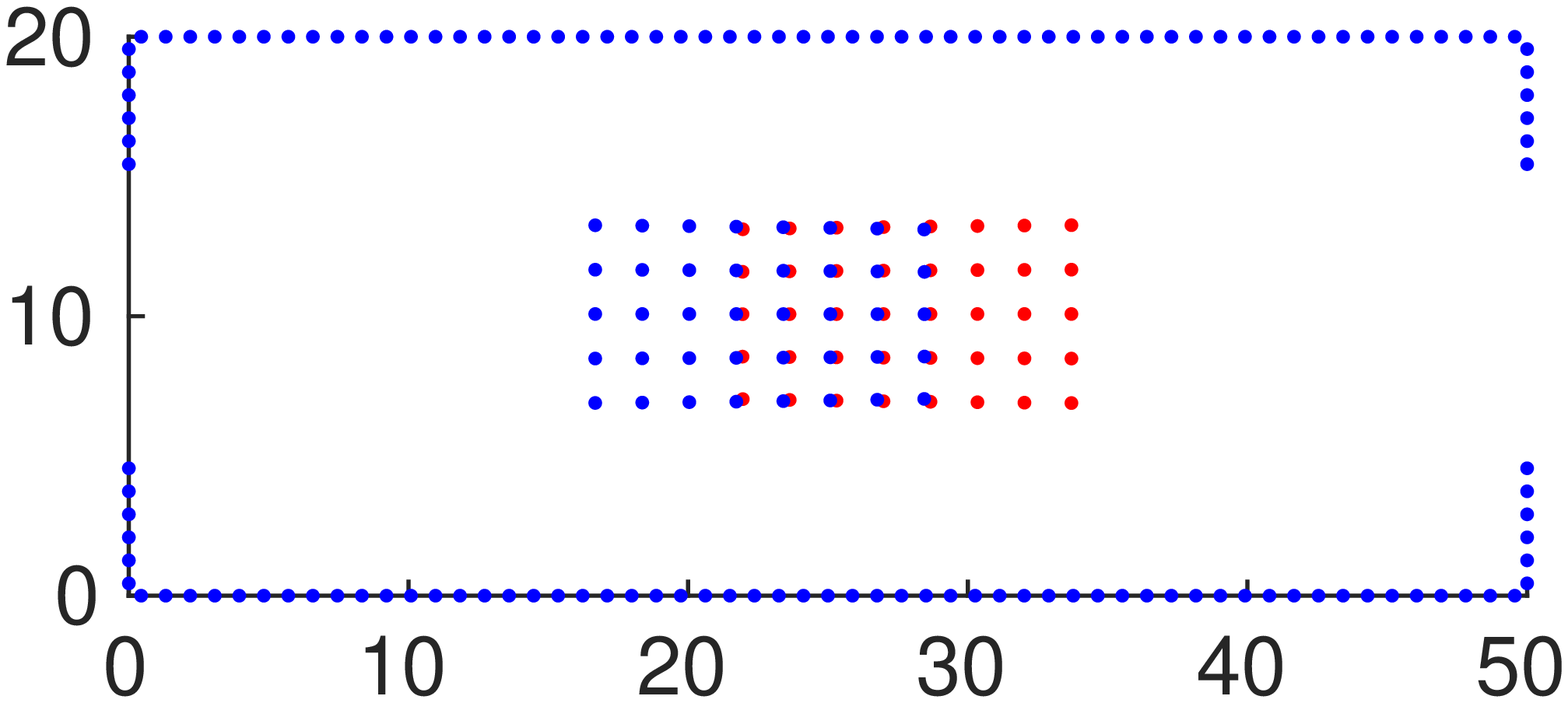}\hfill
	\includegraphics[width=.5\linewidth]{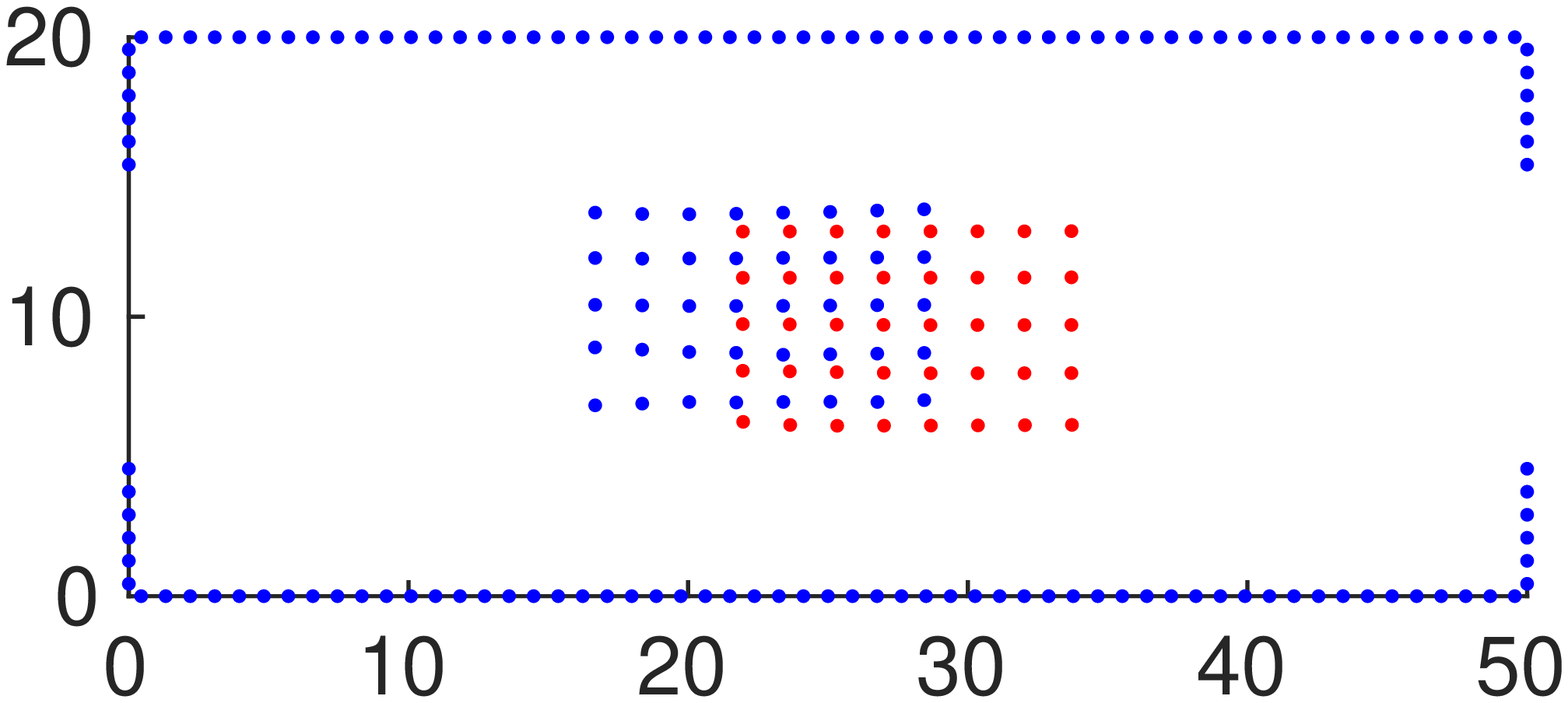}\\
	\includegraphics[width=.5\linewidth]{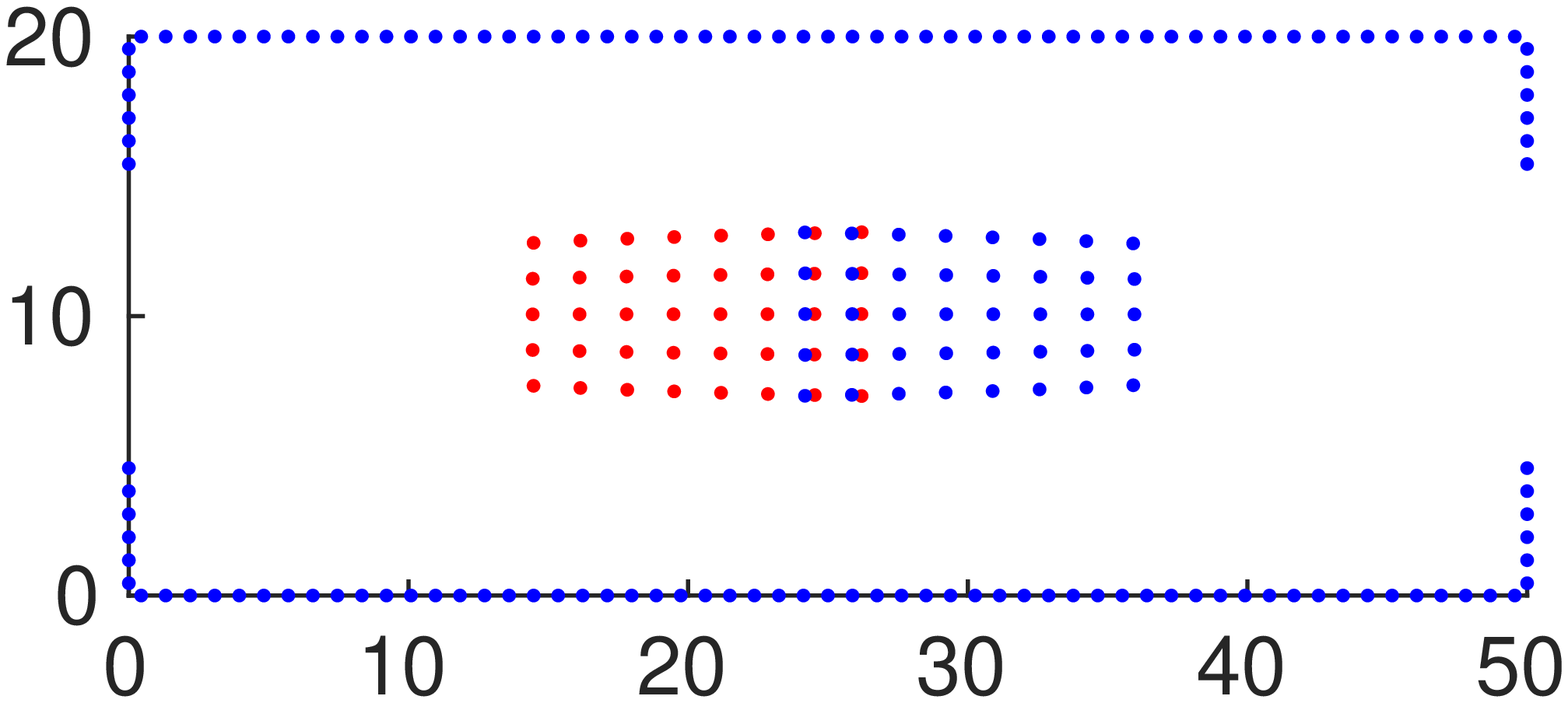}\hfill
	\includegraphics[width=.5\linewidth]{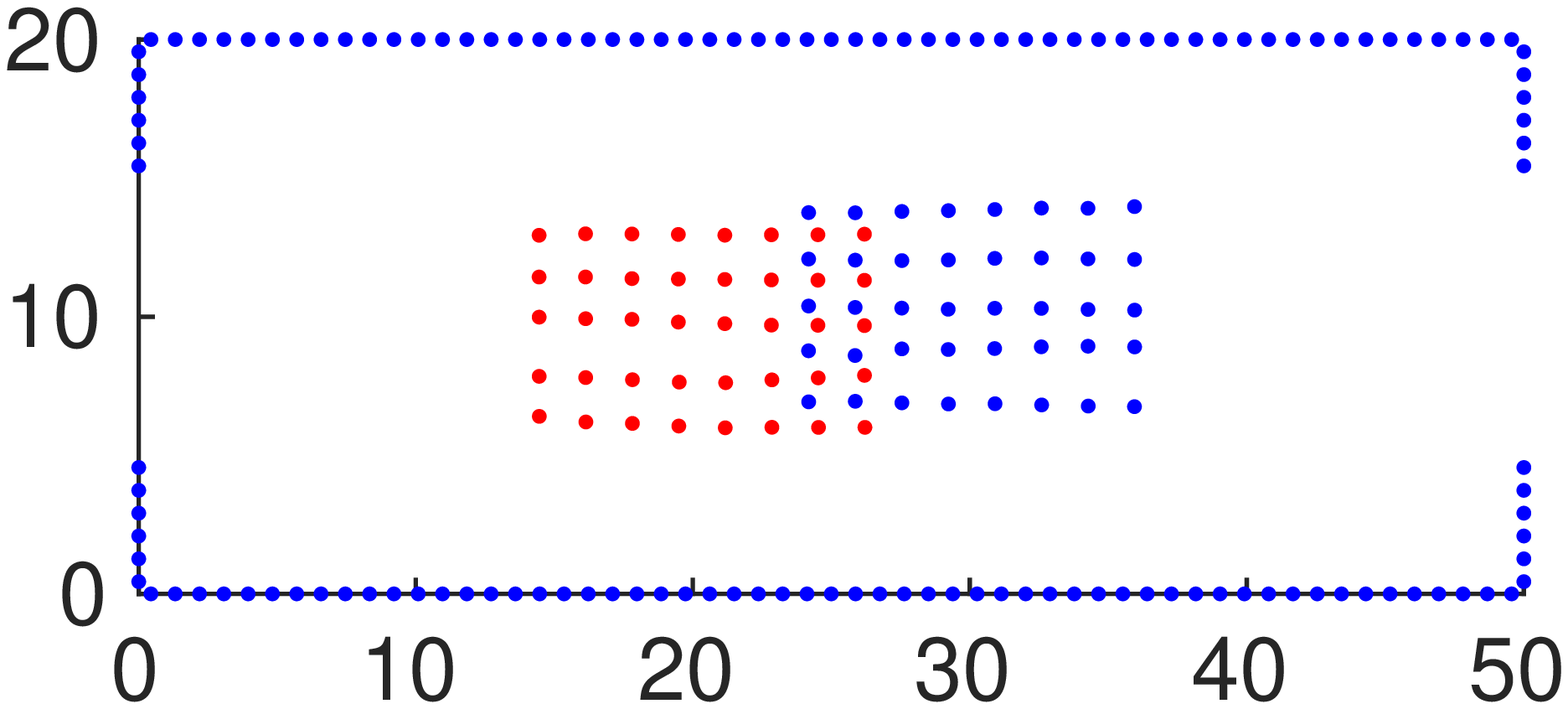}\\
	\includegraphics[width=.5\linewidth]{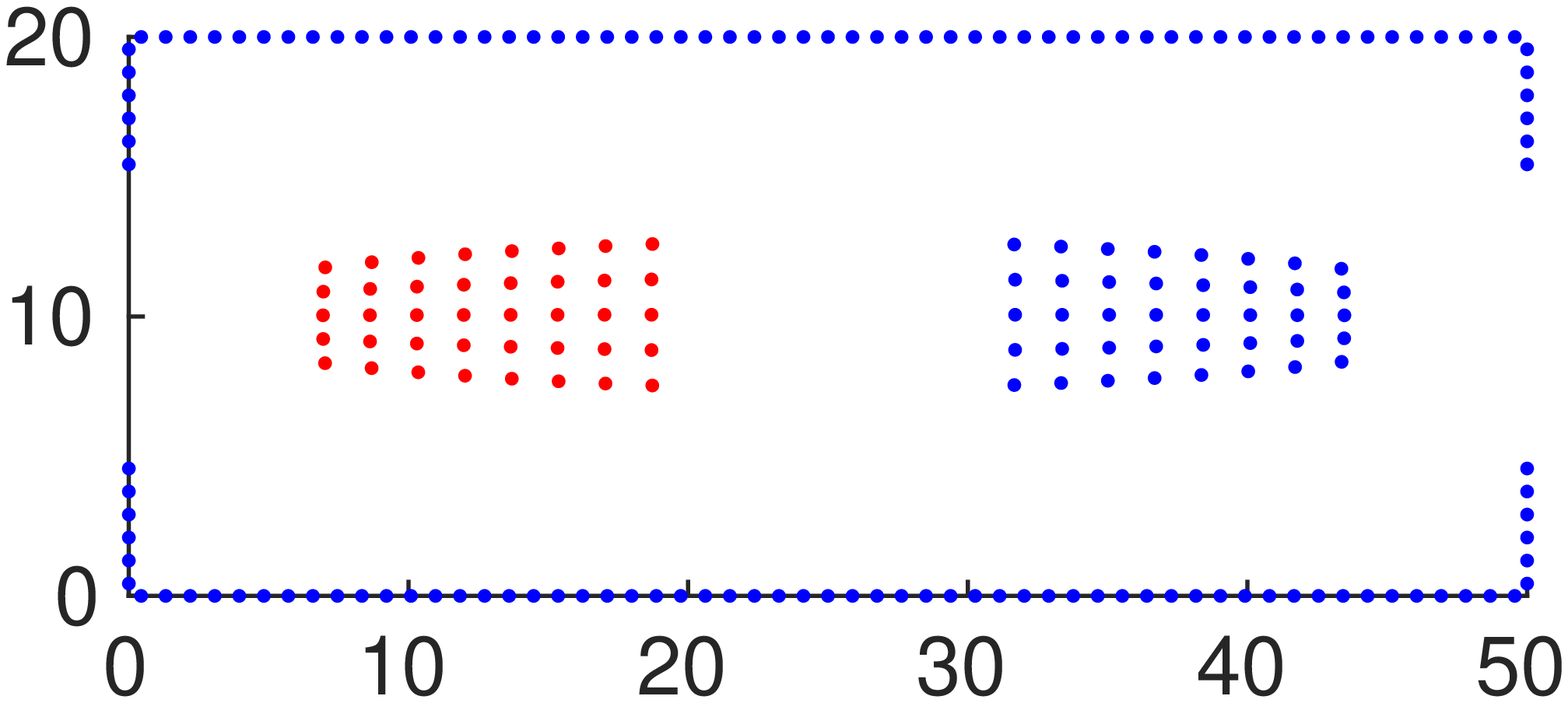}\hfill
	\includegraphics[width=.5\linewidth]{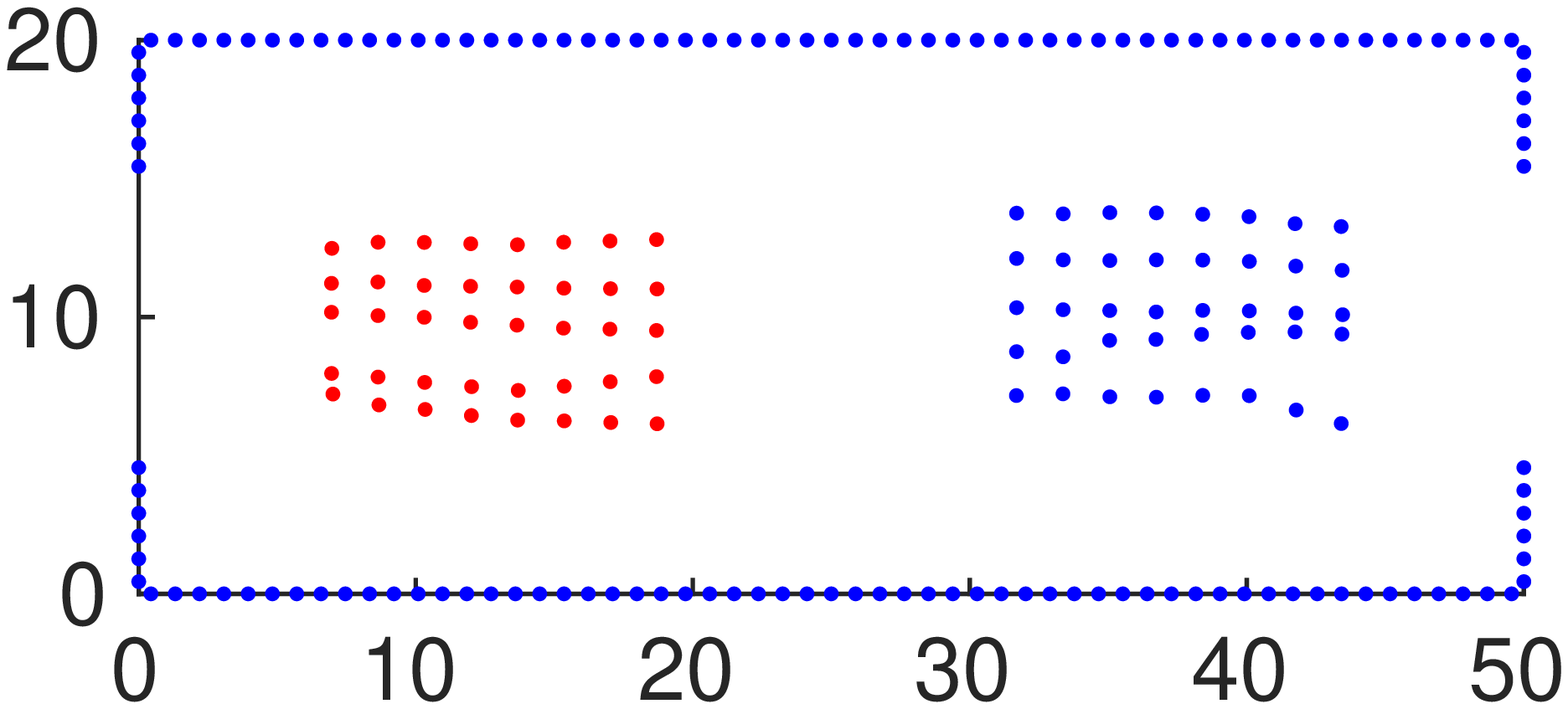}\\
	\includegraphics[width=.5\linewidth]{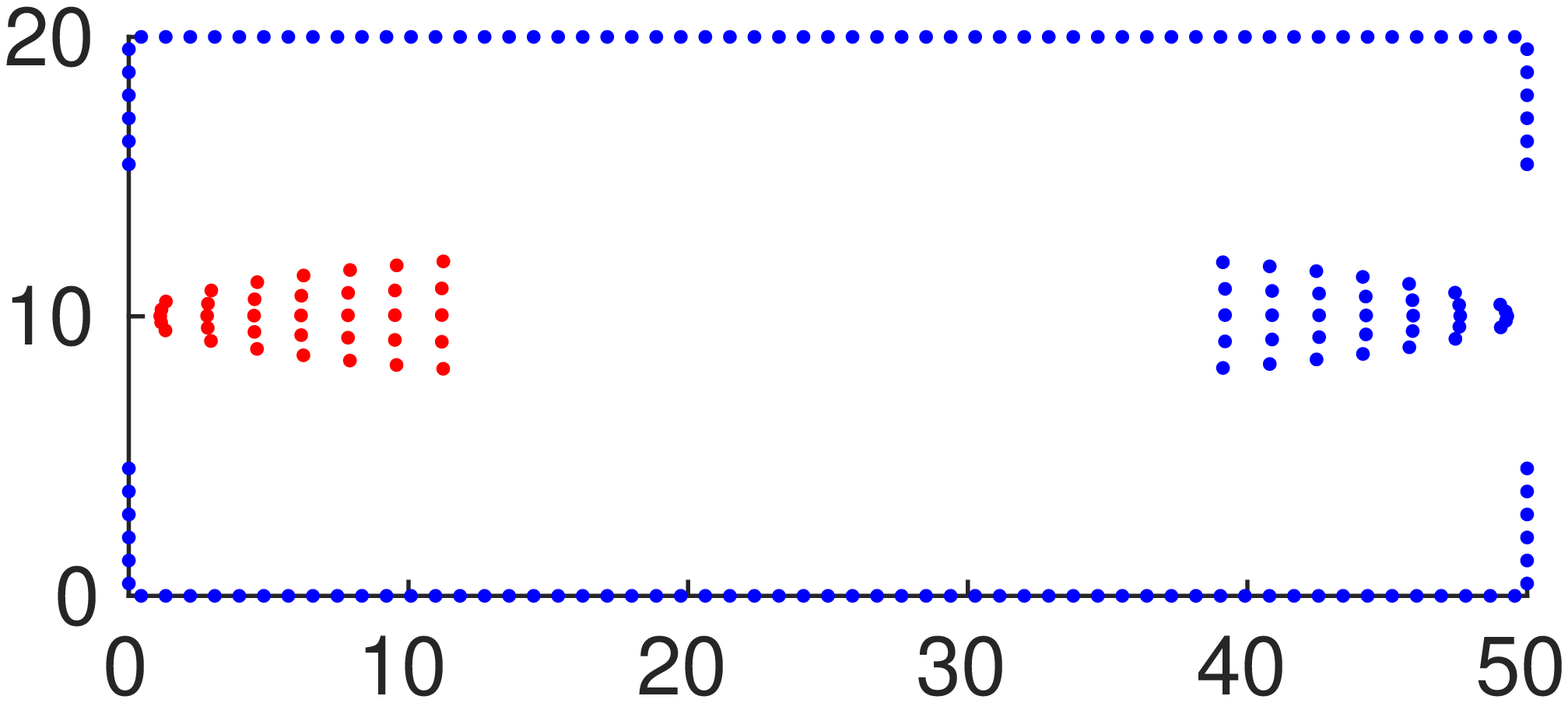}\hfill
	\includegraphics[width=.5\linewidth]{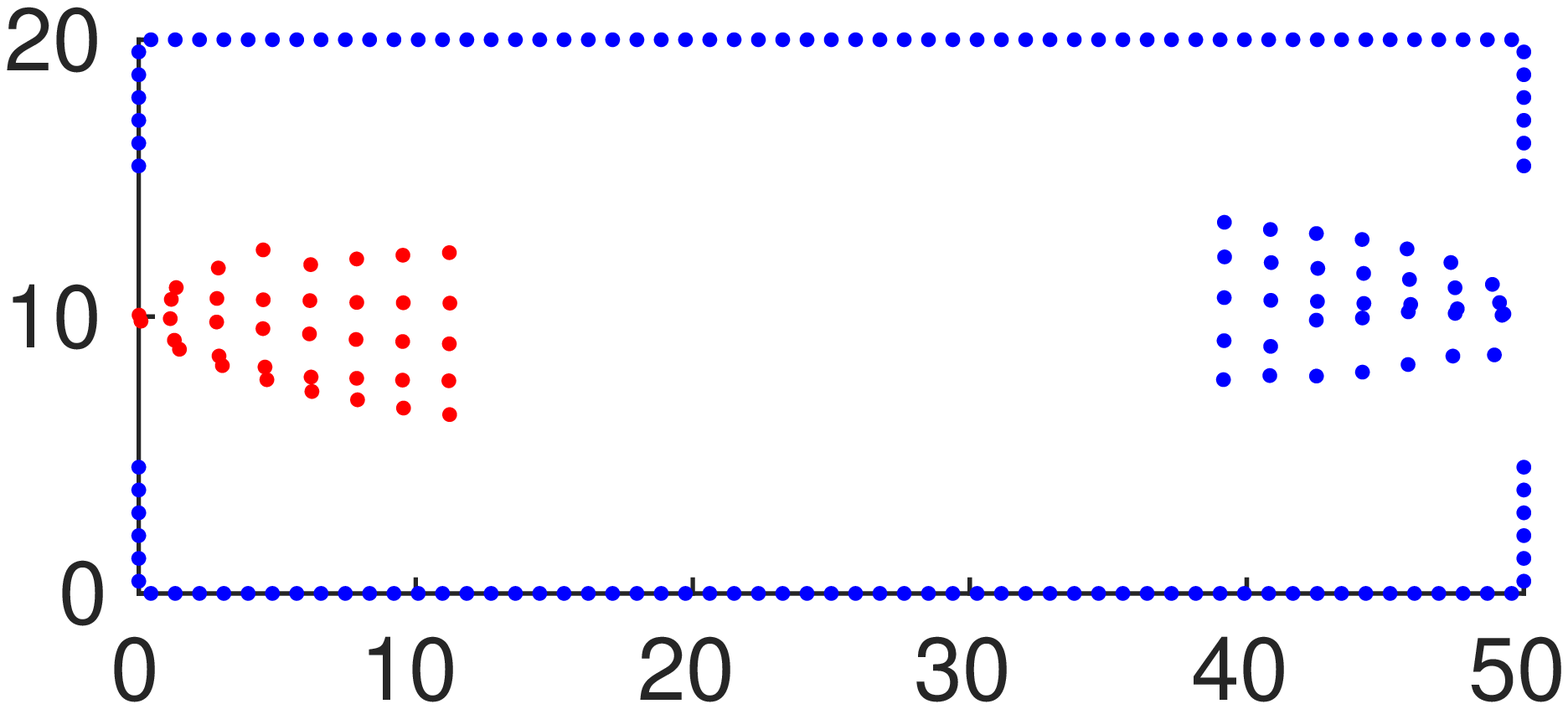}
	\caption{Distribution of grid particles in the model without control on direction (first column) and with control on direction (second column) at time $t = 8,10, 15,20,25$ (top to bottom) respectively.}
	\label{fig:2}
\end{figure}

\begin{figure}
	\includegraphics[scale=0.5]{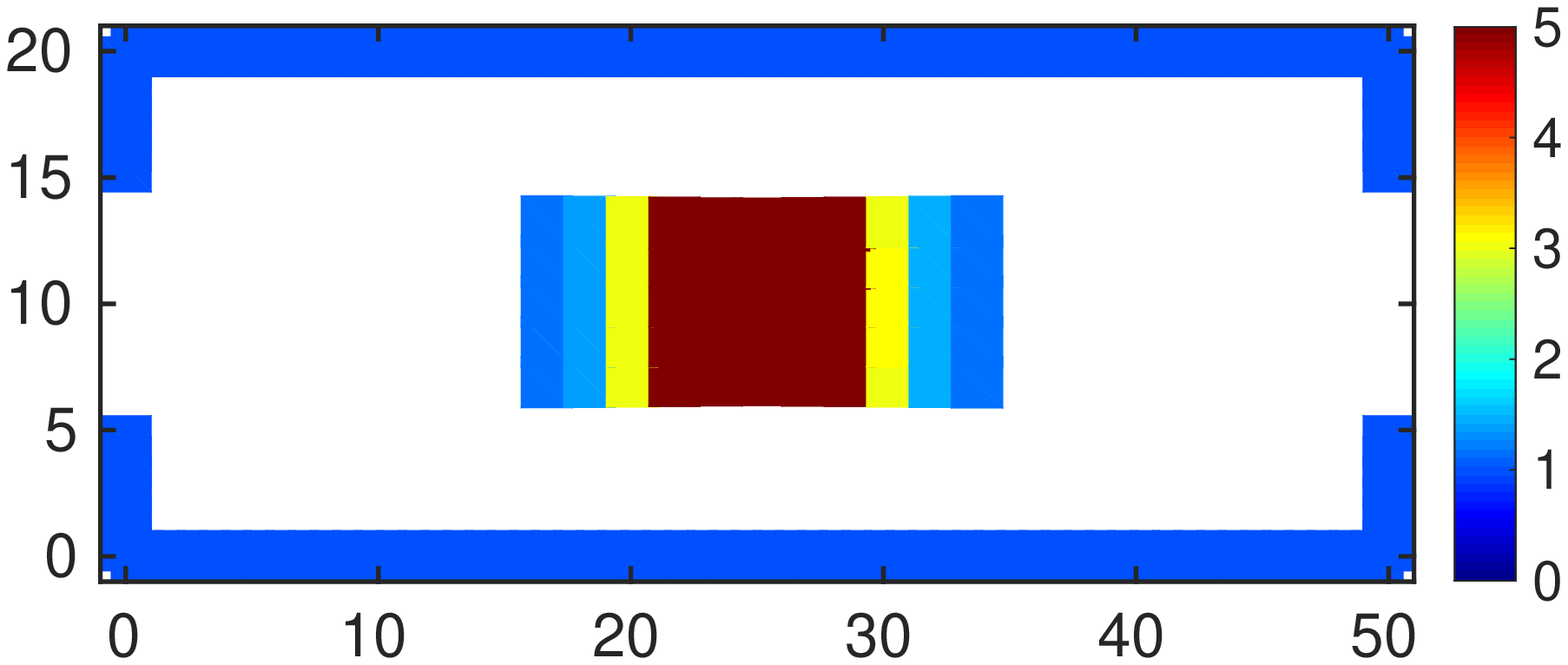}\\
	\includegraphics[scale=0.5]{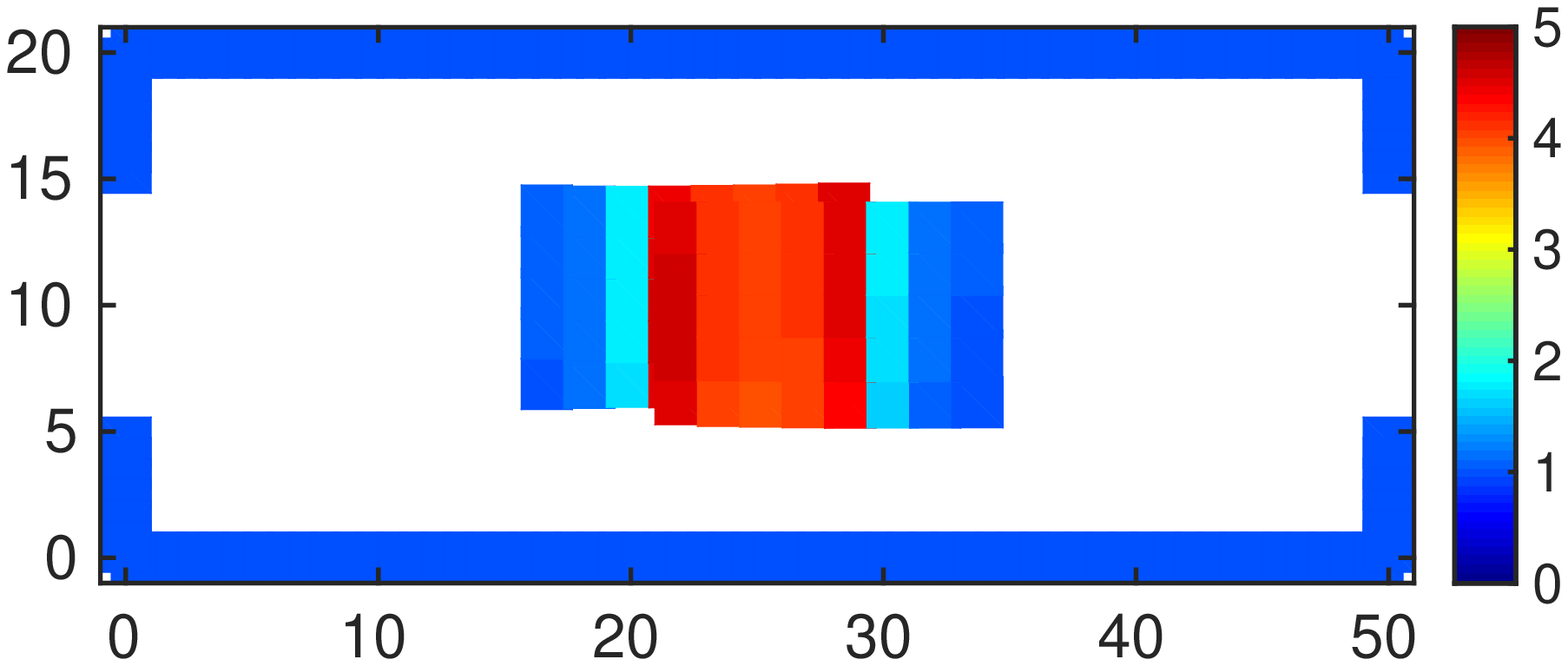}
	\caption{Density of pedestrians in the model without control on direction (top) and with control on direction (bottom) at time $t = 10$.}
	\label{fig:3} 
\end{figure}

Figure \ref{fig:2} shows the time evolution of the grid particles in the model without control on direction (first column) and with control on direction (second column) at time $t = 8$, $t = 10$, $t = 15$, $t = 20$ and $t = 25$. One observes that when grid particles are close to each other then they avoid collision by changing their direction.
The results for the location of the grid points show similarities with those obtined in \cite{ondrej}. Figure \ref{fig:3} shows the corresponding density plots for the time $t = 10$ in the model without control on direction (top) and with control on direction (bottom). The  density is higher in case  there is no control on direction.

\begin{figure}
	\includegraphics[width=.5\linewidth]{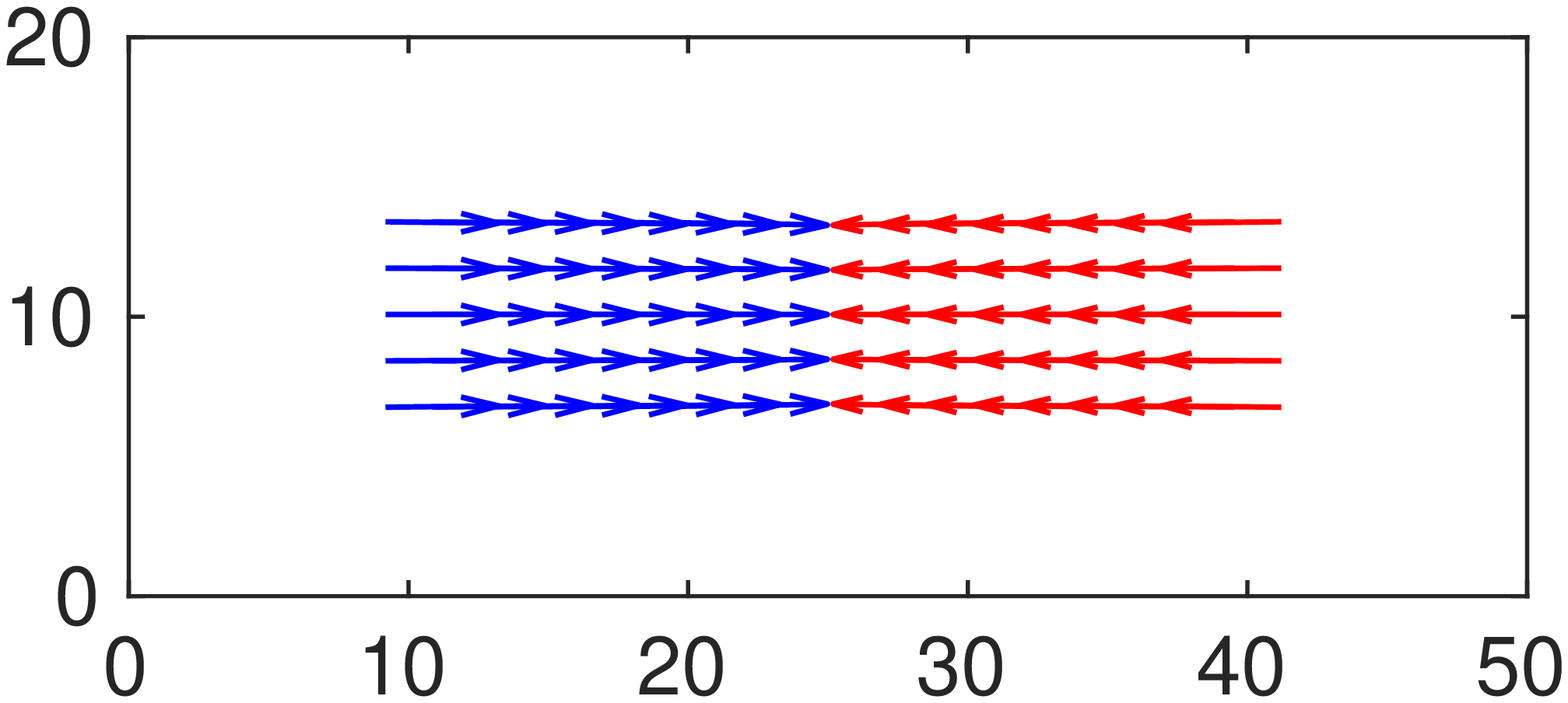}\hfill
	\includegraphics[width=.5\linewidth]{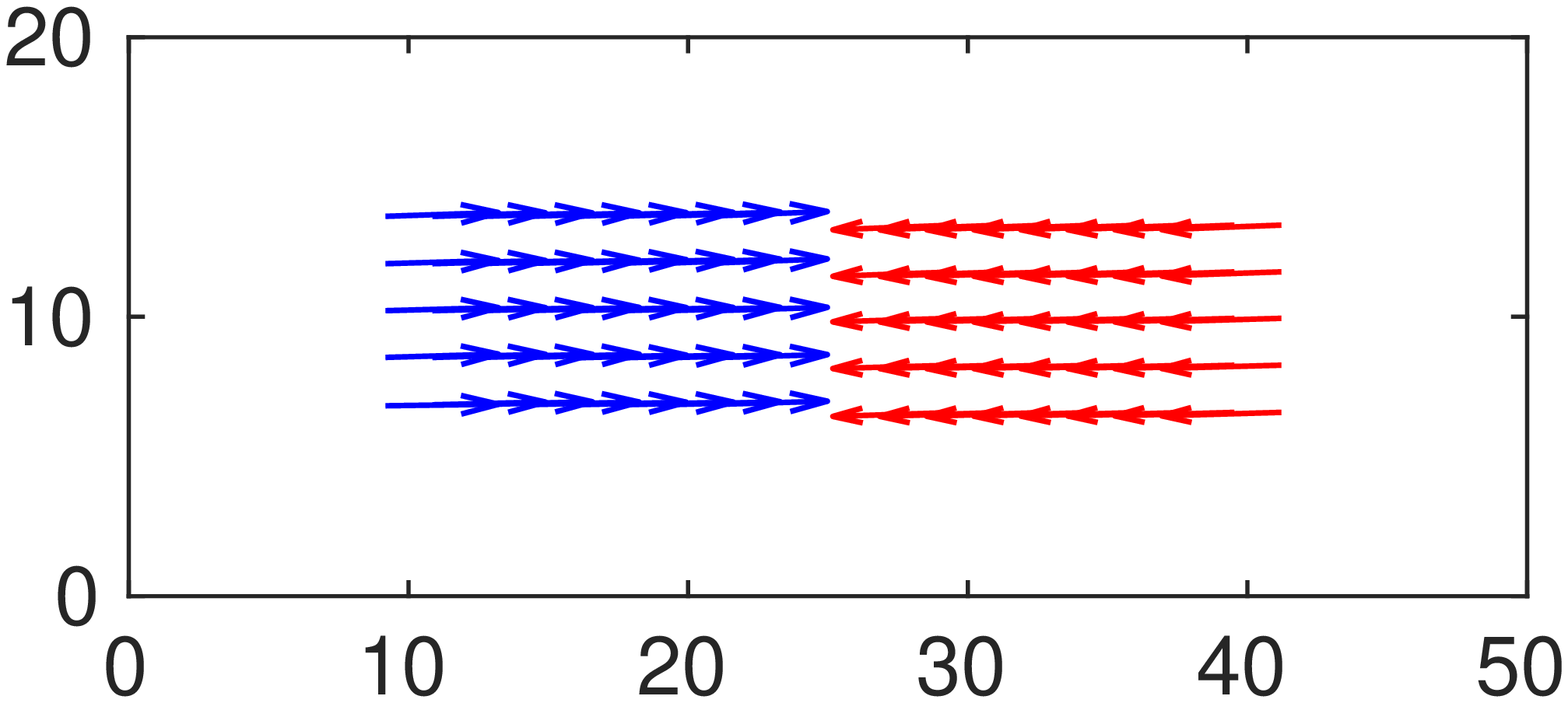}\\
	\includegraphics[width=.5\linewidth]{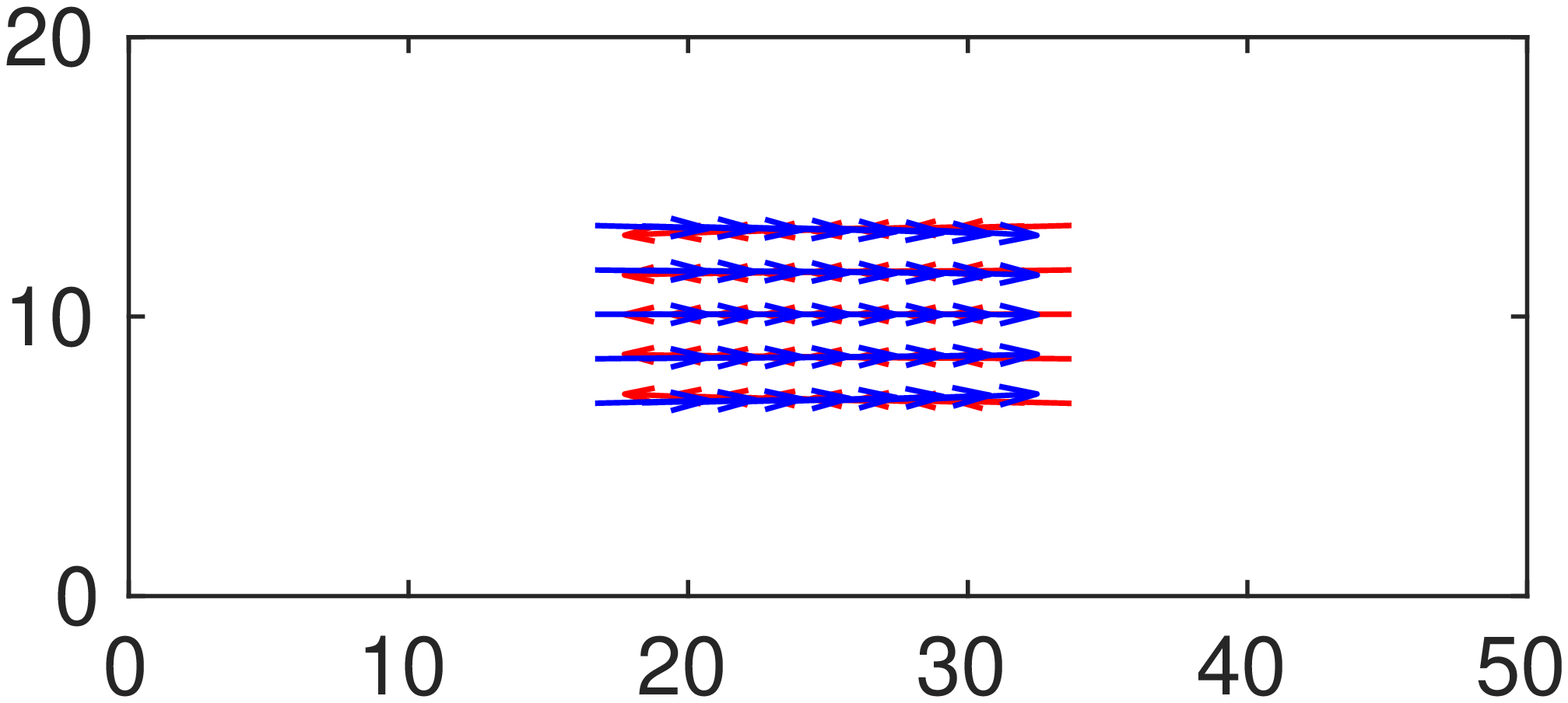}\hfill
	\includegraphics[width=.5\linewidth]{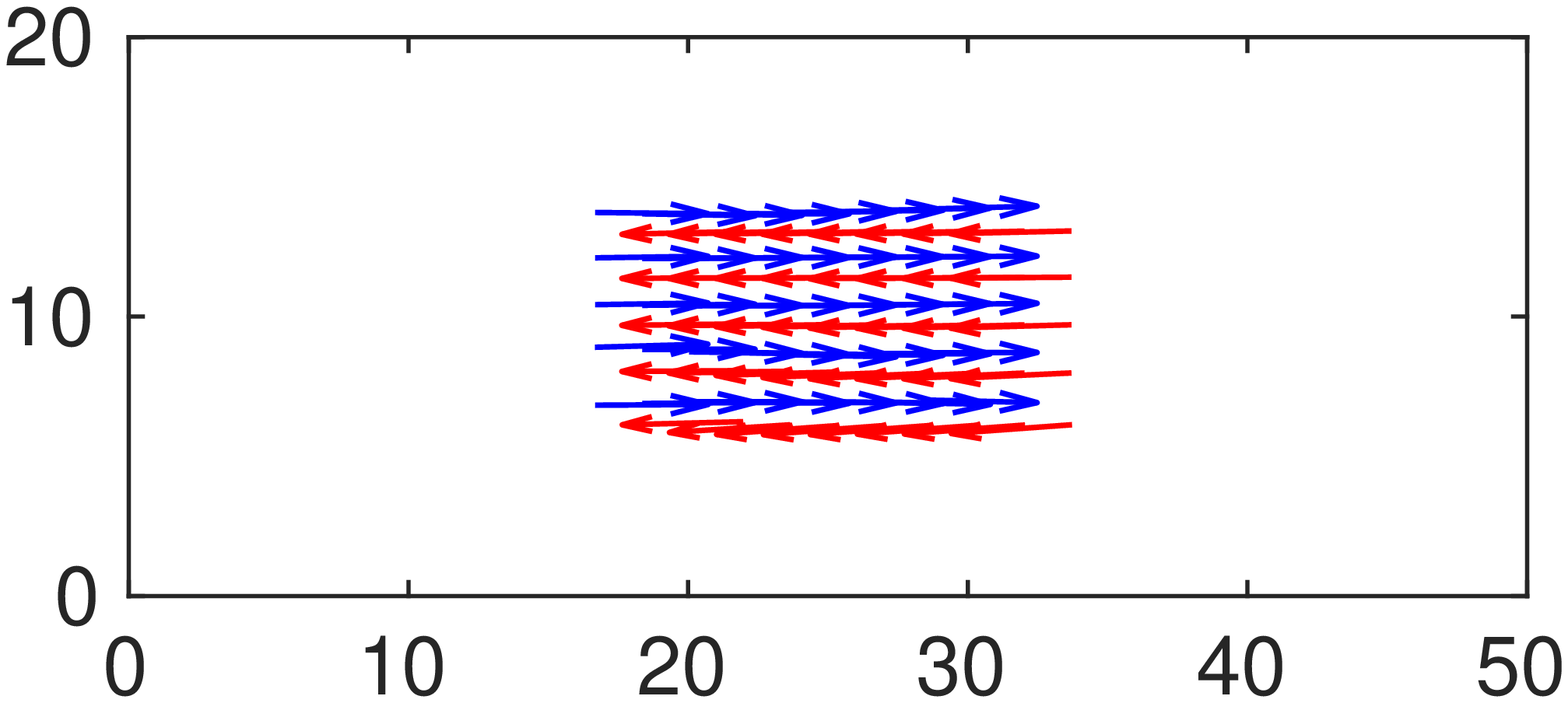}\\
	\includegraphics[width=.5\linewidth]{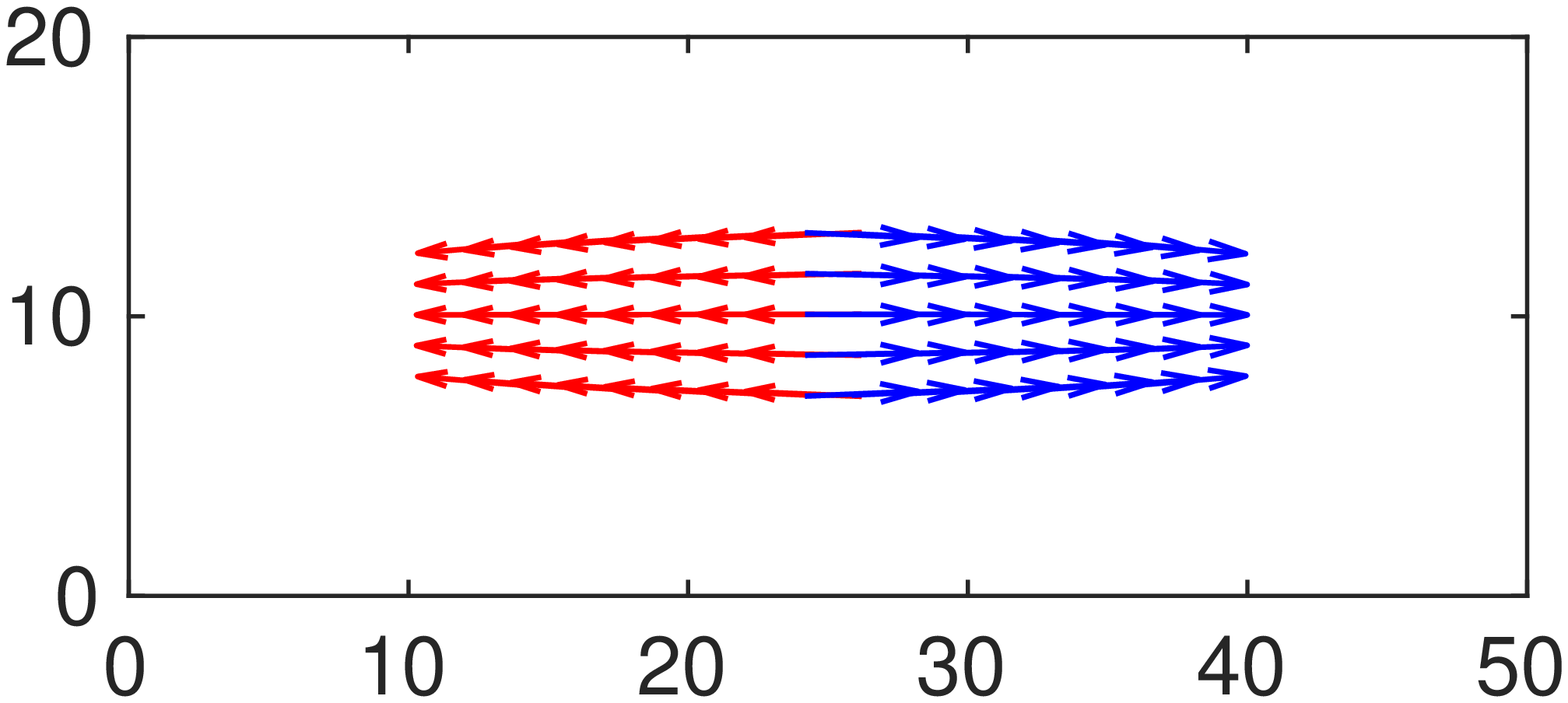}\hfill
	\includegraphics[width=.5\linewidth]{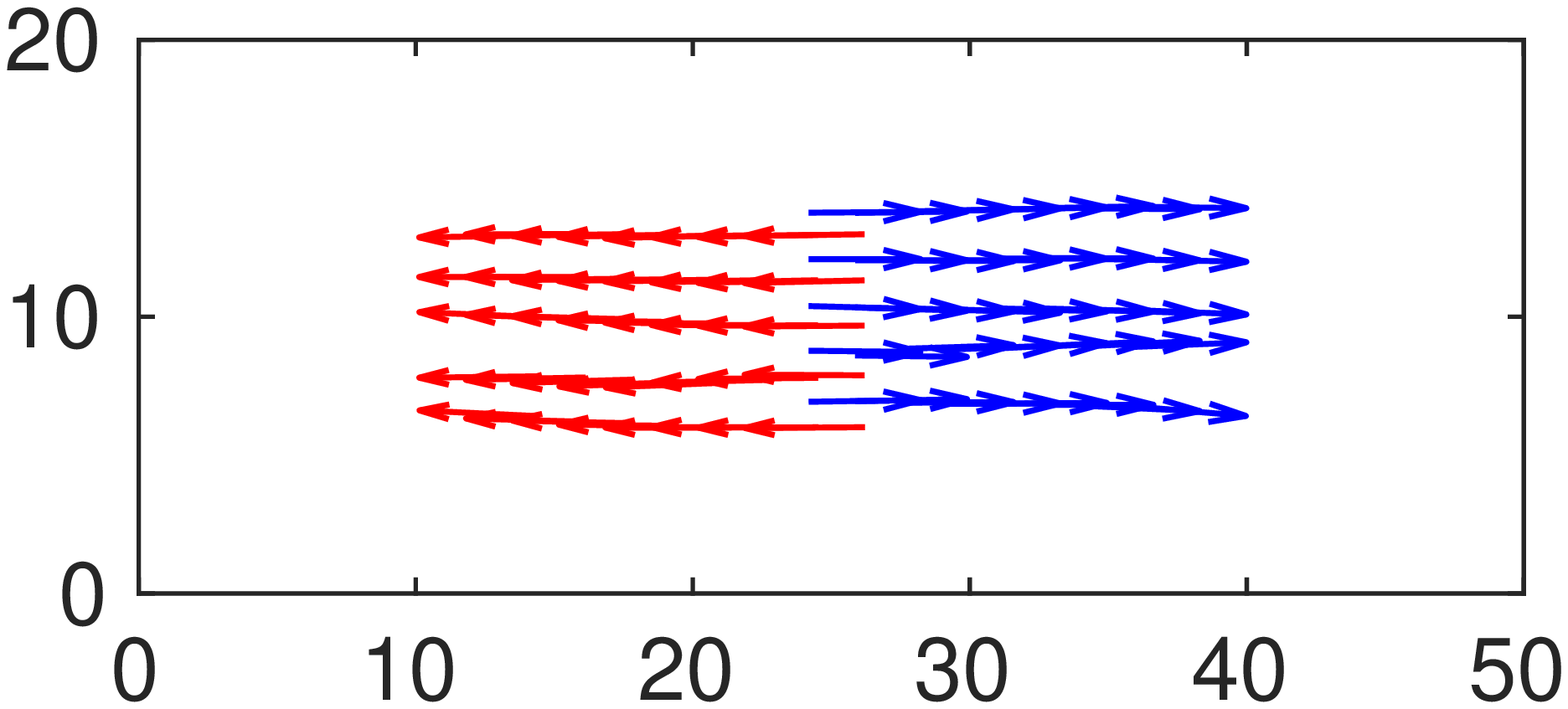}\\
	\includegraphics[width=.5\linewidth]{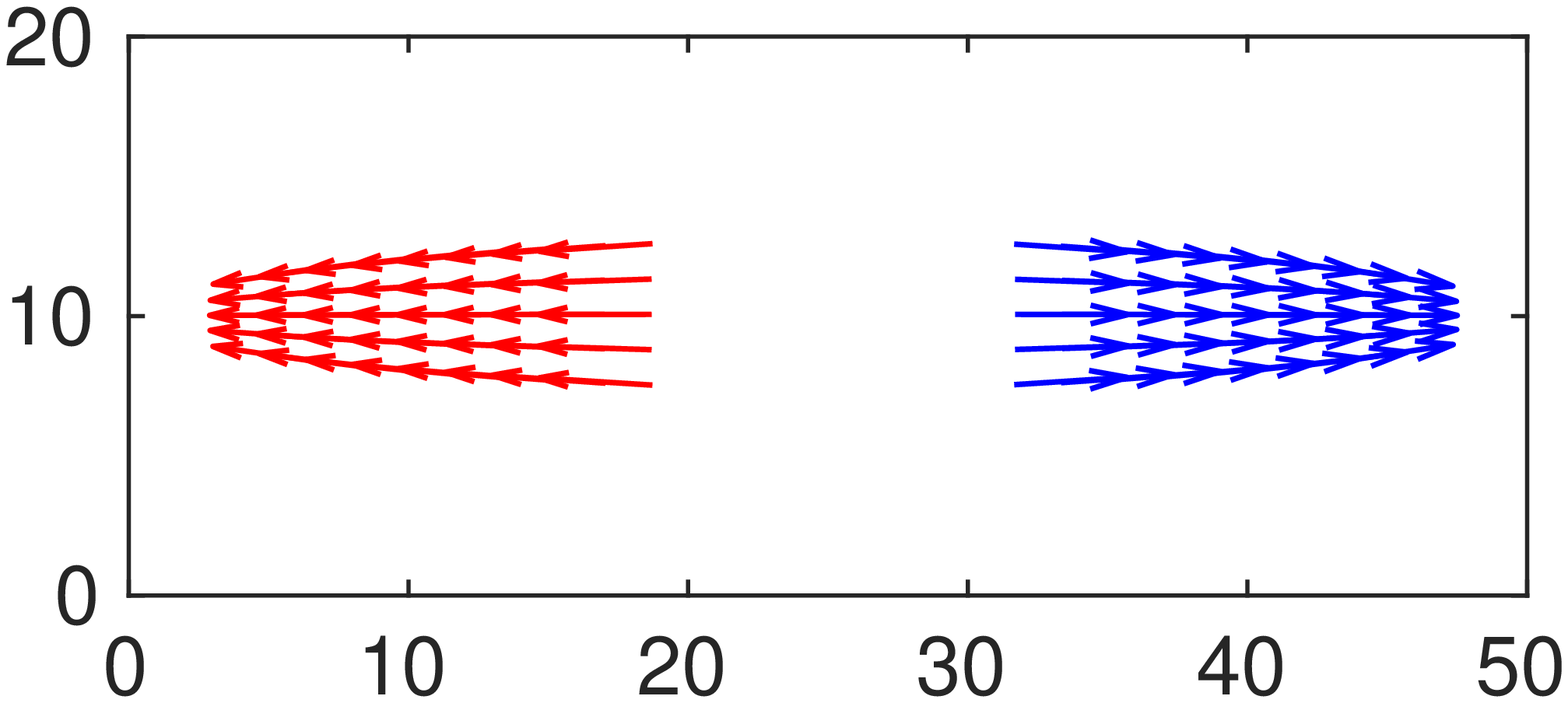}\hfill
	\includegraphics[width=.5\linewidth]{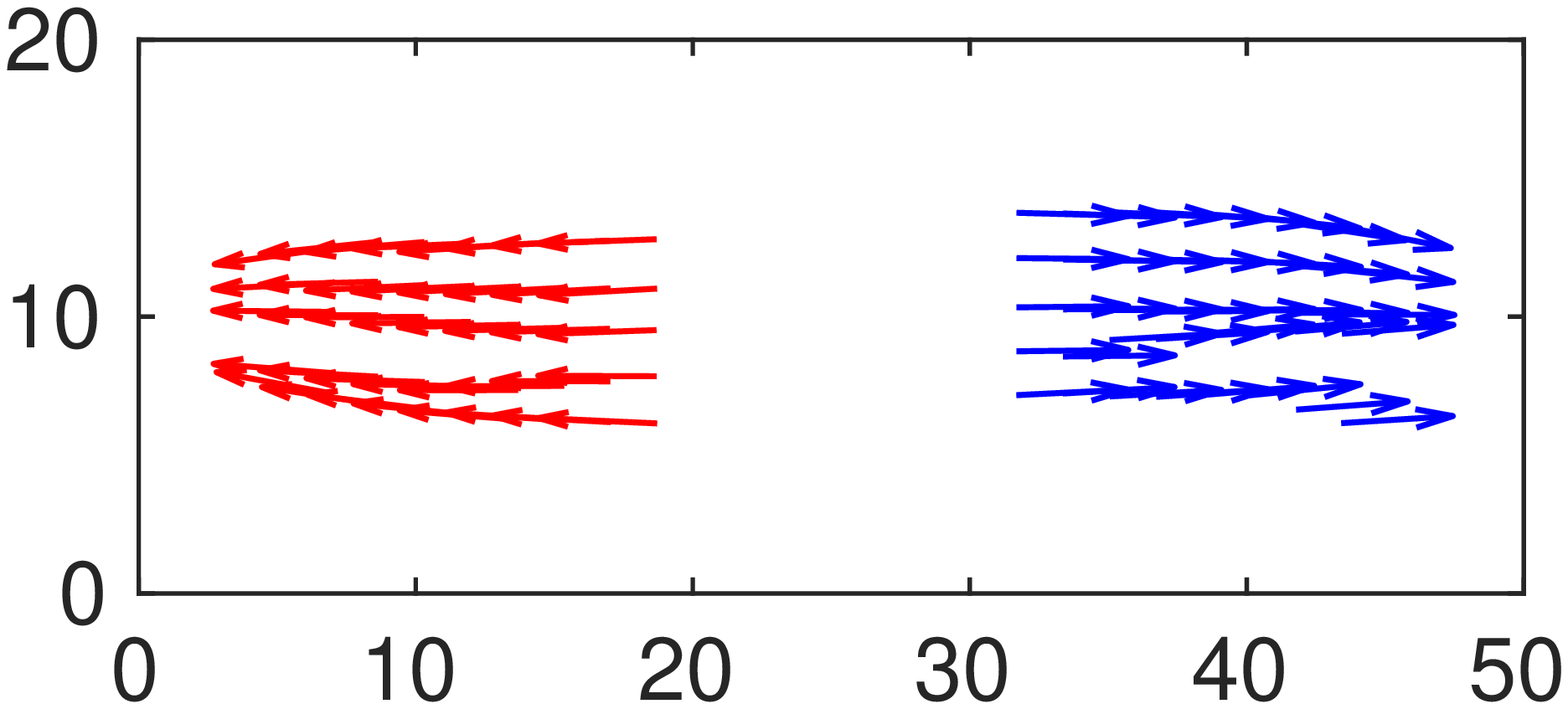}
	\caption{Direction of grid particles in the model without control on direction (first column) and with control on direction (second column) at time $t = 5,10, 15,20$ (top to bottom) respectively.}
	\label{fig:4}
\end{figure}

Figure \ref{fig:4} shows the time evolution of particle paths in the model at time $t = 5$, $t = 10$, $t = 15$ and $t = 20$. One observes the change in the direction of particle paths when pedestrians  are close to each other.

\begin{figure}
	\includegraphics[width=.5\linewidth]{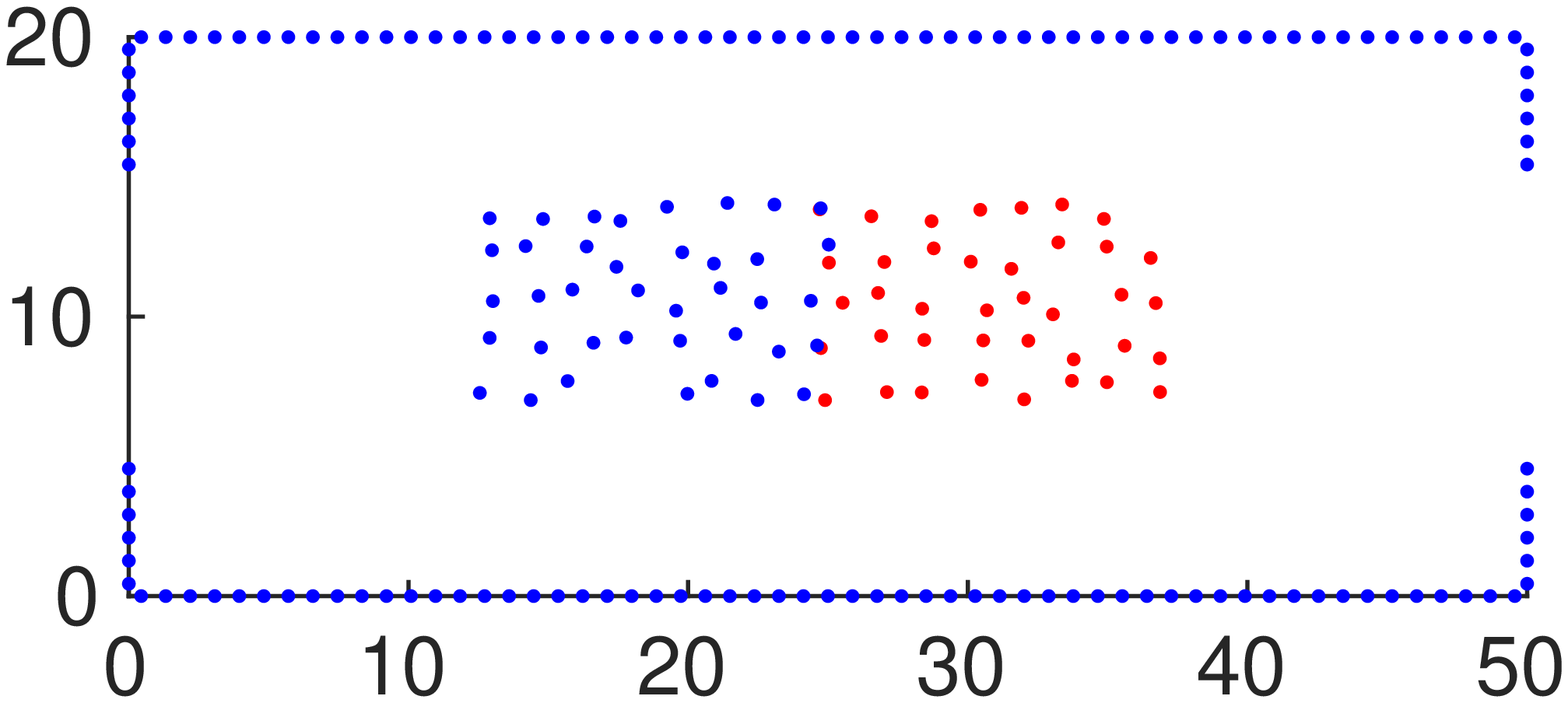}\hfill
	\includegraphics[width=.5\linewidth]{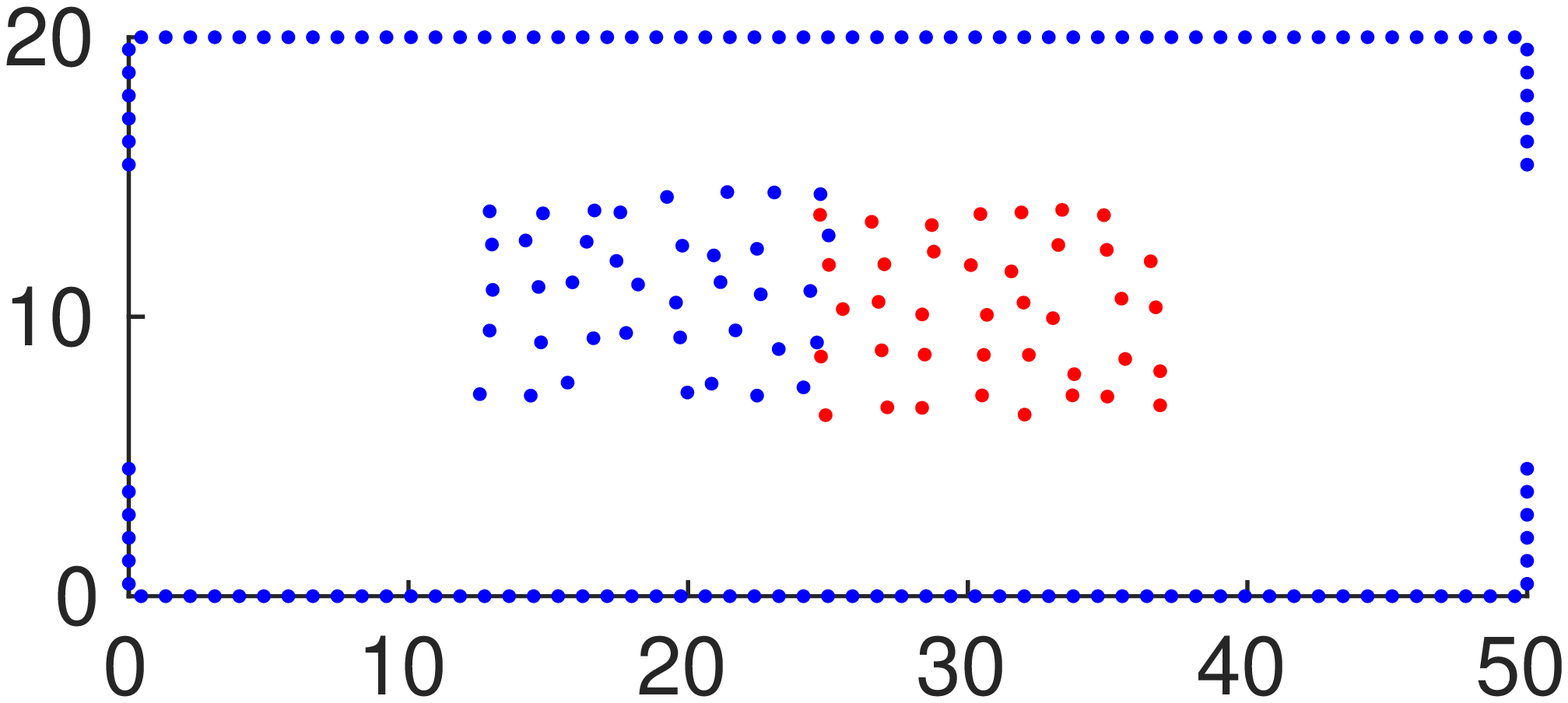}\\
	\includegraphics[width=.5\linewidth]{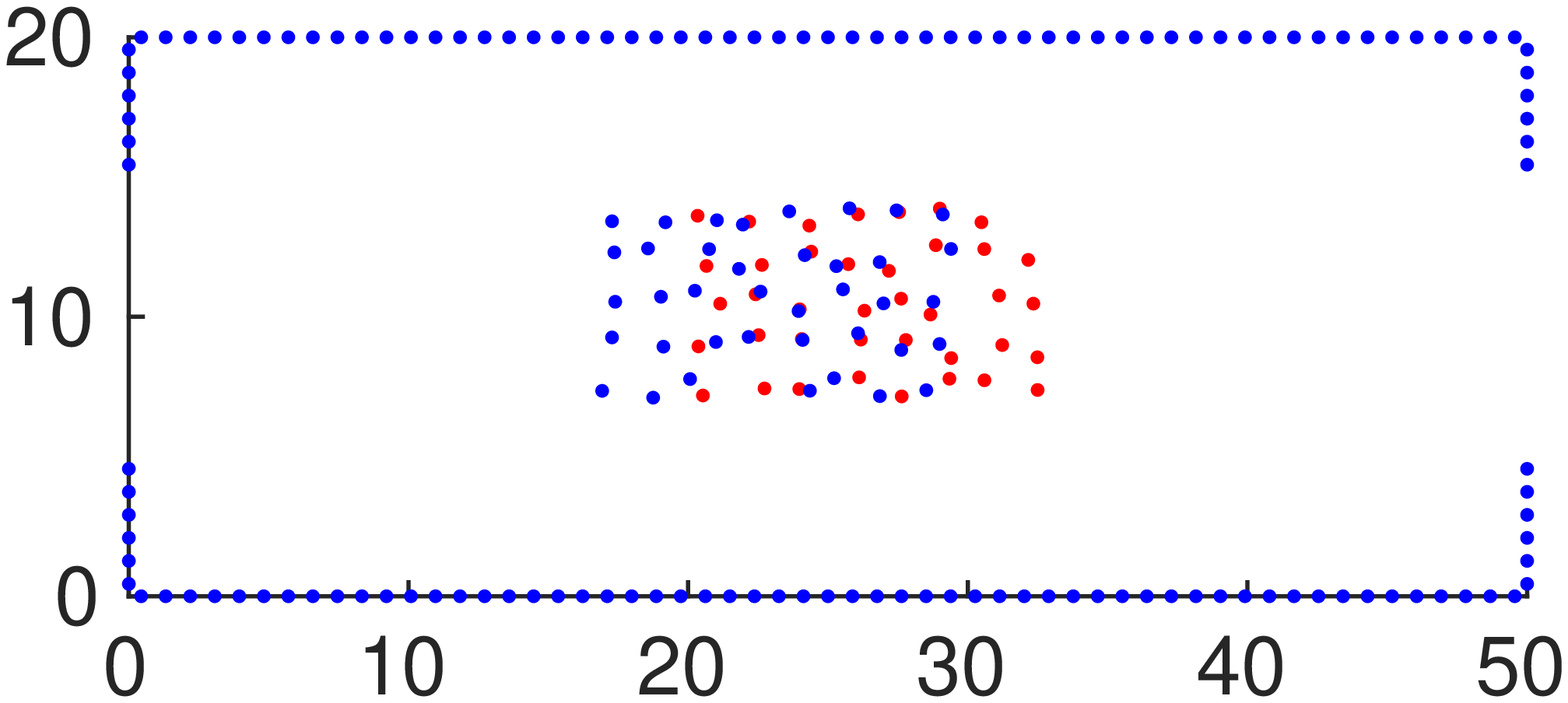}\hfill
	\includegraphics[width=.5\linewidth]{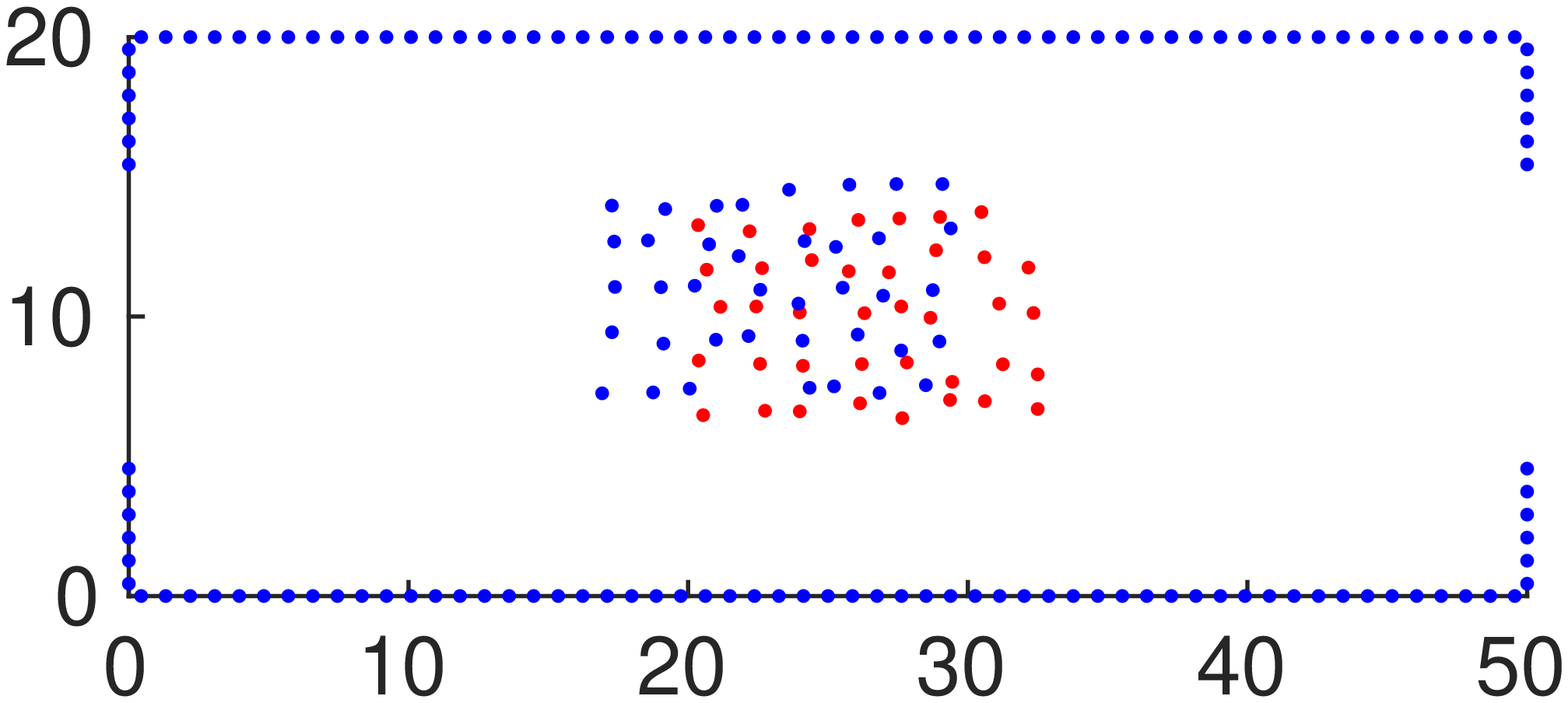}\\
	\includegraphics[width=.5\linewidth]{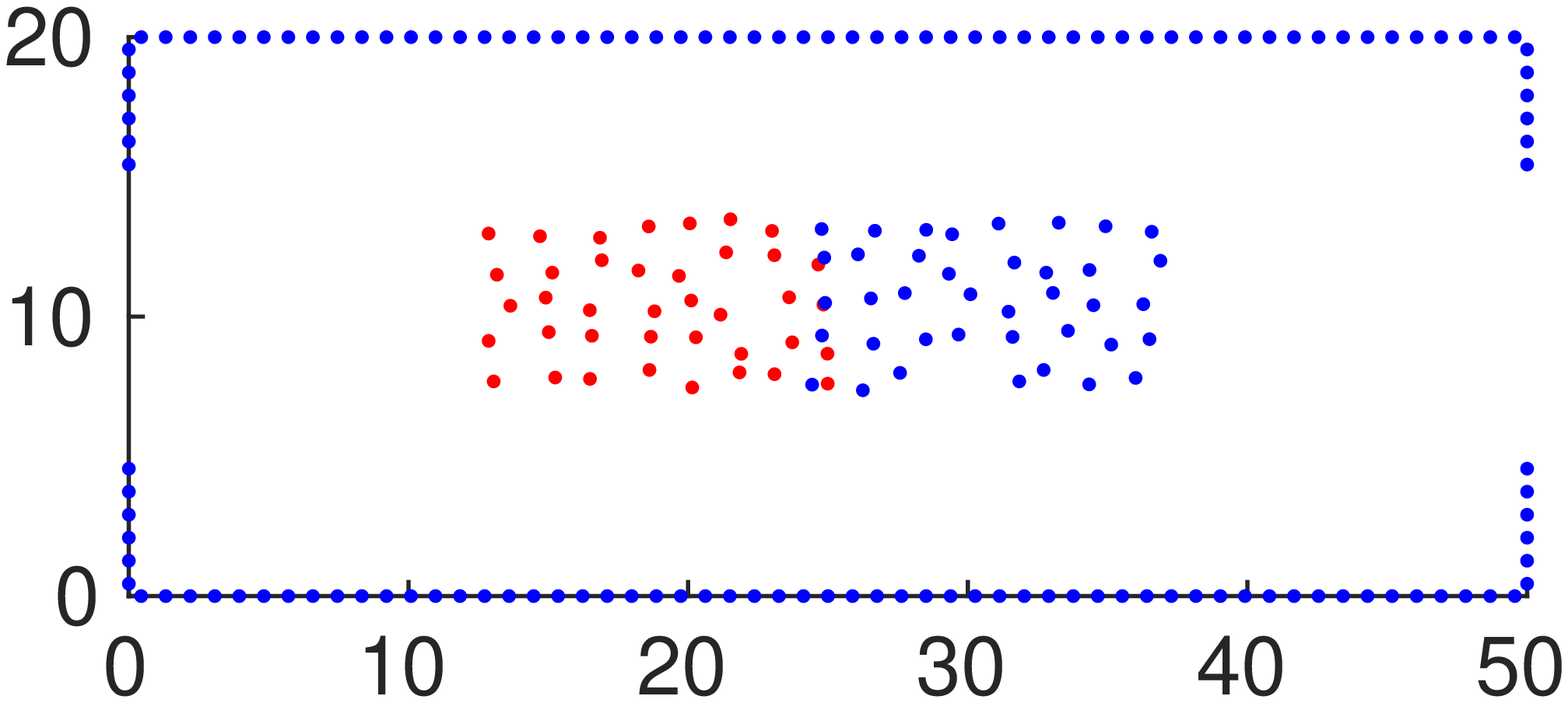}\hfill
	\includegraphics[width=.5\linewidth]{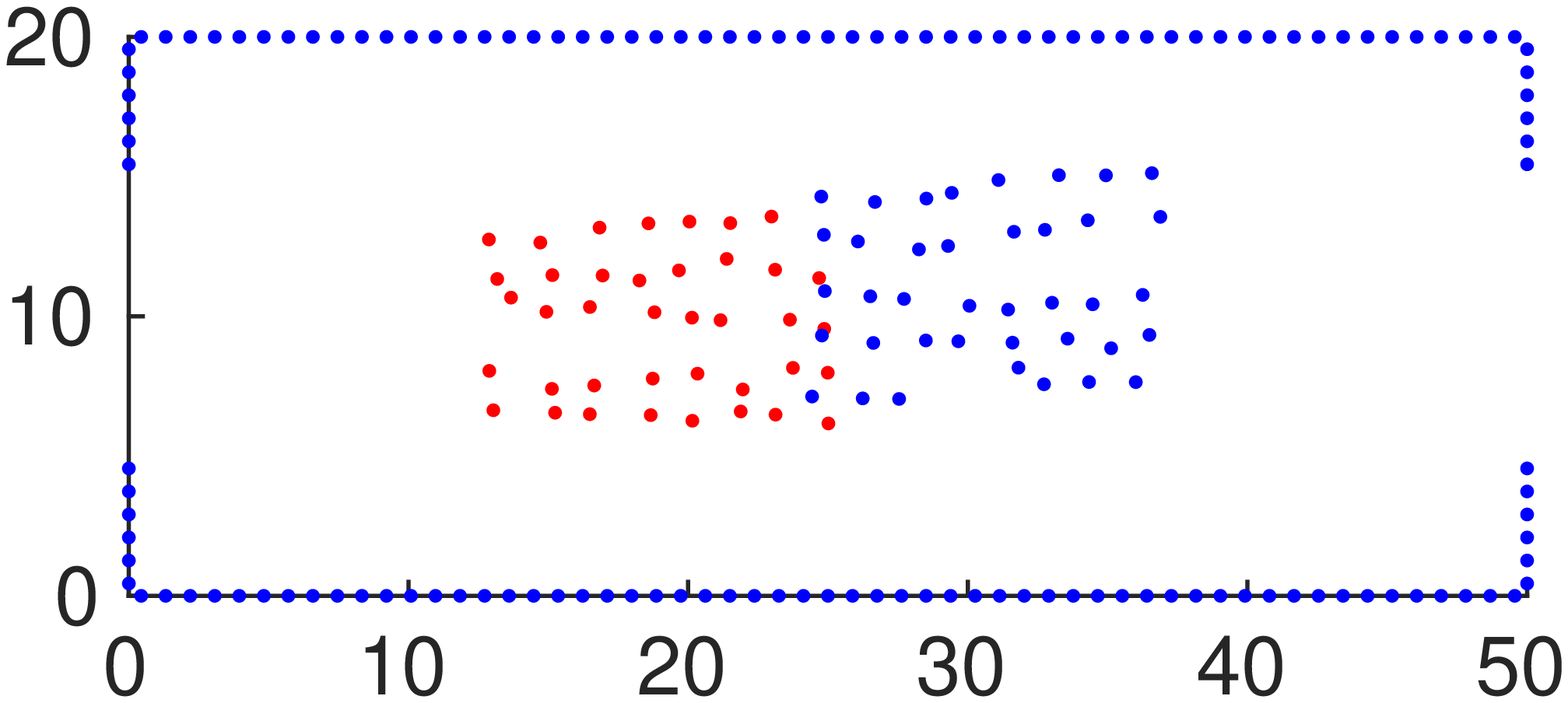}\\
	\includegraphics[width=.5\linewidth]{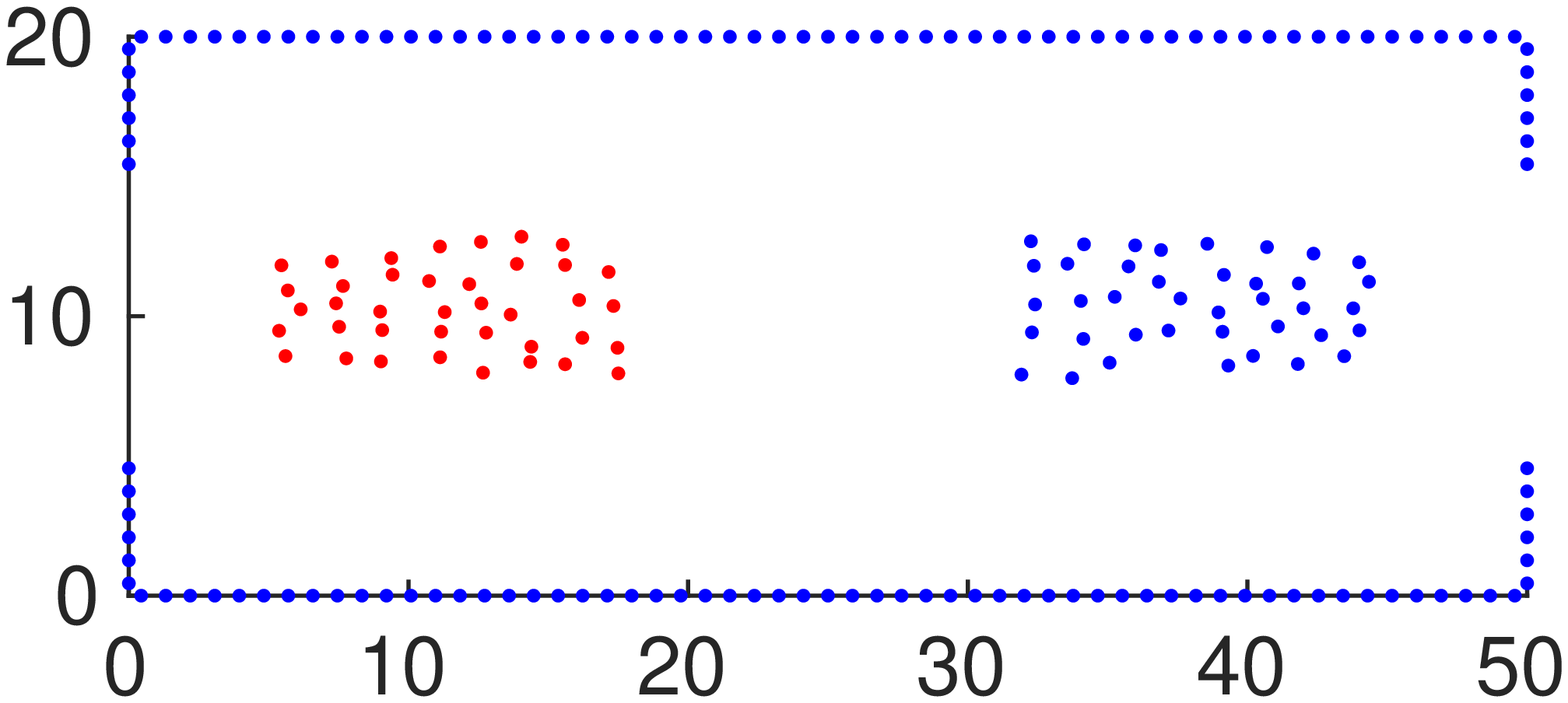}\hfill
	\includegraphics[width=.5\linewidth]{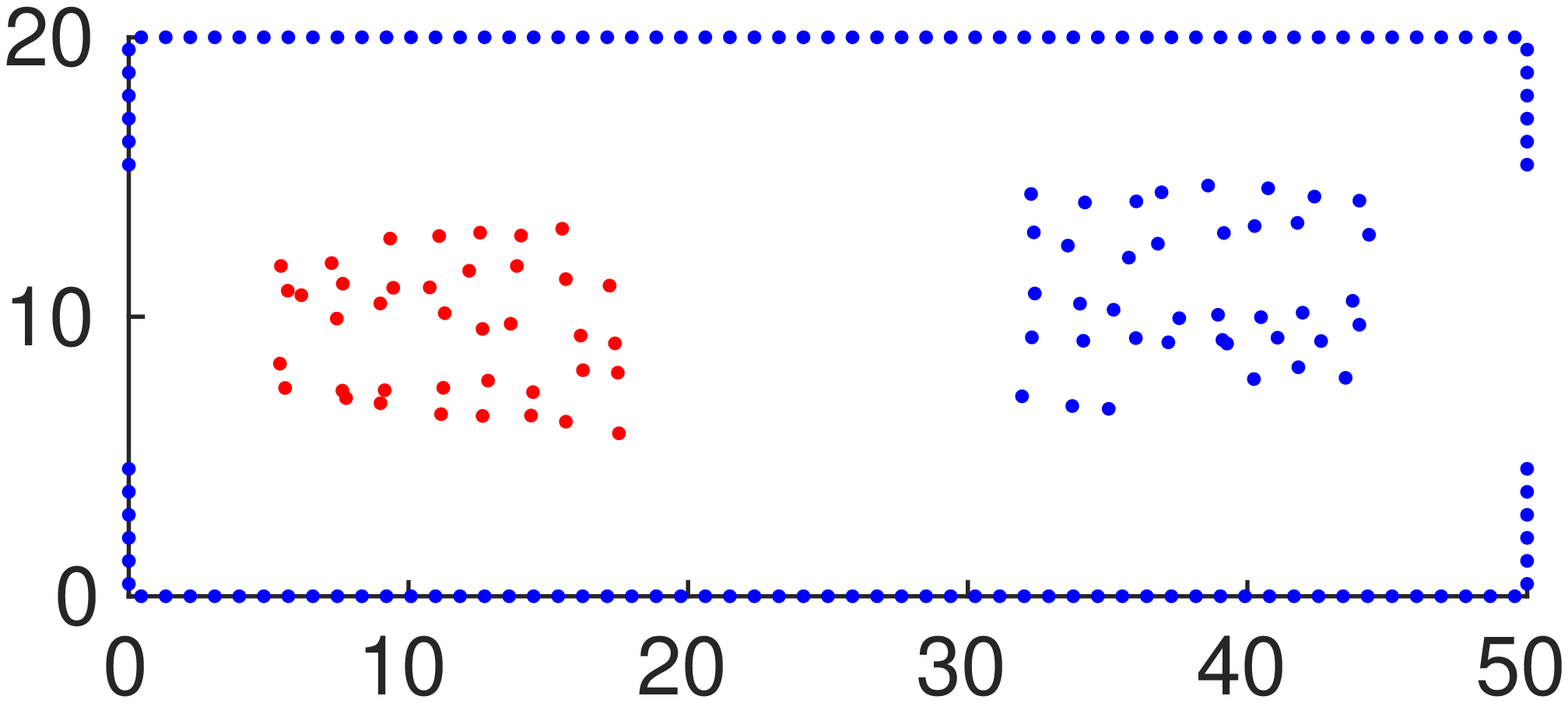}\\
	\includegraphics[width=.5\linewidth]{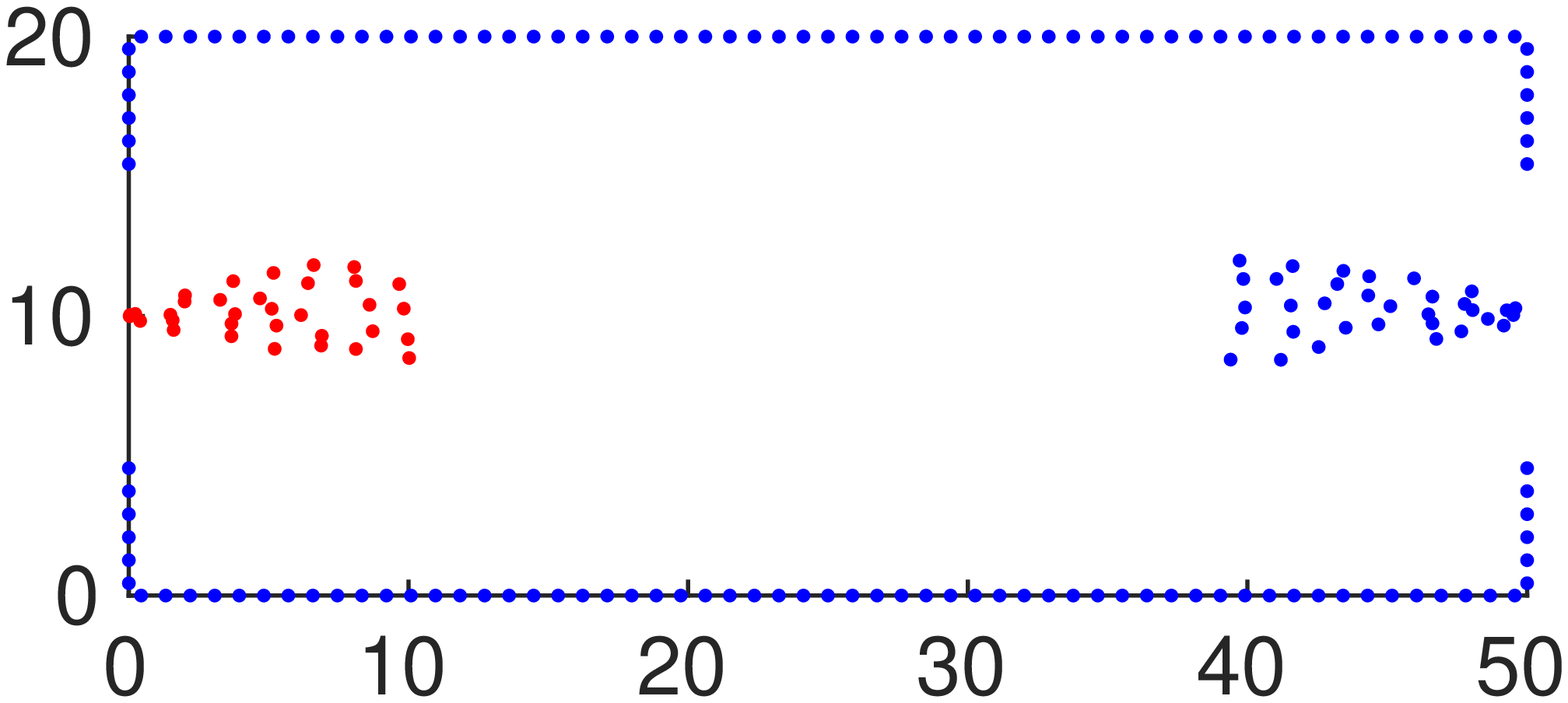}\hfill
	\includegraphics[width=.5\linewidth]{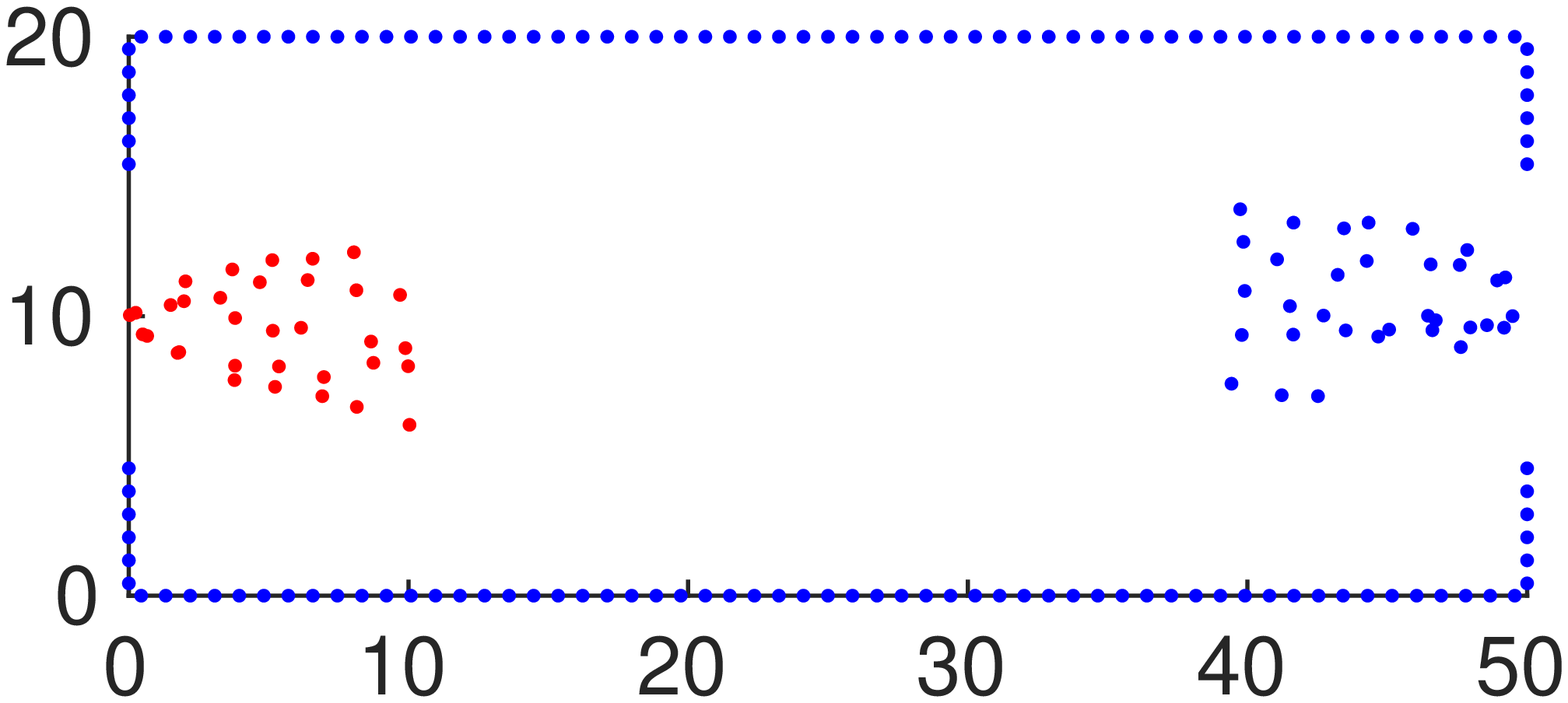}
	\caption{Random distribution of grid particles  in the model without control on direction (first column) and with control on direction (second column) at time $t =  7.14, 10, 15, 20, 25$ (top to bottom) respectively.}
	\label{fig:5}
\end{figure}

In the earlier examples, we have initialized particles in a regular lattice. In real situations, the initial positions can be randomly distributed. Therefore, in figure \ref{fig:5} we consider the time evolution of randomly distributed grid particles in the model without control on direction (first column) and with control on direction (second column) at time $t = 7.14$, $t = 10$, $t = 15$, $t=20$ and $t = 25$. Again one  observes the change in  direction  to avoid collision.

\begin{figure}[htbp]
	\centering
	\includegraphics[width=0.7\linewidth]{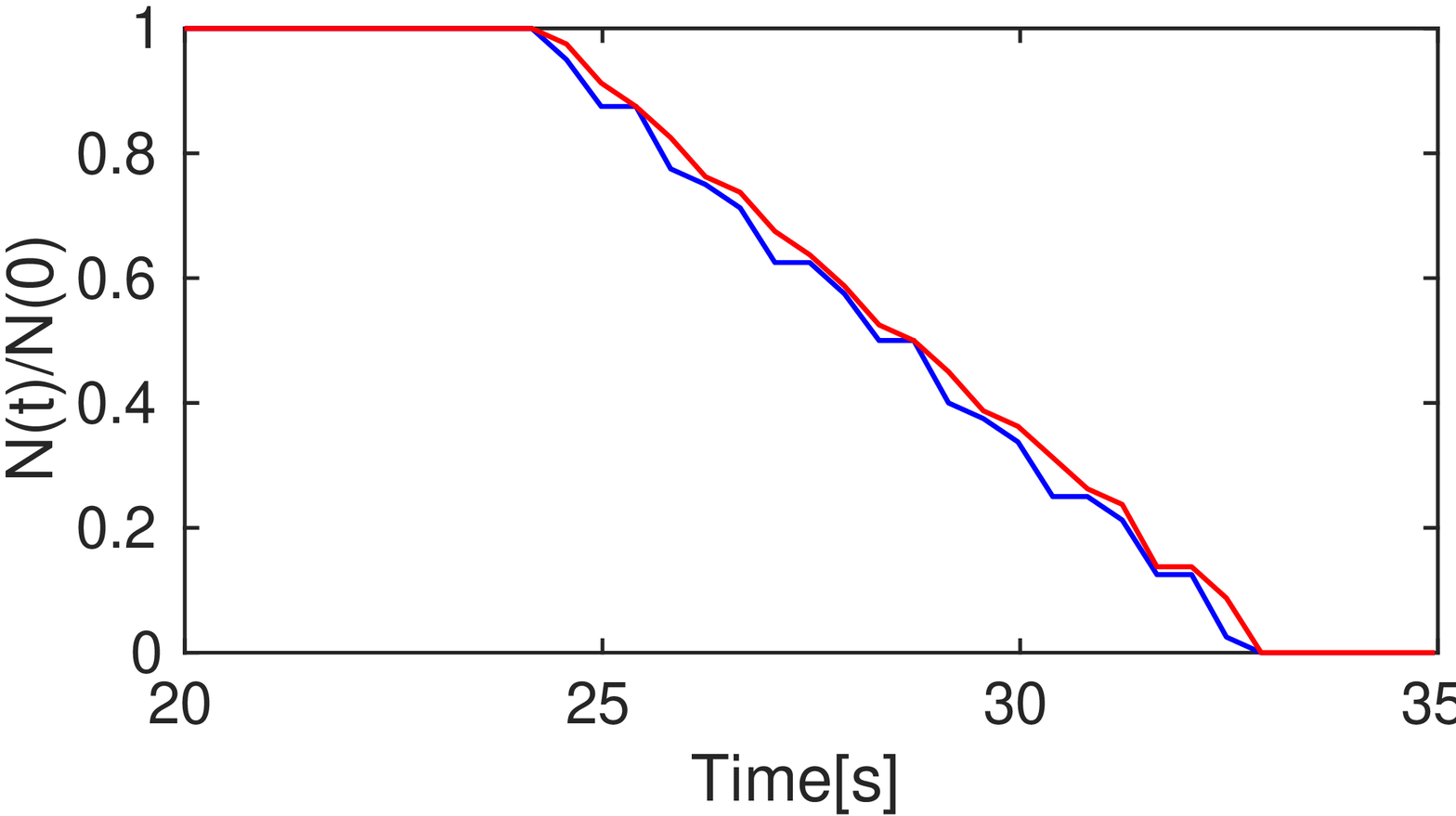}
	\caption{Ratio of initial and actual grid particles with respect to time in the model without control on direction (blue) and with control on direction (red).}
	\label{fig:6}
\end{figure}

Figure \ref{fig:6} shows the percentage of grid particles being in the computational domain for the model without control on direction and with control on direction with respect to time. One observes that the evacuation time is almost similar in this simple case. However, it may  differ for complicated and larger geometries. 

\subsection{Improvement of collision avoidance}\label{sec:5.2}
If we decrease the initial spacing of the particles, then the particles might collide after interaction with other particles. So, in this subsection we improve collision avoidance between particles by adding some small extra repulsive force as in the social force model between particle in the vision-based model, for details of extra repulsive force we refer \cite{etikyala}.

Here, we consider 100 grid points (50 gridpoints  on each side having target at the opposite end). We use the constant time step $\Delta t = 0.00021$, initial distance of grid points $\Delta x = 0.84$ and $R = 1.0$.  For repulsive force between particles, we use the following values of parameters as $k_n = 1.0$, $\gamma_n = 0.01$ and $\gamma_t = 0.01$. 

\begin{figure}
	\includegraphics[width=.5\linewidth]{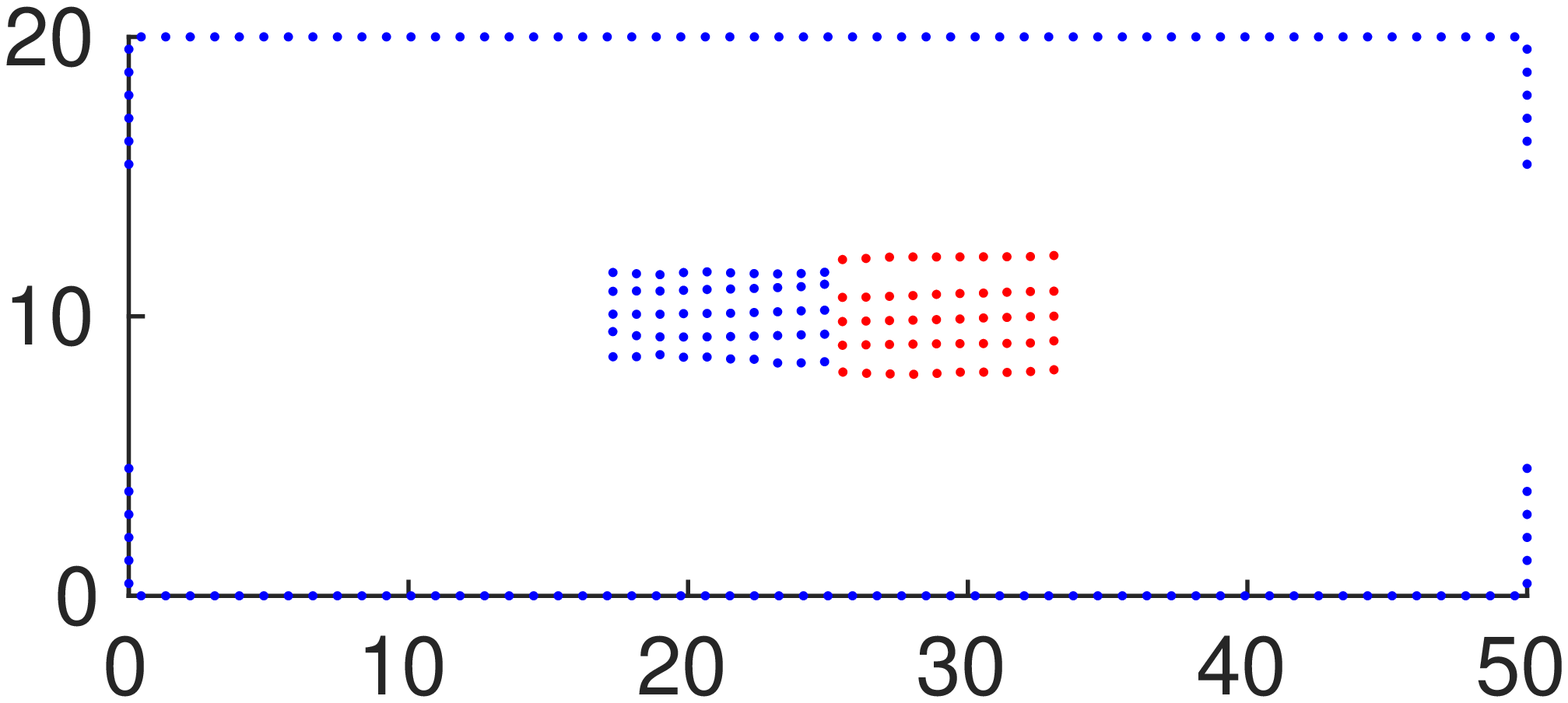}\hfill
	\includegraphics[width=.5\linewidth]{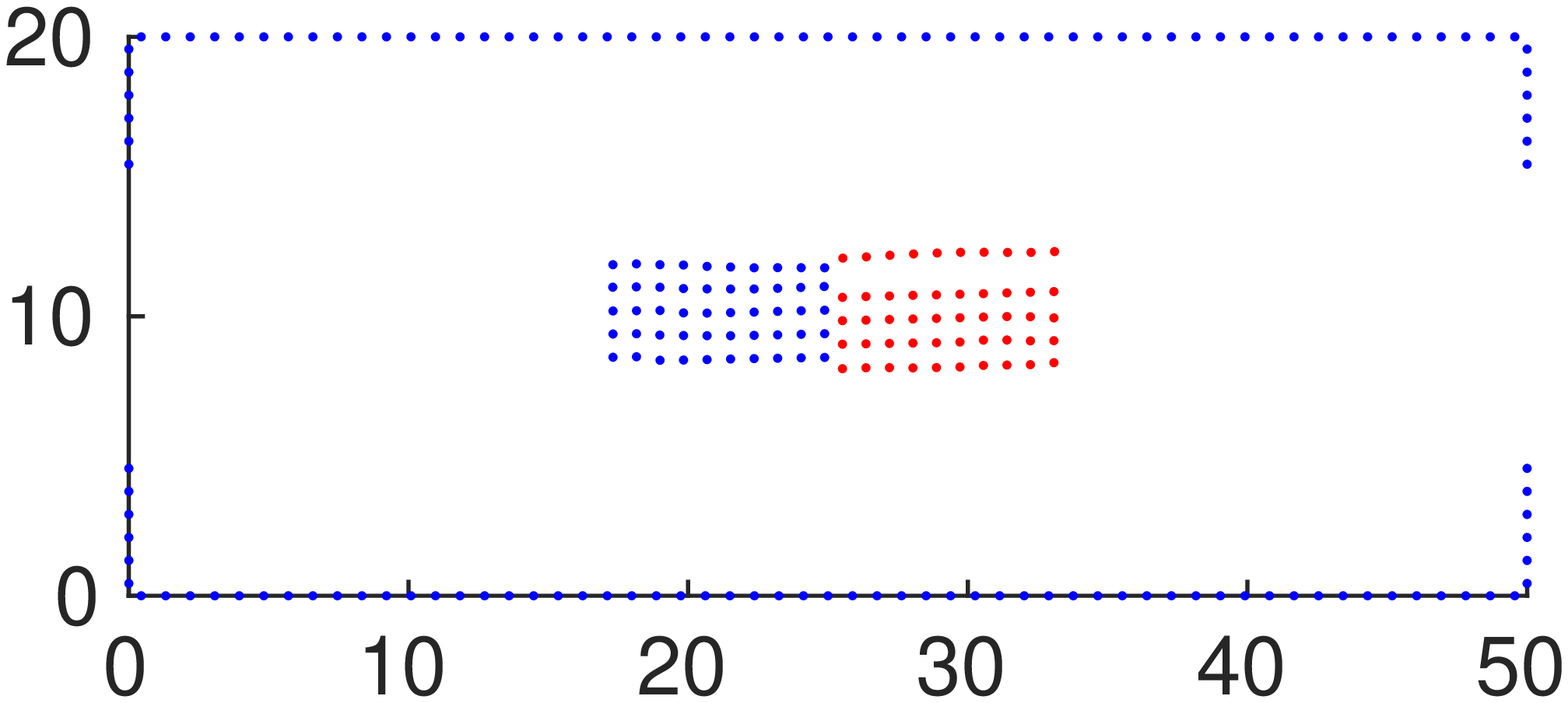}\\
	\includegraphics[width=.5\linewidth]{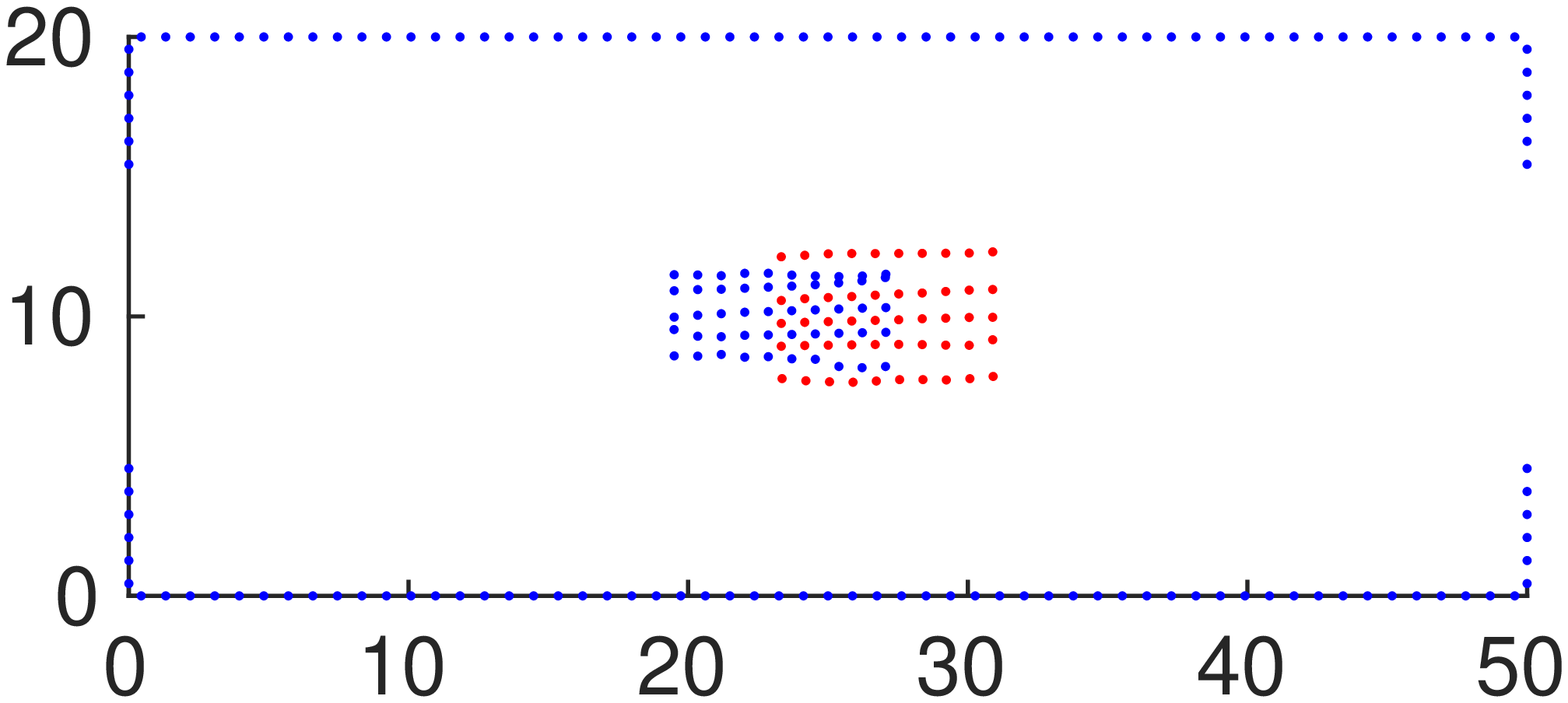}\hfill
	\includegraphics[width=.5\linewidth]{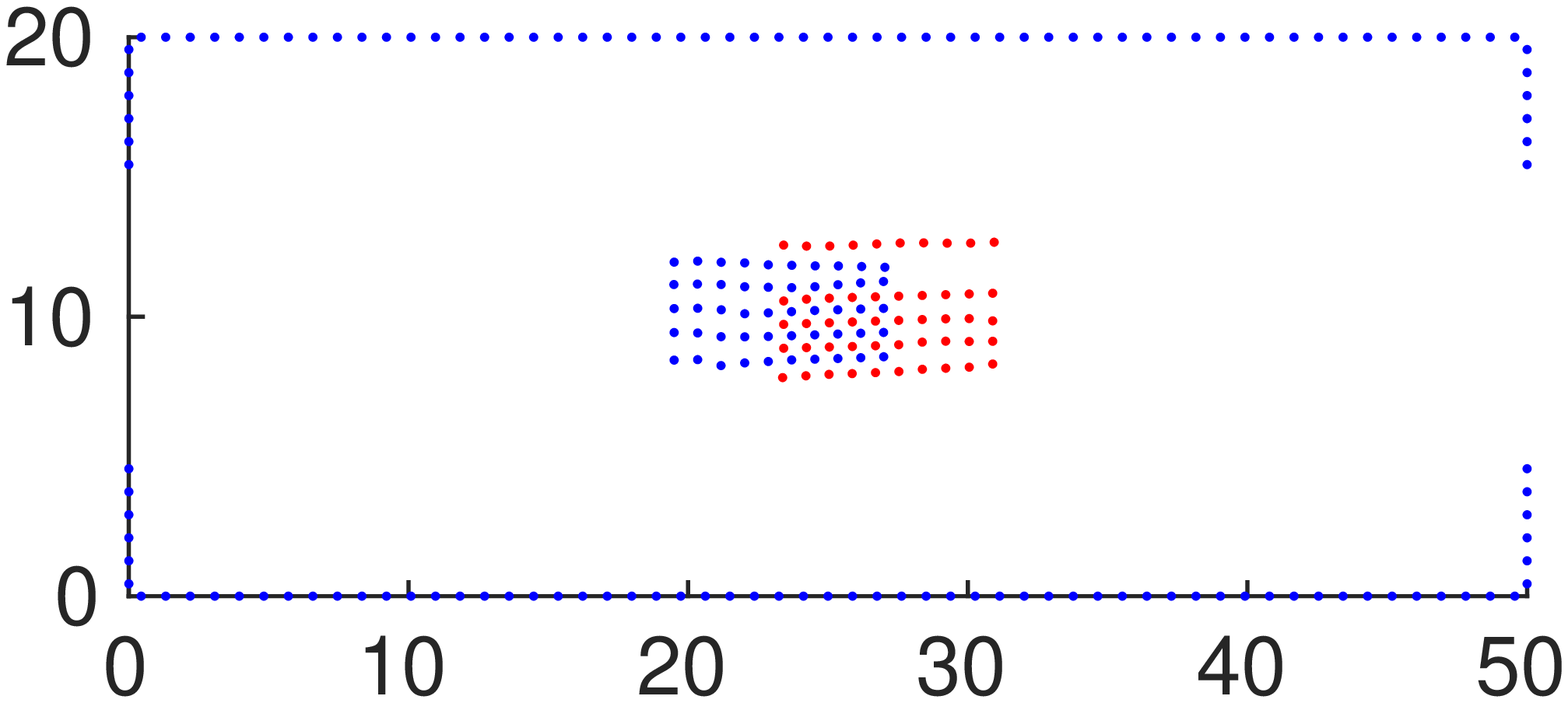}\\
	\includegraphics[width=.5\linewidth]{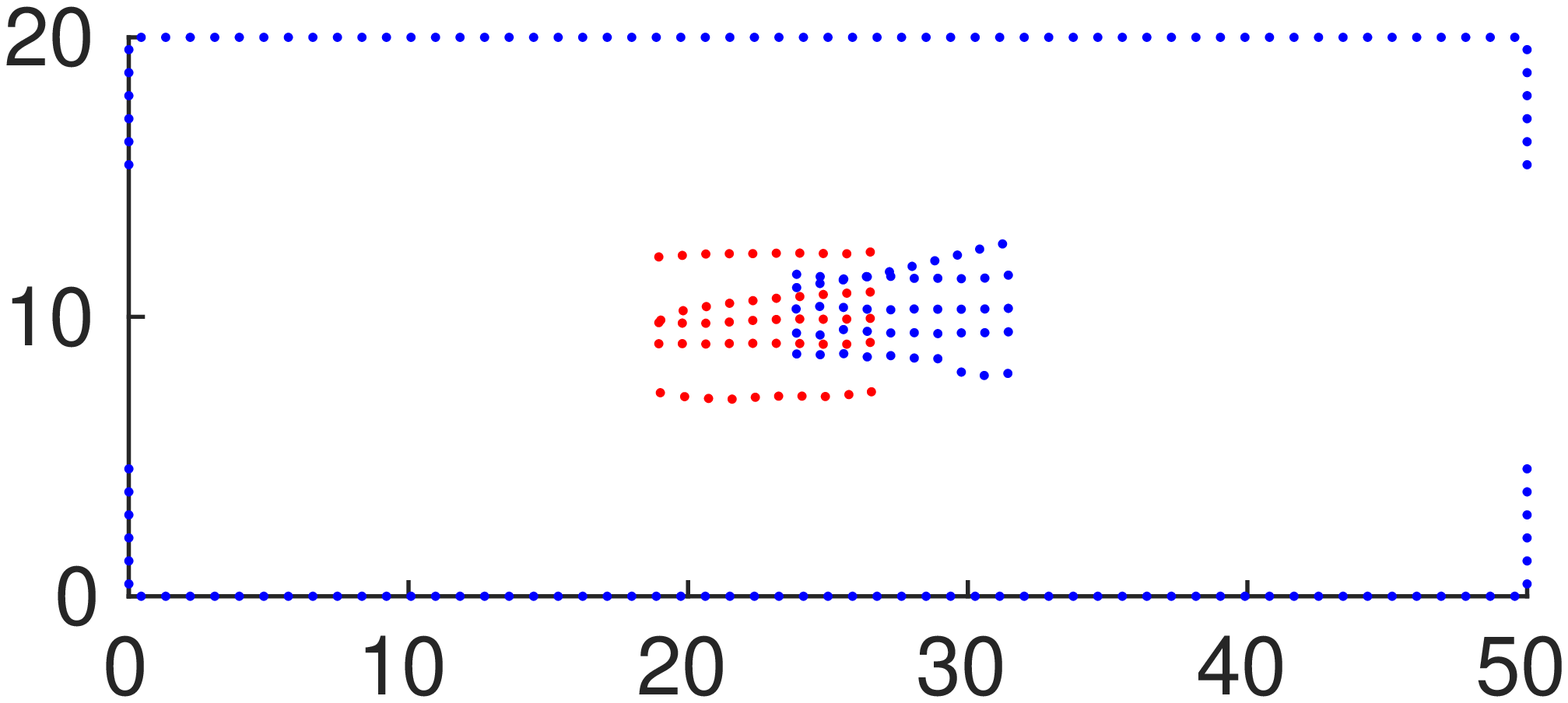}\hfill
	\includegraphics[width=.5\linewidth]{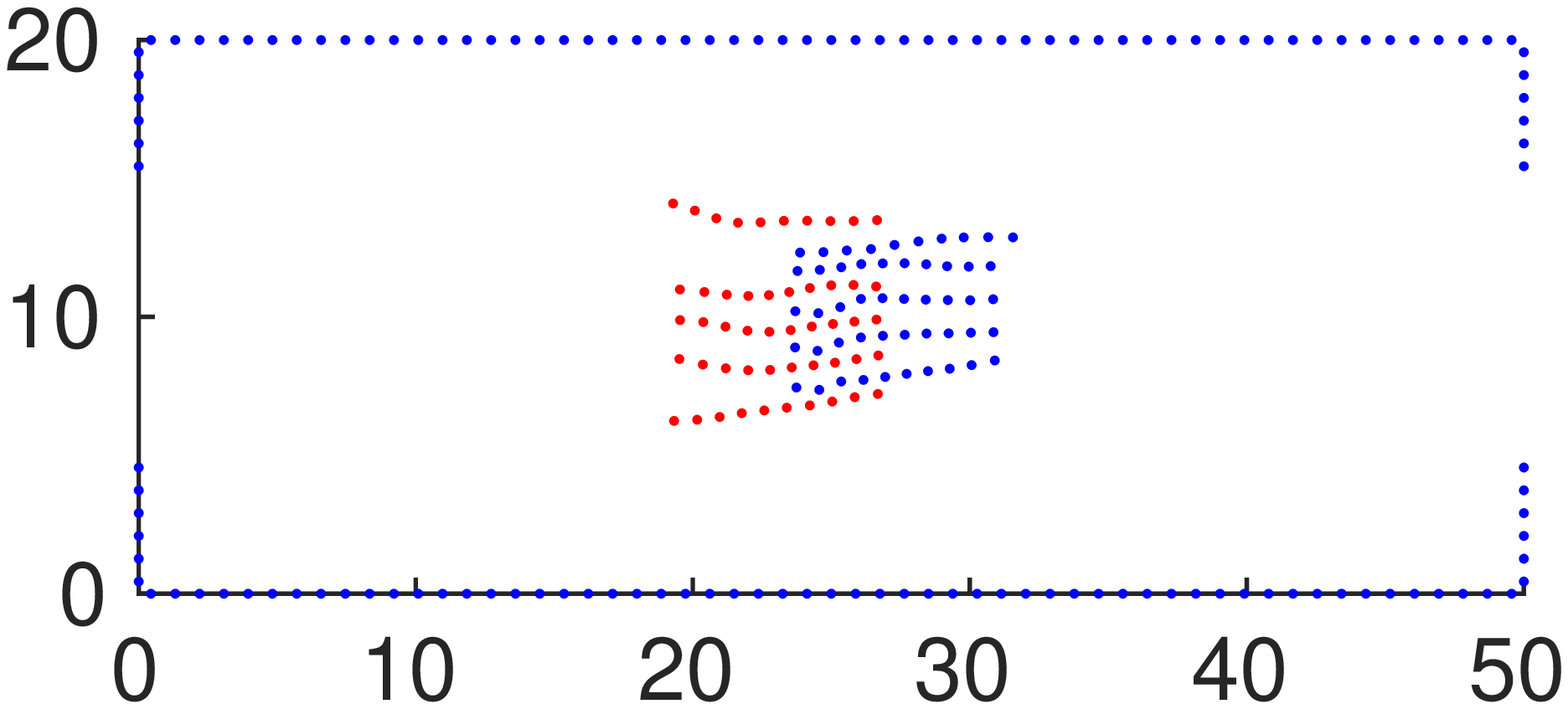}\\
	\includegraphics[width=.5\linewidth]{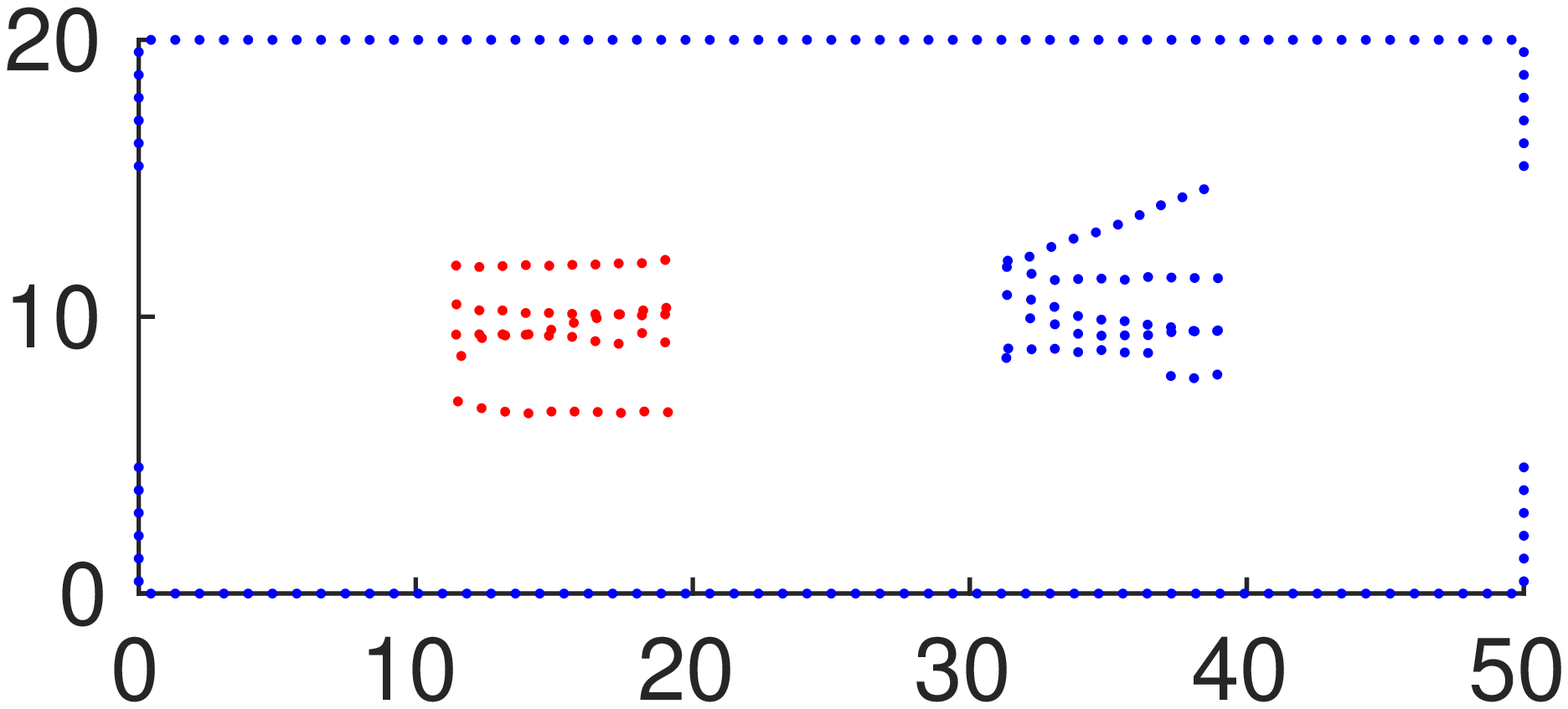}\hfill
	\includegraphics[width=.5\linewidth]{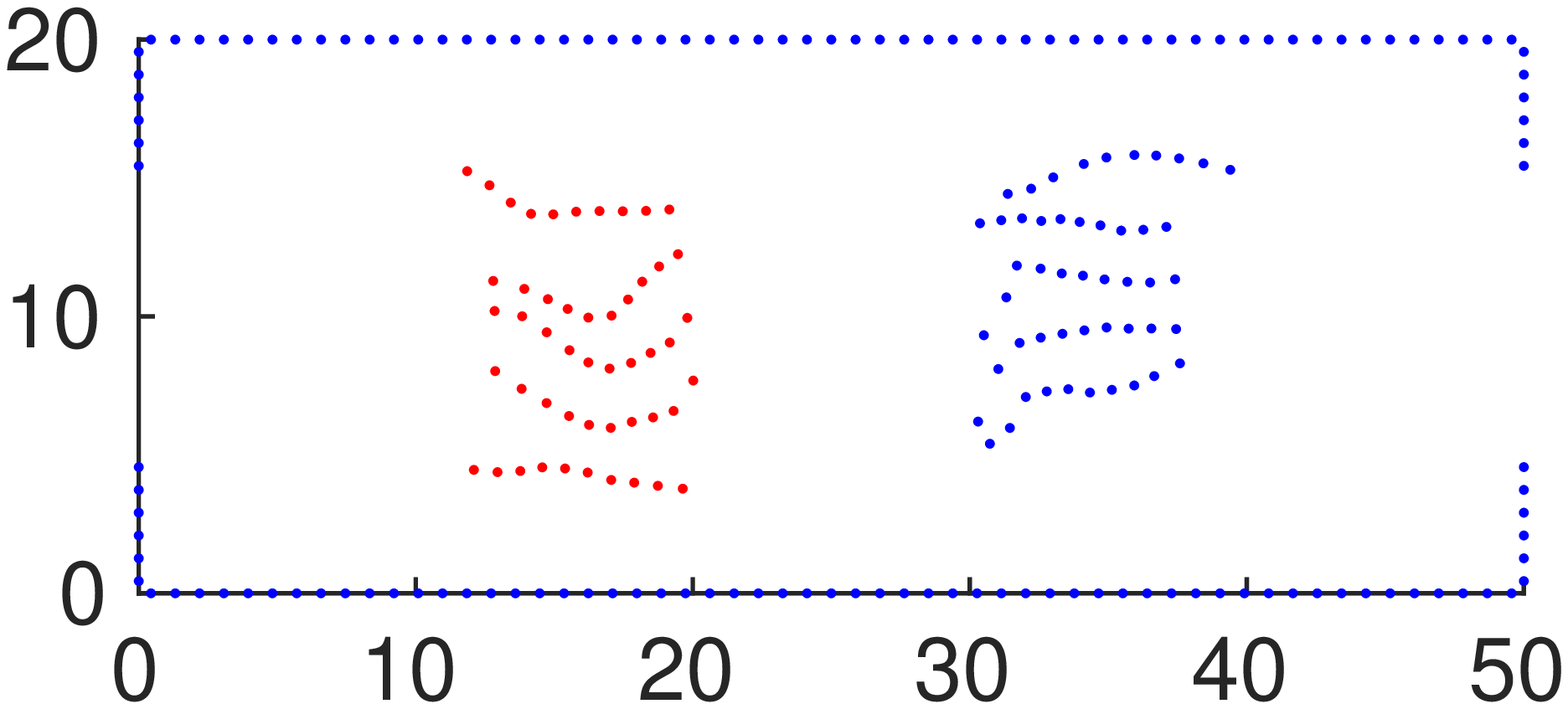}\\
	\includegraphics[width=.5\linewidth]{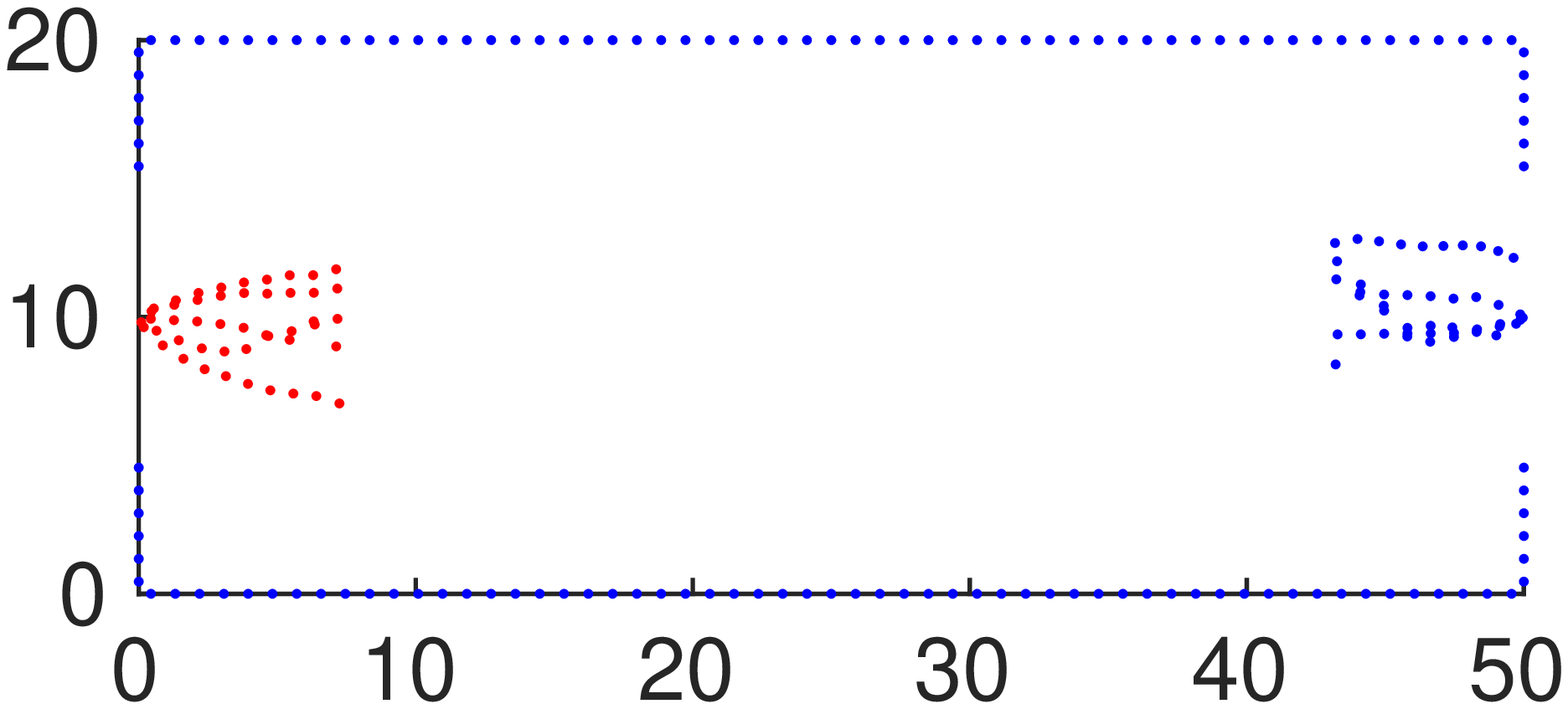}\hfill
	\includegraphics[width=.5\linewidth]{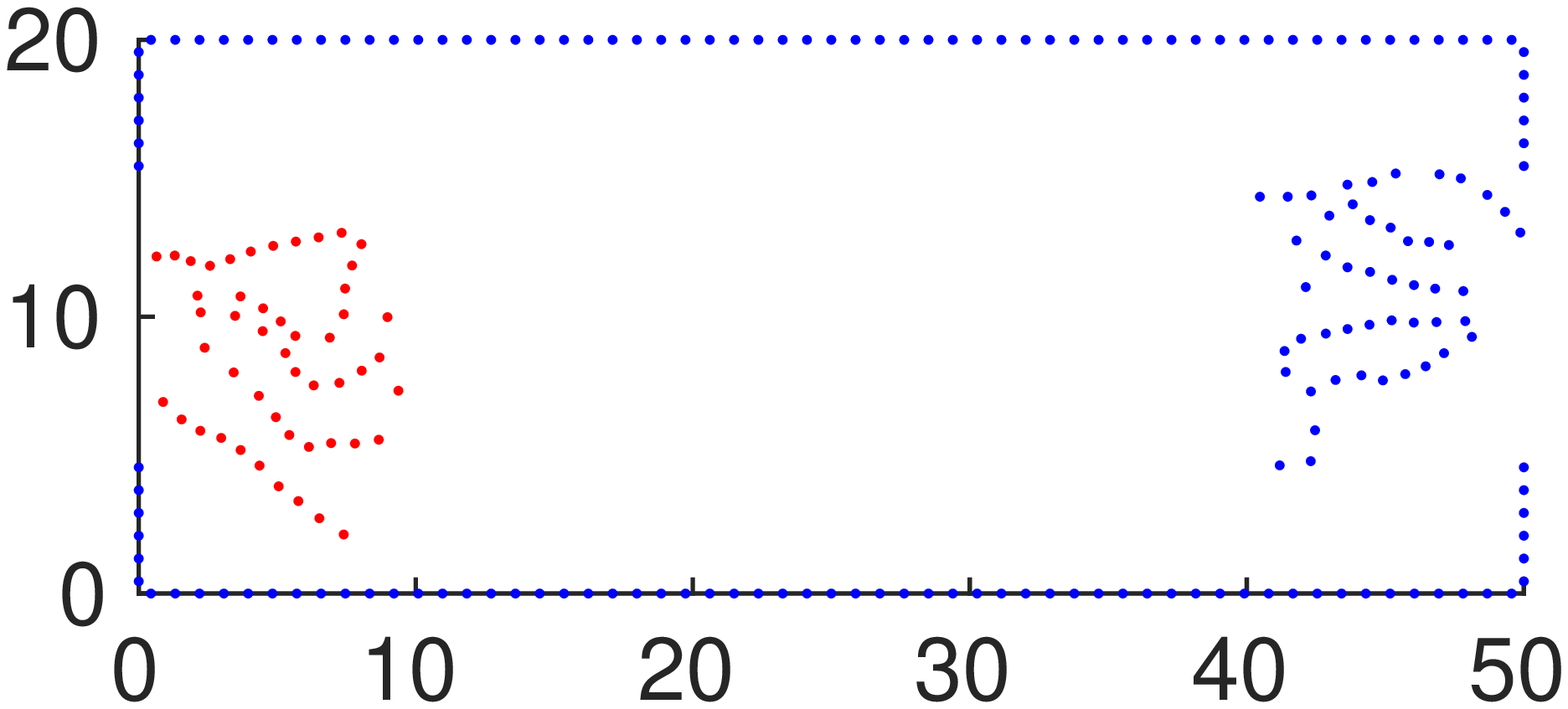}	
	\caption{Distribution of grid particles in the vision-based model (left column) and vision-based model with extra repulsion force (right column) at time $t = 10.5,12, 15, 20, 28$ (top to bottom).}
	\label{fig:7}
\end{figure}

Figure \ref{fig:7} shows the time evolution of the grid particles in the model without extra repulsive force between grid particles (left column) and with  extra repulsive force between grid particles (right column) at time $t = 10.5$, $t = 12$, $t = 15$, $t = 20$ and $t = 28$. If we do not use an additional repulsive force   we can clearly observe that some particles collide with each other. 


\subsection{Comparison of non-local and local approximation}\label{sec:5.3}
In this subsection we compare the  non-local approximation  with the local approximation for different $\lambda$. We use the same value of parameters as in subsection \ref{sec:5.1}.

\begin{figure}
	\includegraphics[width=.33\linewidth]{sys_with_dir_t_8.eps}\hfill
	\includegraphics[width=.33\linewidth]{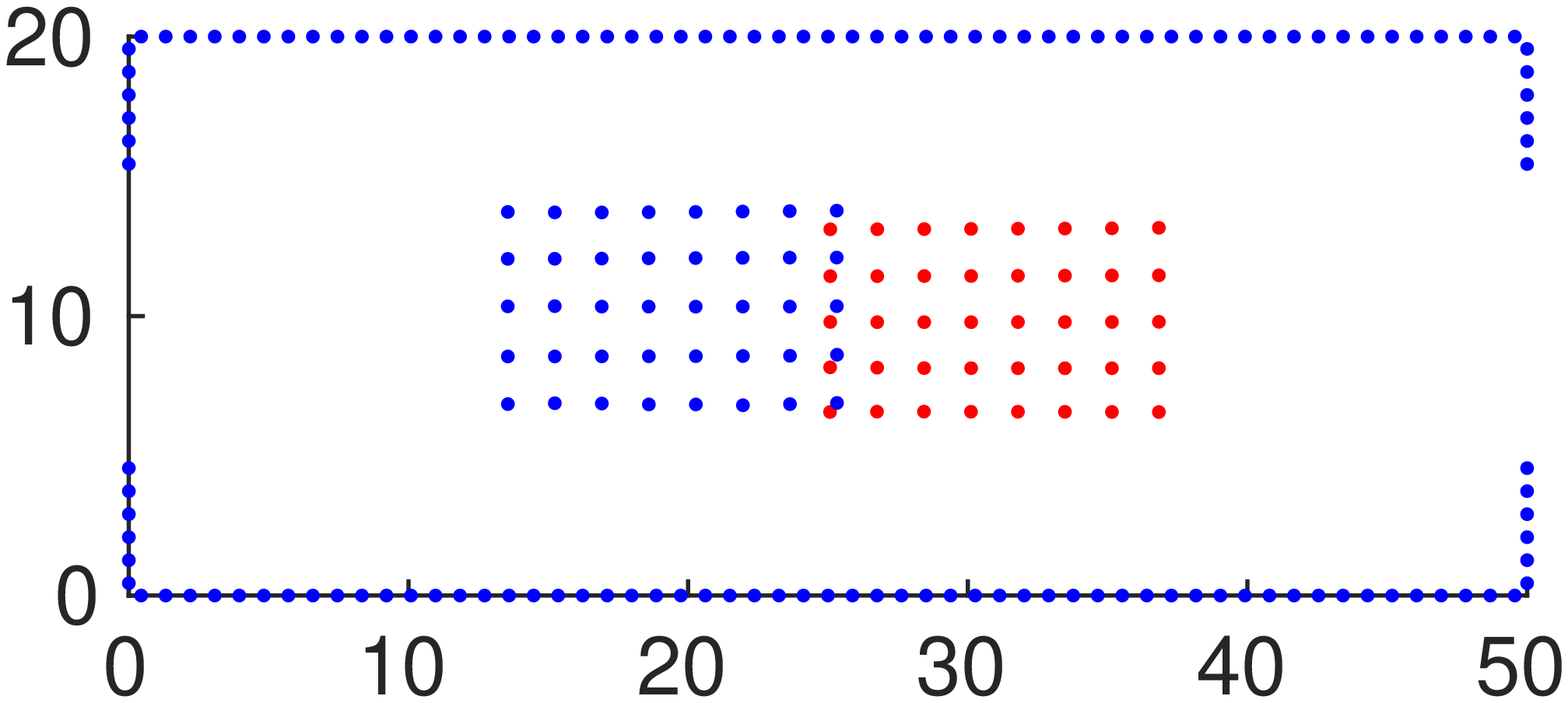}\hfill
	\includegraphics[width=.33\linewidth]{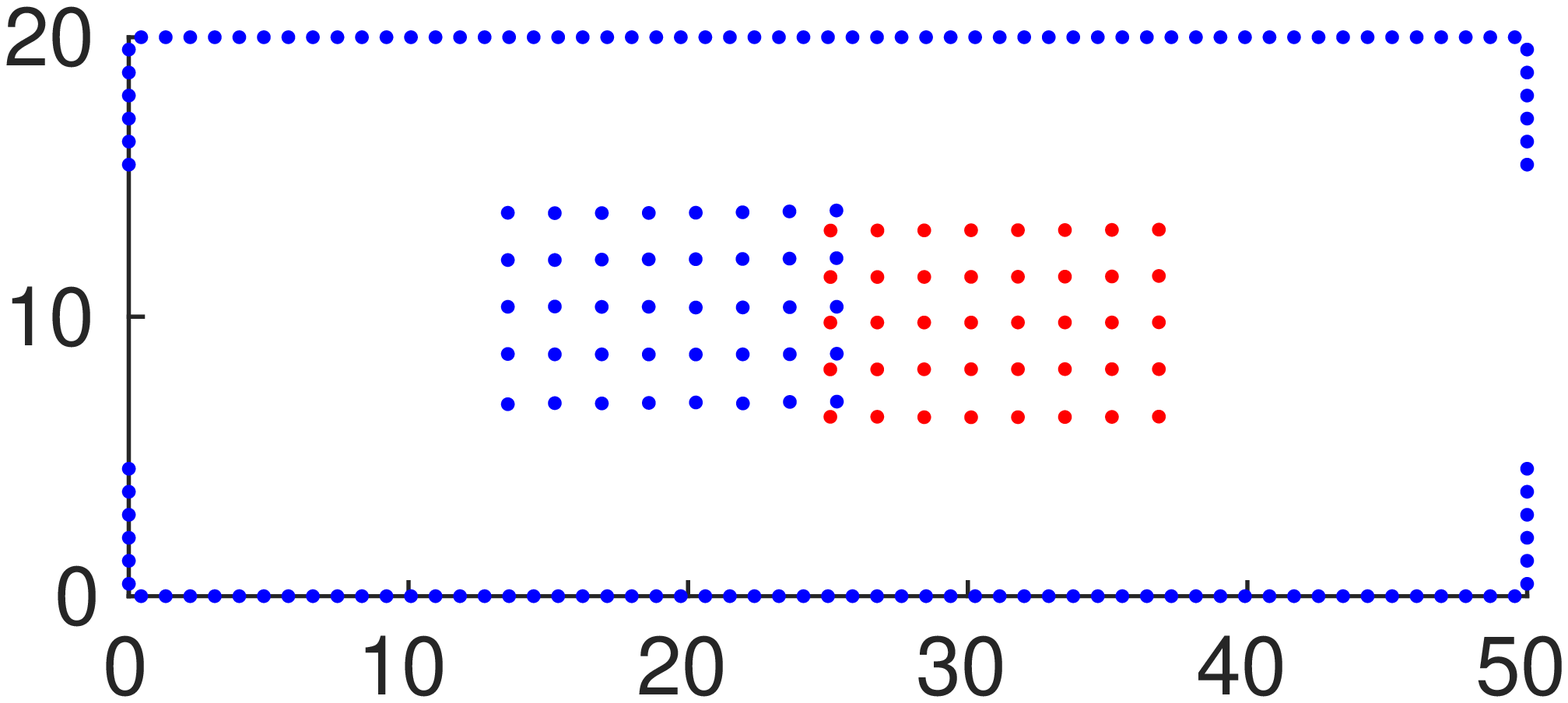}\\
	\includegraphics[width=.33\linewidth]{sys_with_dir_t_10.eps}\hfill
	\includegraphics[width=.33\linewidth]{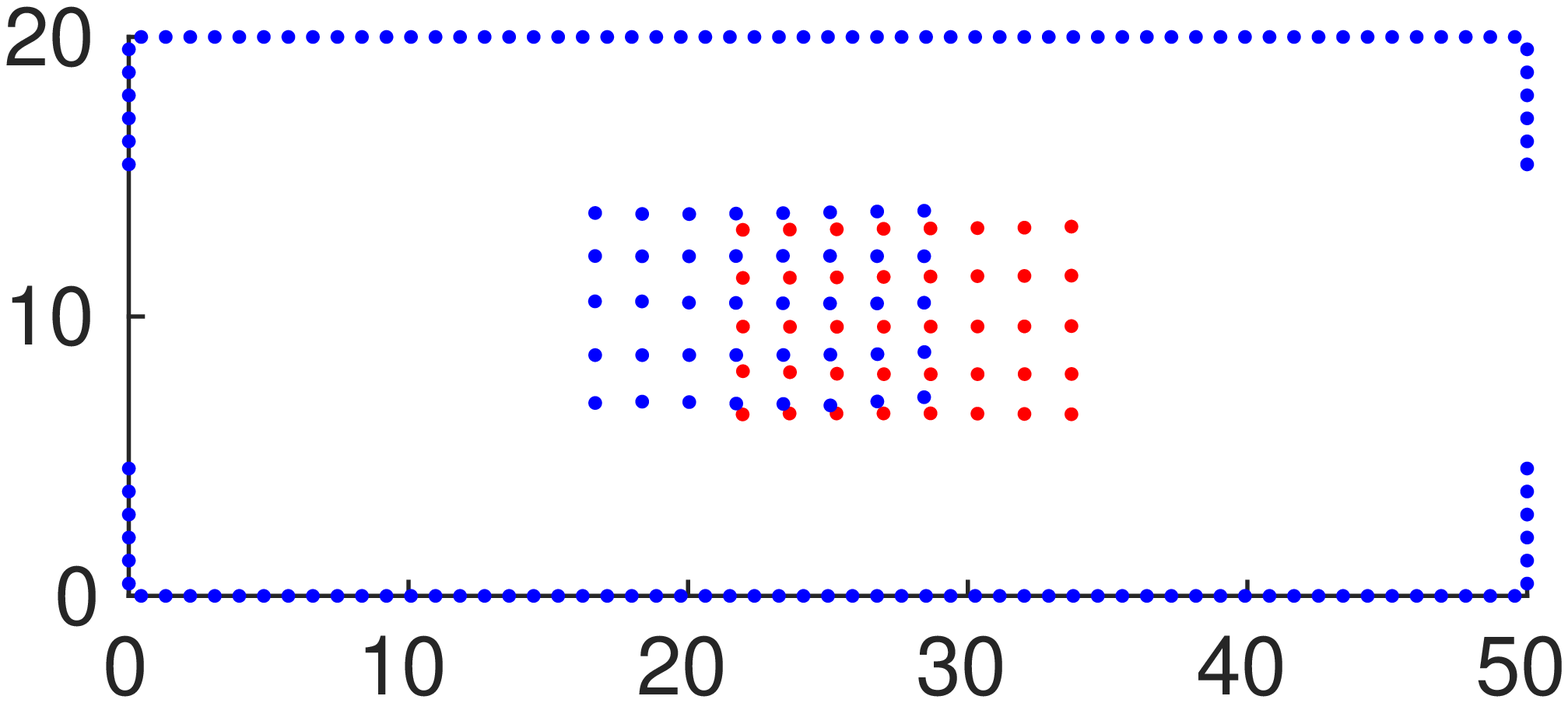}\hfill
	\includegraphics[width=.33\linewidth]{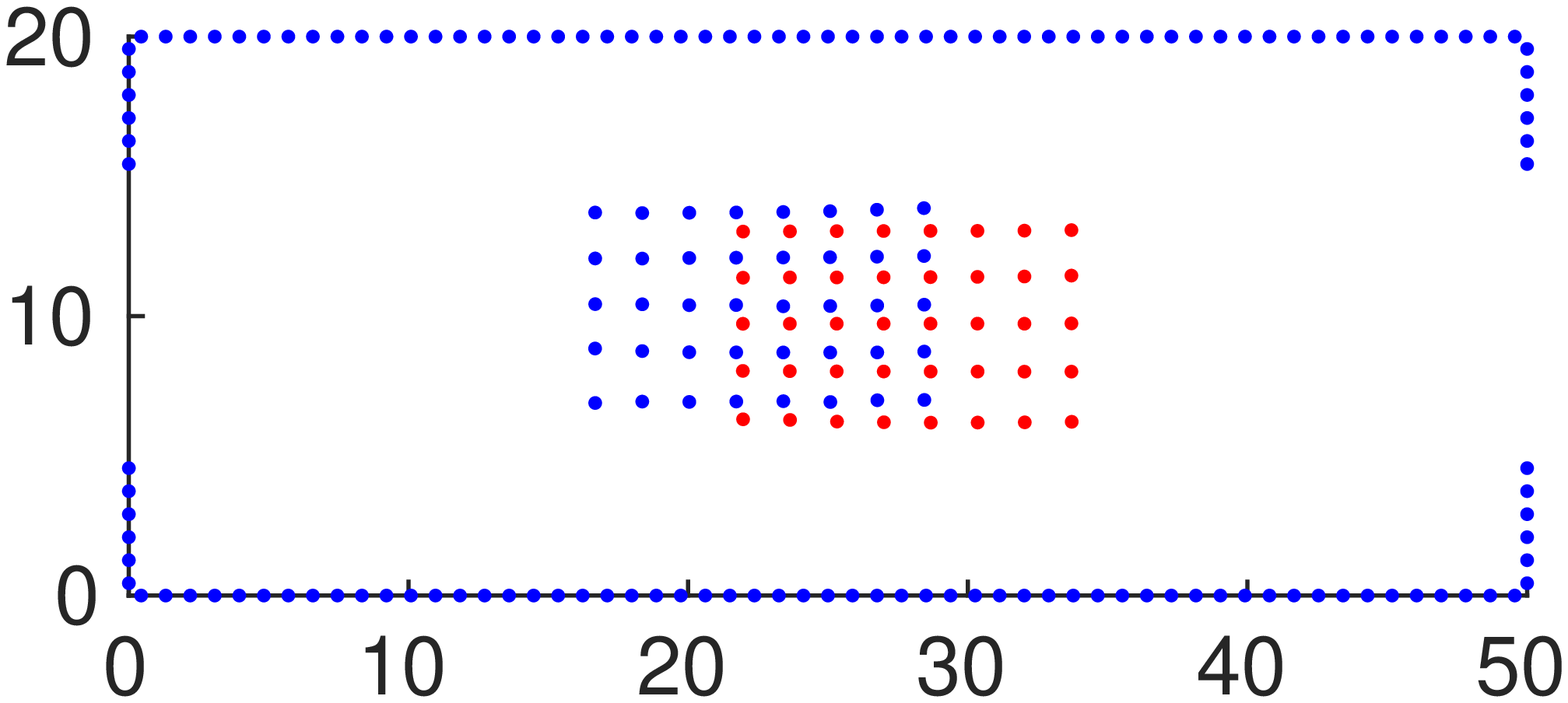}\\
	\includegraphics[width=.33\linewidth]{sys_with_dir_t_15.eps}\hfill
	\includegraphics[width=.33\linewidth]{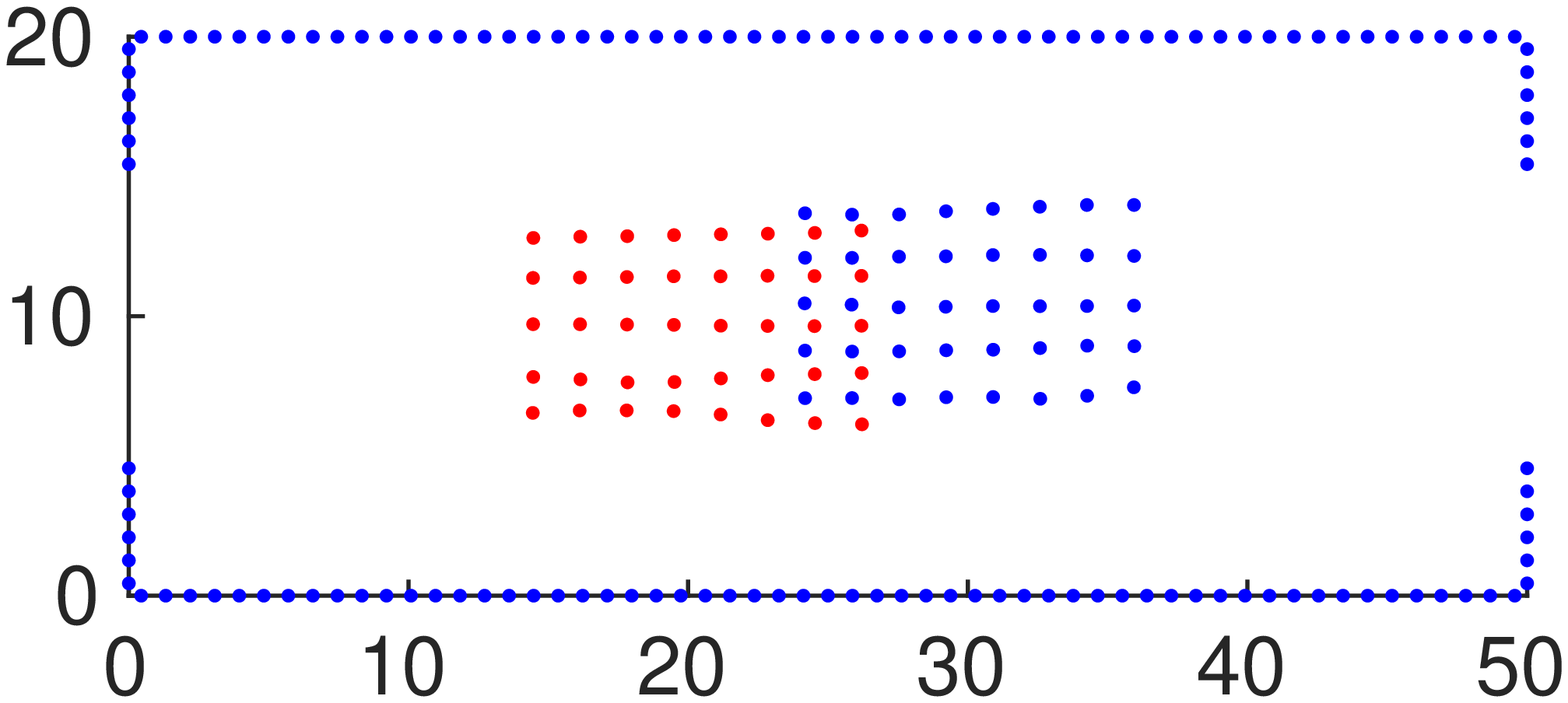}\hfill
	\includegraphics[width=.33\linewidth]{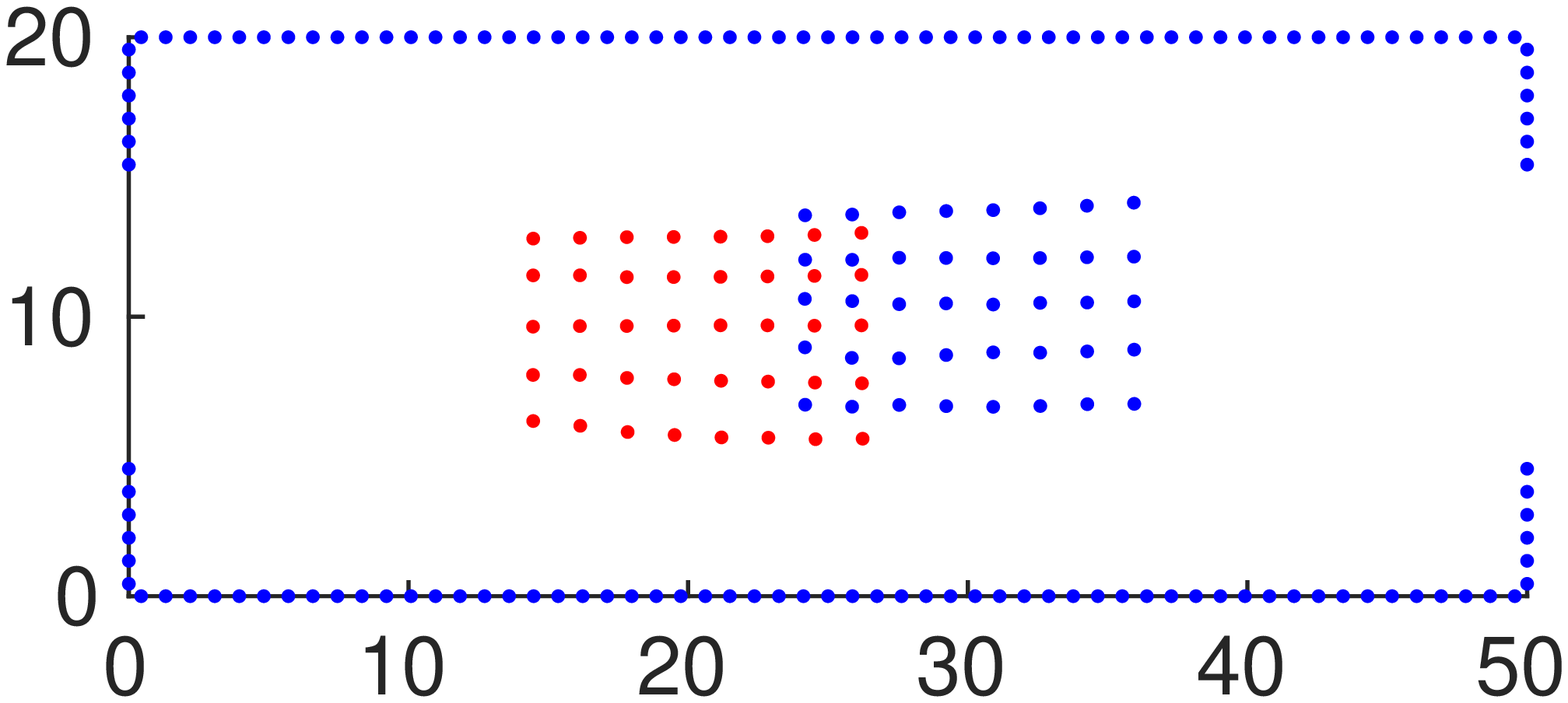}\\
	\includegraphics[width=.33\linewidth]{sys_with_dir_t_20.eps}\hfill
	\includegraphics[width=.33\linewidth]{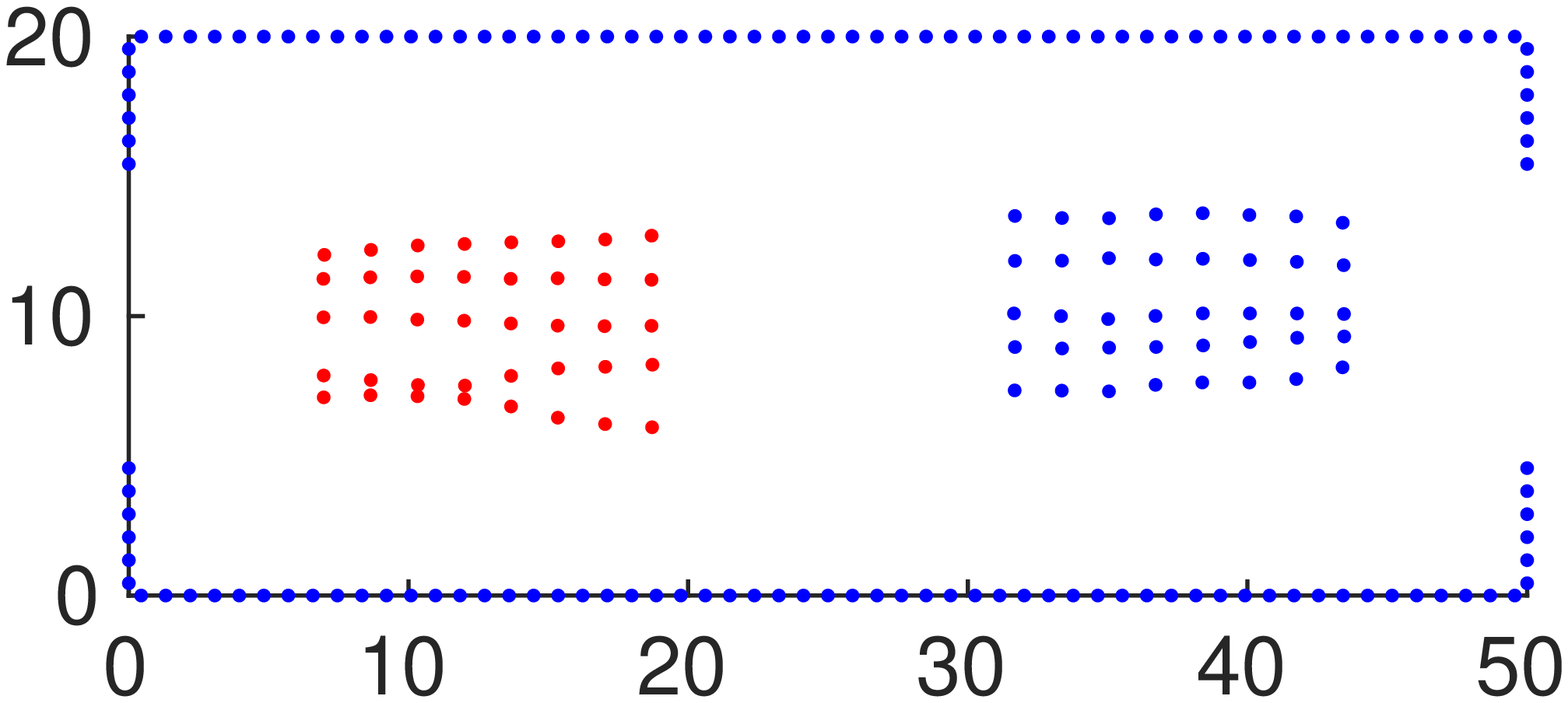}\hfill
	\includegraphics[width=.33\linewidth]{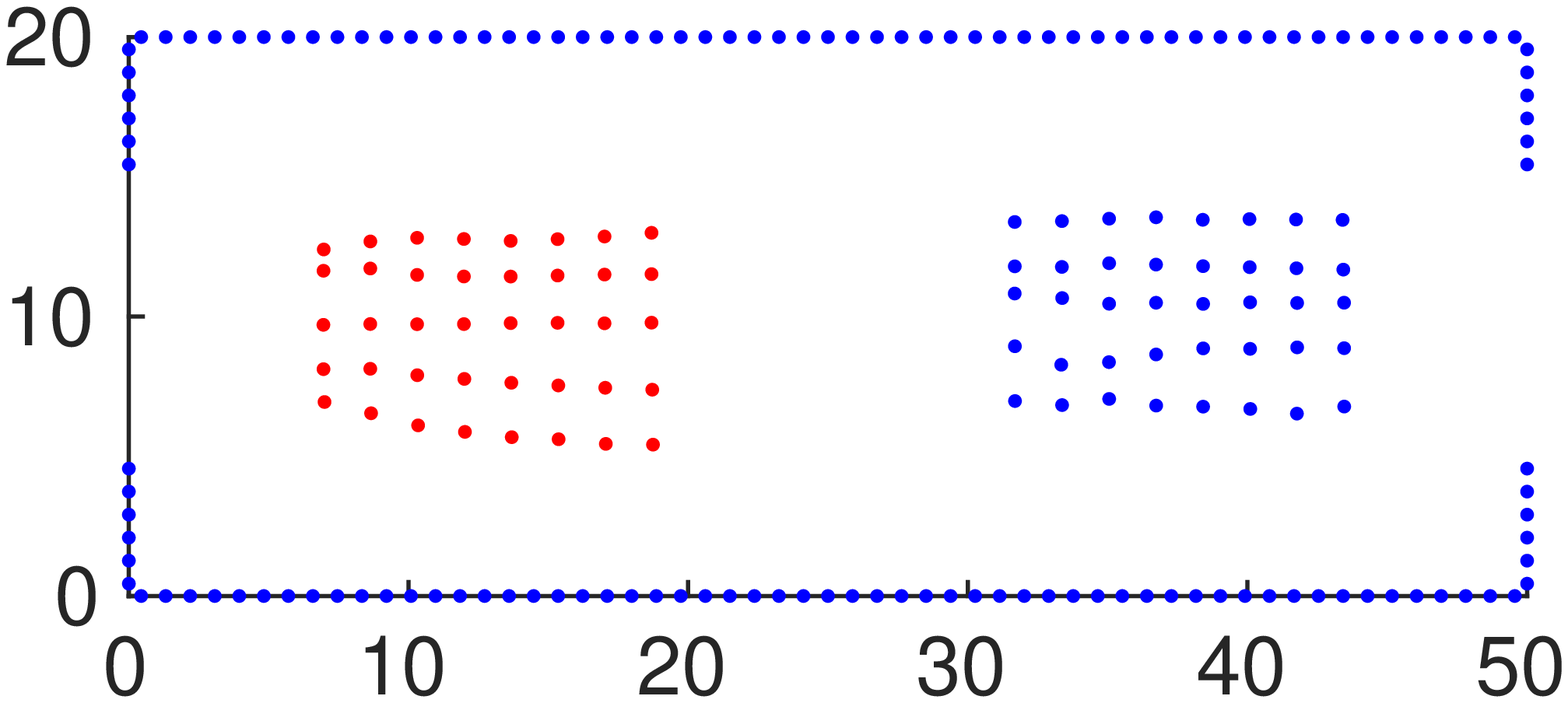}\\		
	\includegraphics[width=.33\linewidth]{sys_with_dir_t_25.eps}\hfill
	\includegraphics[width=.33\linewidth]{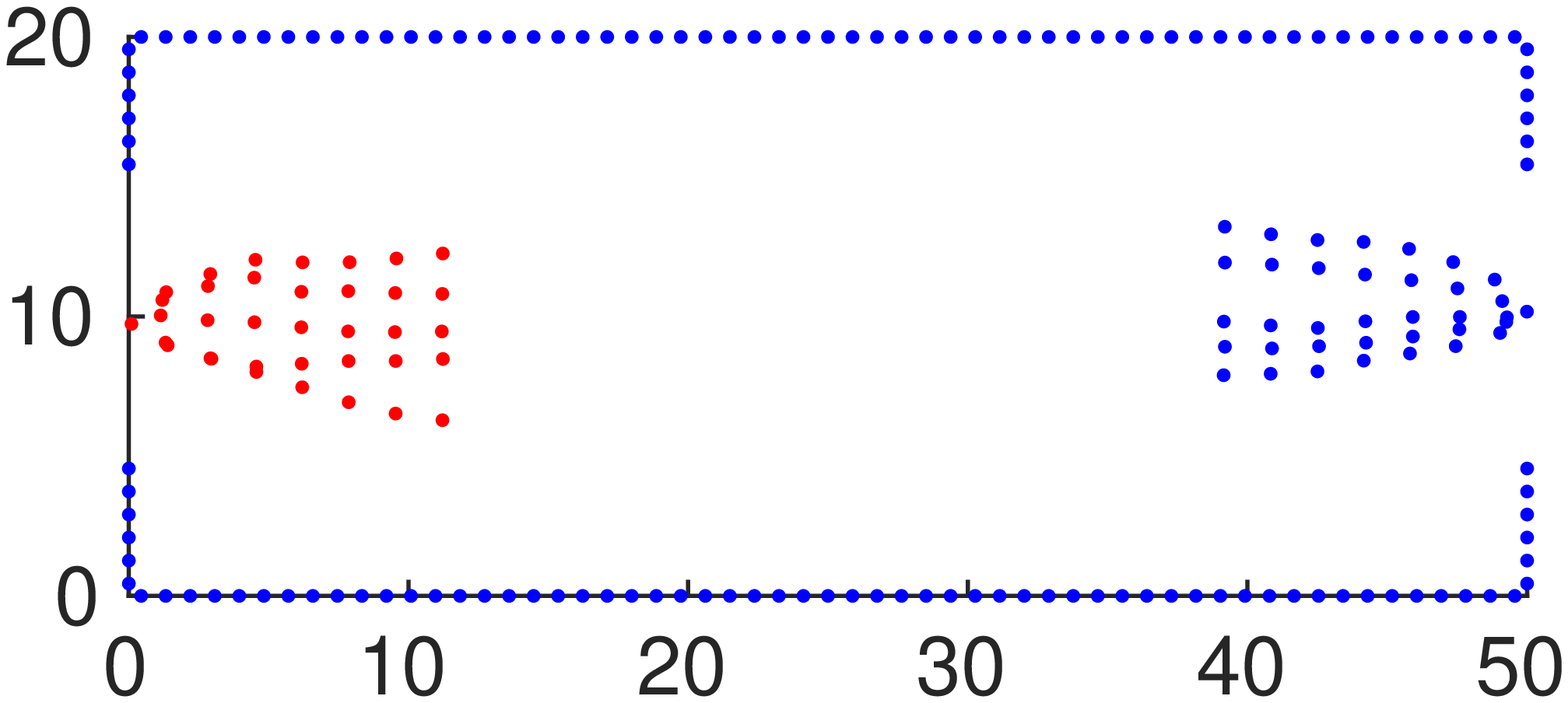}\hfill
	\includegraphics[width=.33\linewidth]{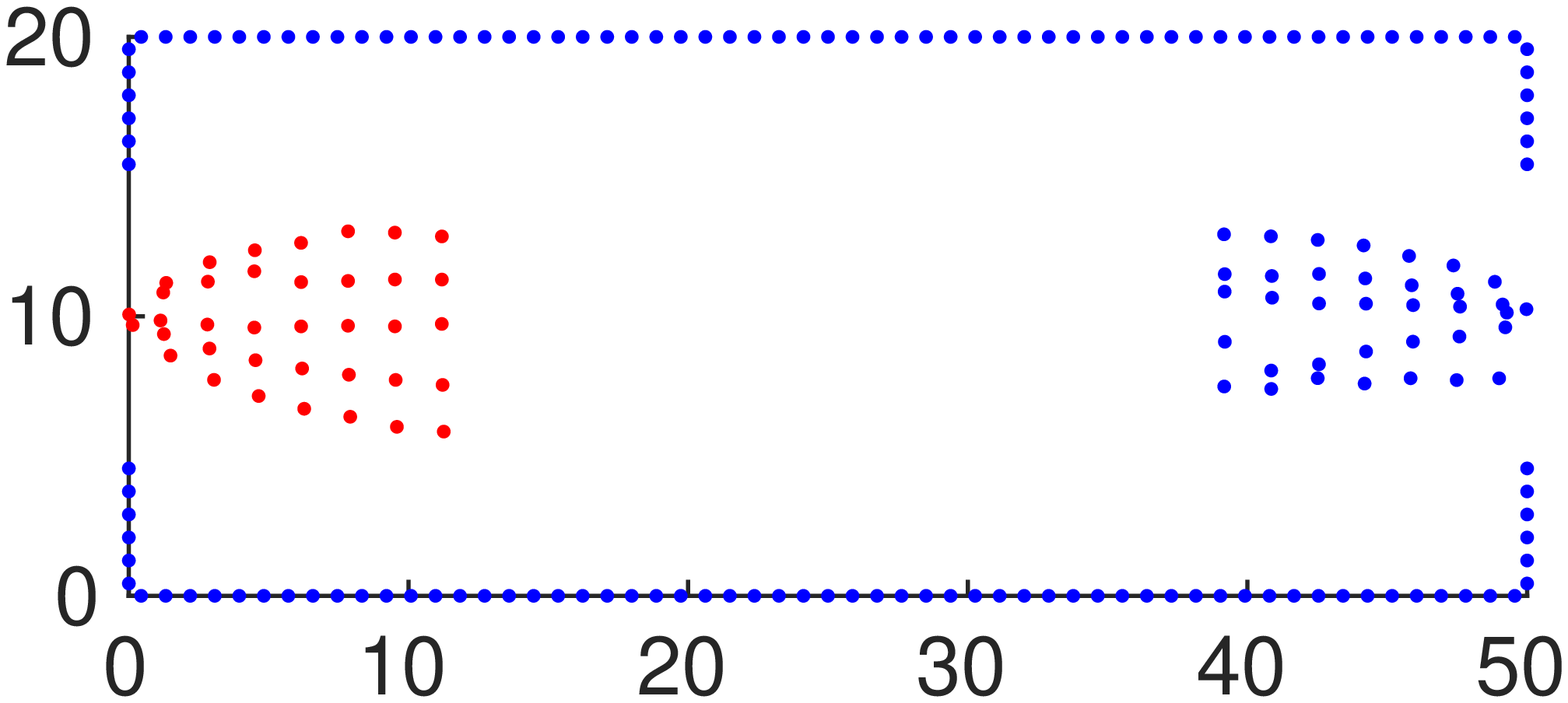}
	
	\caption{Distribution of grid particles in non-local (first column) and local for  $\lambda = 10^{-2}$ (second column),  $\lambda = 10^{-5}$ (third column) at time $t = 8, 10, 15, 20, 25$ (top to bottom) respectively.}
	\label{fig:8}
\end{figure}

\begin{figure}
	\includegraphics[scale=0.5]{density_with_dir_t_10.eps}\\
	\includegraphics[scale=0.5]{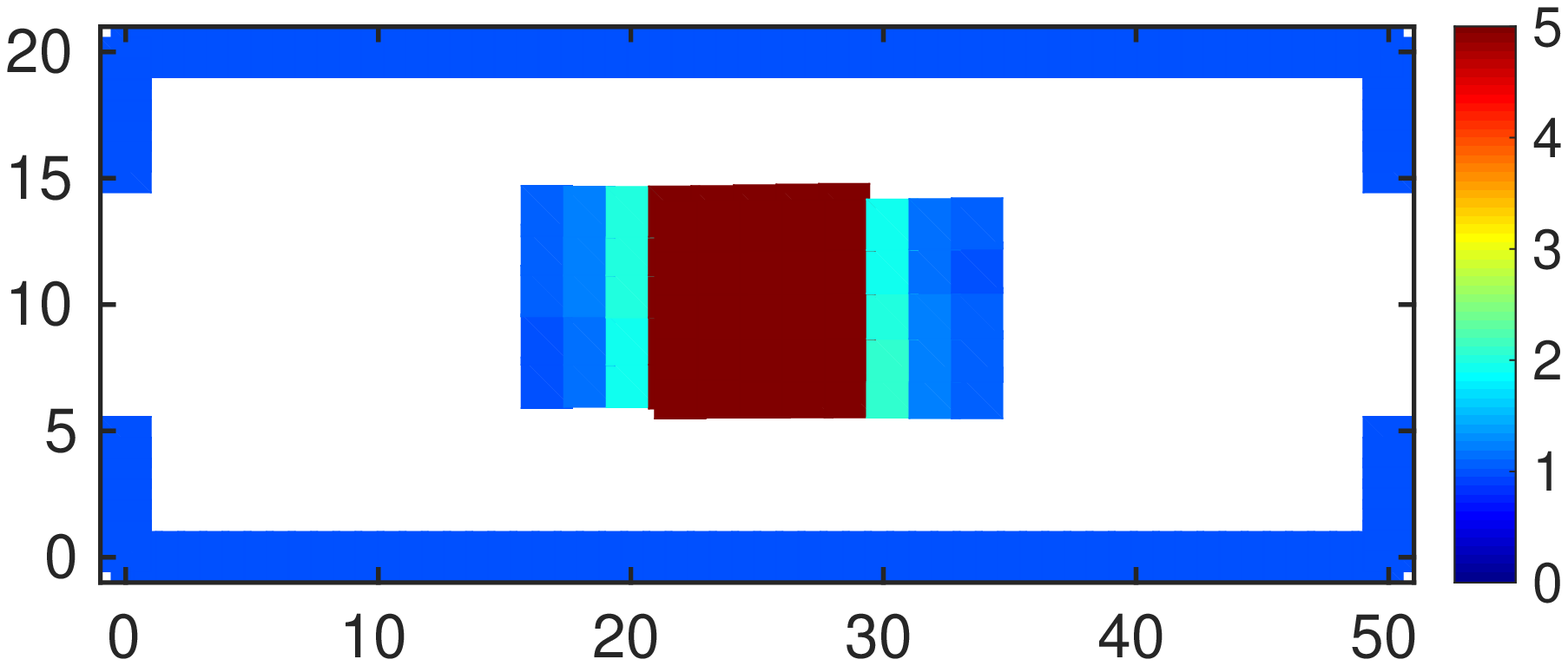}\\
	\includegraphics[scale=0.5]{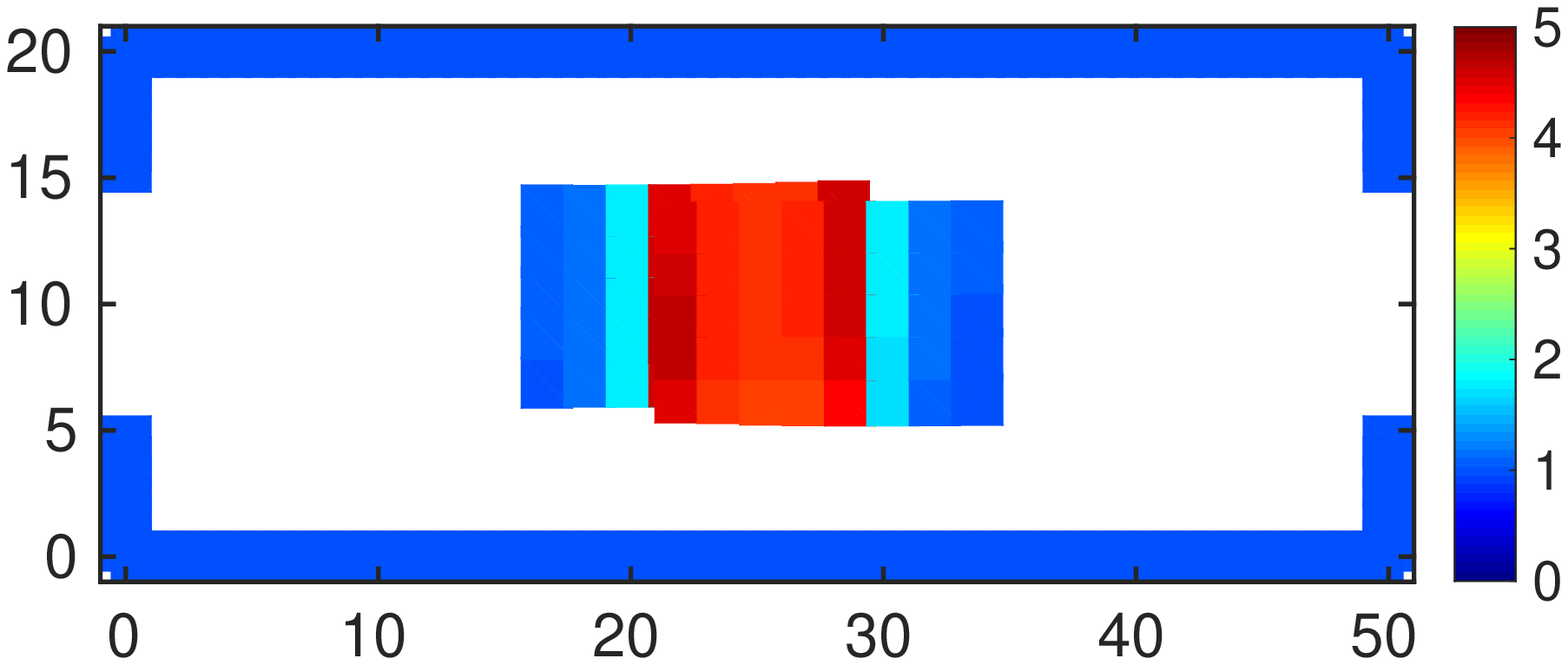}
	\caption{Density of pedestrians in the model for non-local (top) and local approximations for  $\lambda = 10^{-2}$ (middle),  $\lambda = 10^{-5}$ (bottom) at time $t = 10$.}
	\label{fig:9} 
\end{figure}

Figure \ref{fig:8} shows the time evolution of the grid particles in non-local (first column), local for $\lambda = 10^{-2}$ (second column) and local for $\lambda = 10^{-5}$ (third column) at time $t = 8$, $t = 10$, $t = 15$, $t = 20$ and $t = 25$. One observes that for  $\lambda$
going to $0$, the  the nonlocal  approximation approaches the local model.  Figure \ref{fig:9} shows the corresponding density plots at time $t = 10$ in the models for non-local and local approximations, where the higher density in the center of geometry indicates more collisions of pedestrians.

\begin{figure}[htbp]
	\centering
	\includegraphics[width=0.7\linewidth]{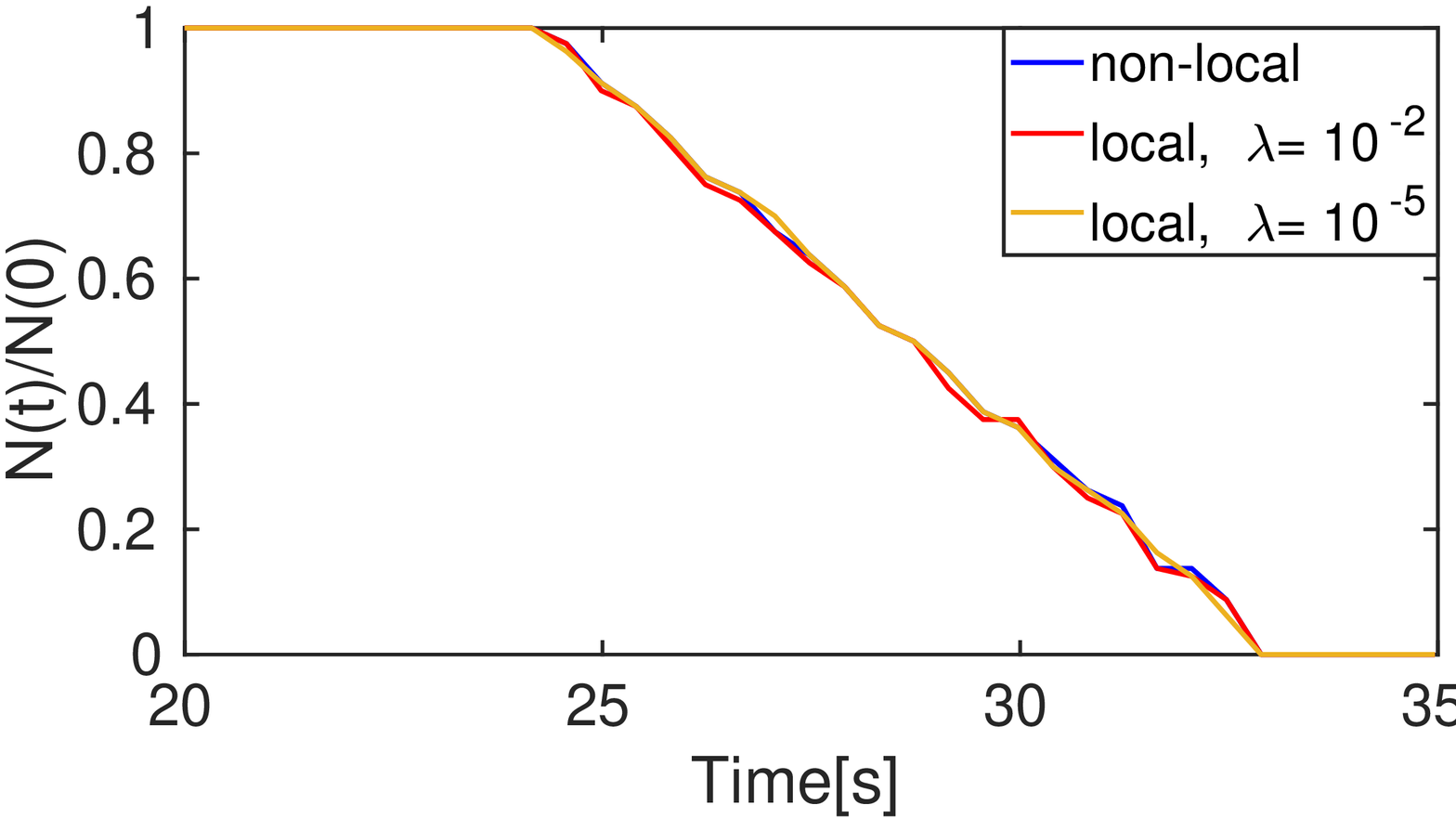}
	\caption{Ratio of initial and actual grid particles with respect to time in non-local and local approximation of the model with different $\lambda$.}
	\label{fig:10}
\end{figure}

Figure \ref{fig:10} shows the percentage of grid particles being in the computational domain in non-local and local approximation of the model with different $\lambda$ with respect to time. In this case also, one observes that the evacuation time is similar in all models.

\subsection{Comparison between vision-based pedestrian model and social force pedestrian model}\label{sec:5.4}
In this subsection we consider the same parameters for the vision-based models as in section \ref{sec:5.1}. We compare here the numerical simulation of the macroscopic vision-based pedestrian model and a  social force based macroscopic pedestrian model which is coupled with the Eikonal equation. For details of the macroscopic  social force model, we refer to \cite{etikyala}. Compare the results here also with the results obtained in \cite{ondrej} for  the microscopic vision-based models. 

For the social force model, we use the following values of parameters: $\Delta t = 0.00042$, $\gamma_n = 1.0$, $\gamma_t = 0.2$, $T = 0.001$, $u_{\max} = 1.5$ $ m/s$ and different values of repulsive force coefficient $k_n$ as 100, 500, 1000, 2000. 

\begin{figure}
	\includegraphics[width=.33\linewidth]{sys_with_dir_t_8.eps}\hfill
	\includegraphics[width=.33\linewidth]{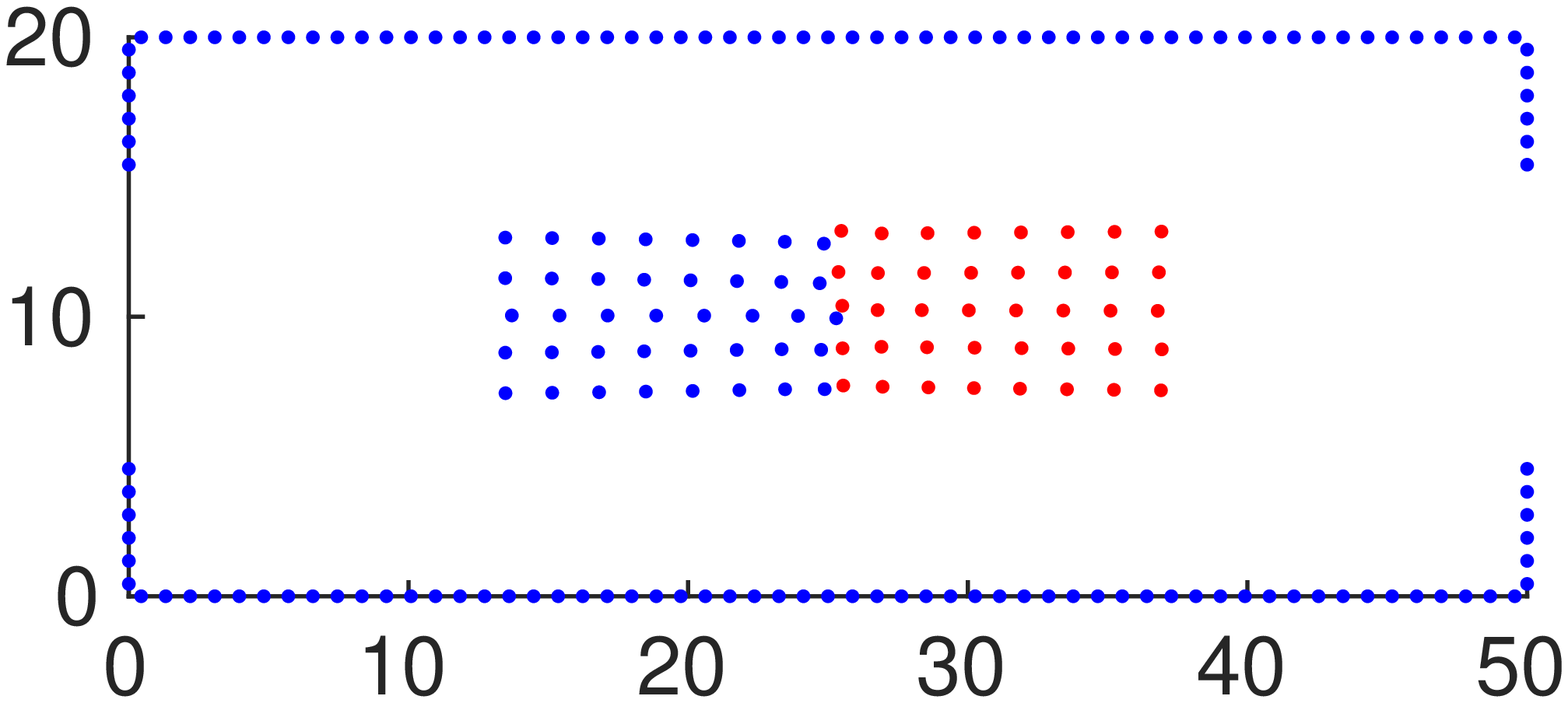}\hfill
	\includegraphics[width=.33\linewidth]{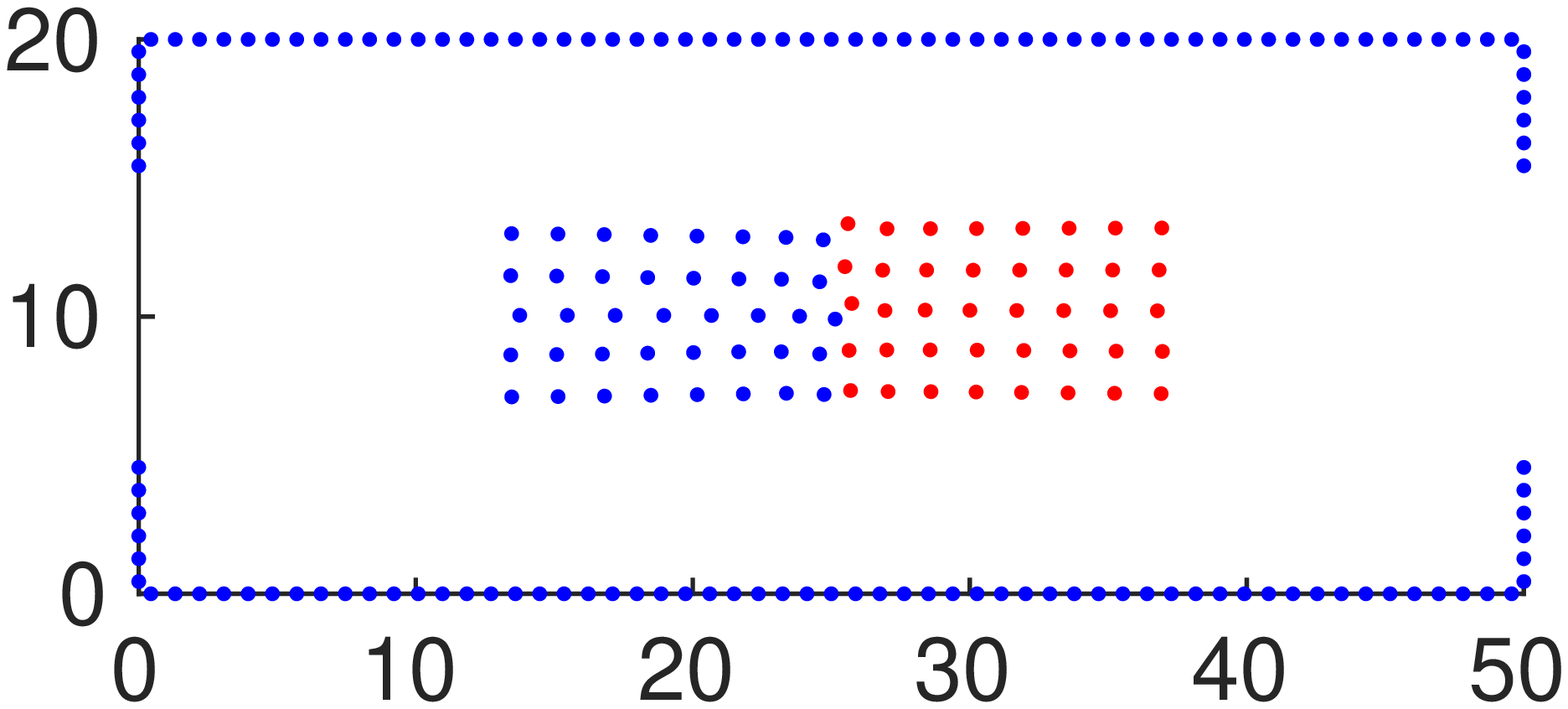}\\
	\includegraphics[width=.33\linewidth]{sys_with_dir_t_10.eps}\hfill
	\includegraphics[width=.33\linewidth]{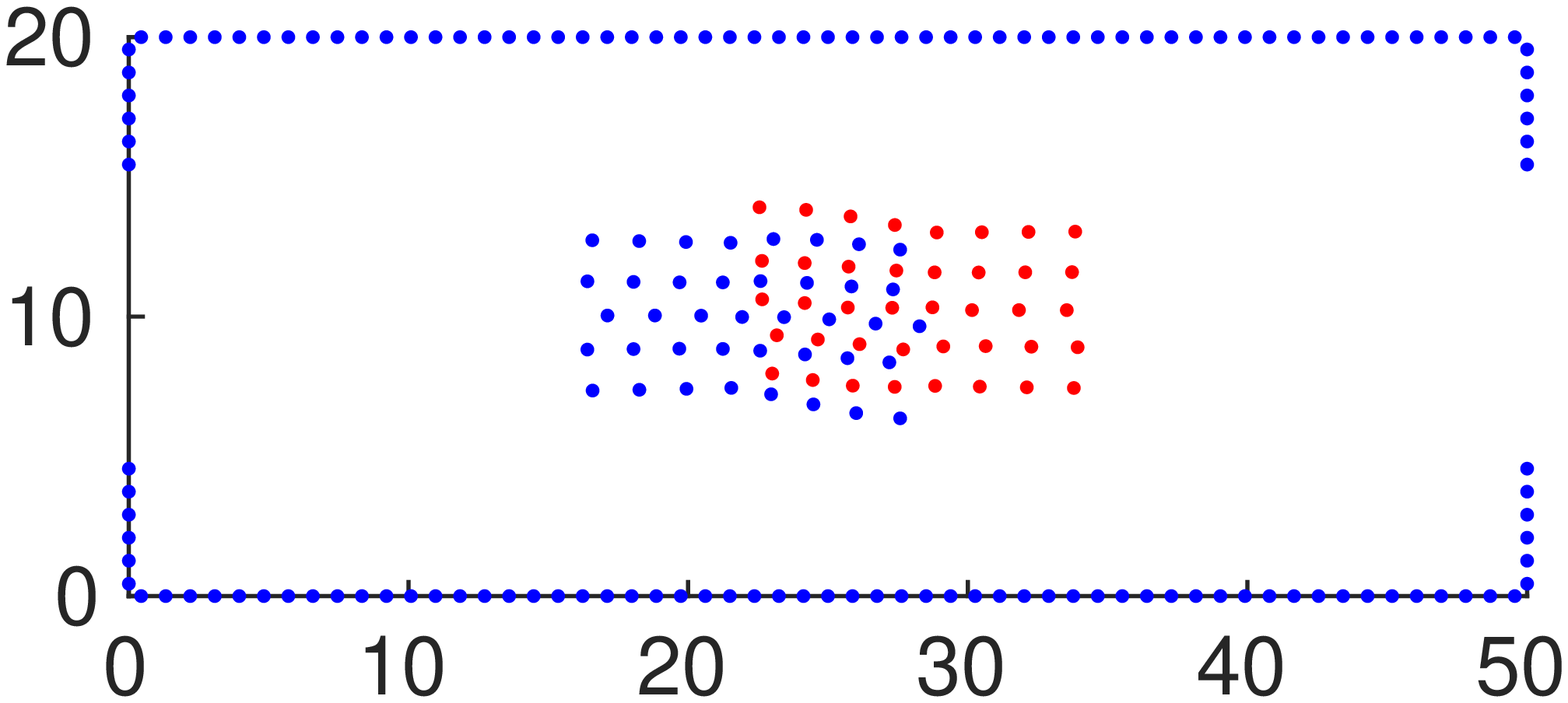}\hfill
	\includegraphics[width=.33\linewidth]{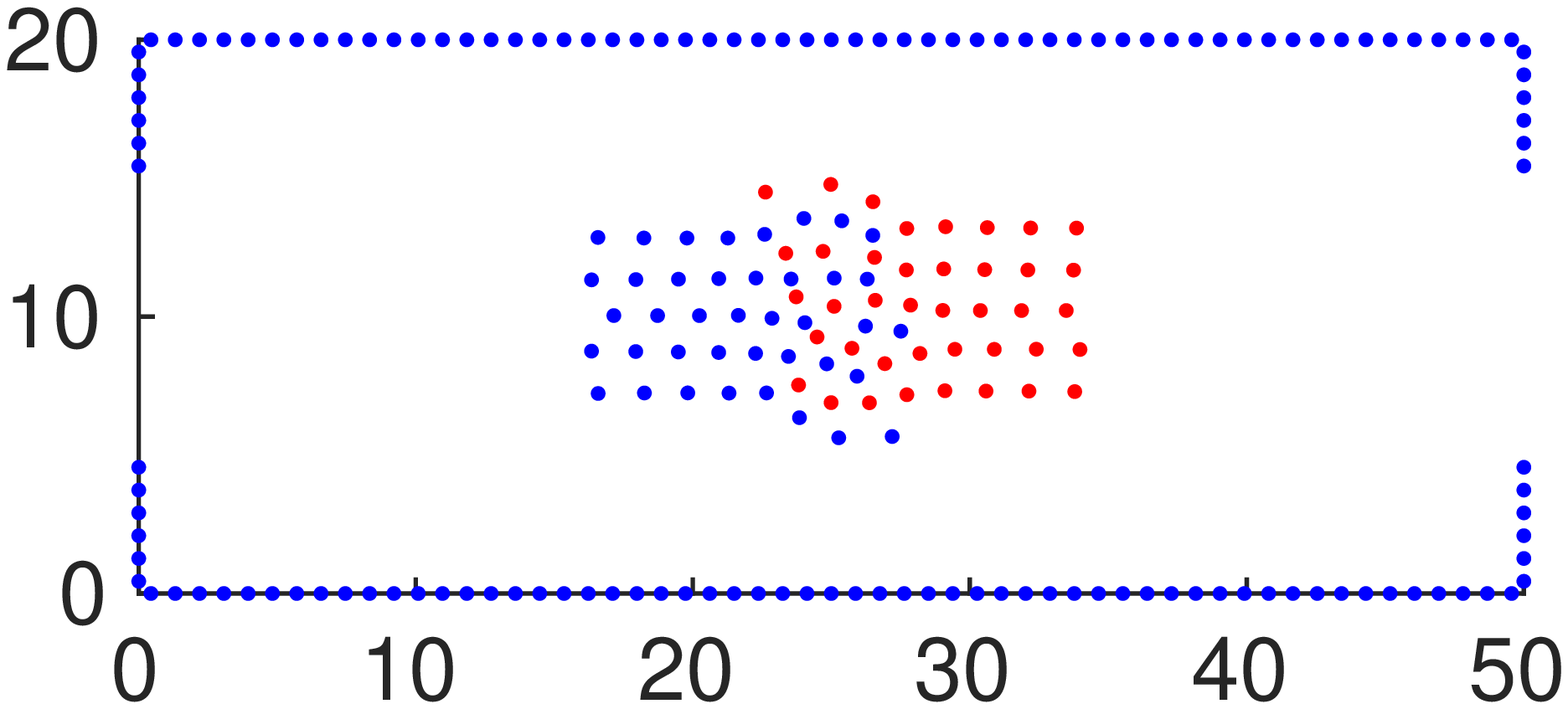}\\
	\includegraphics[width=.33\linewidth]{sys_with_dir_t_15.eps}\hfill
	\includegraphics[width=.33\linewidth]{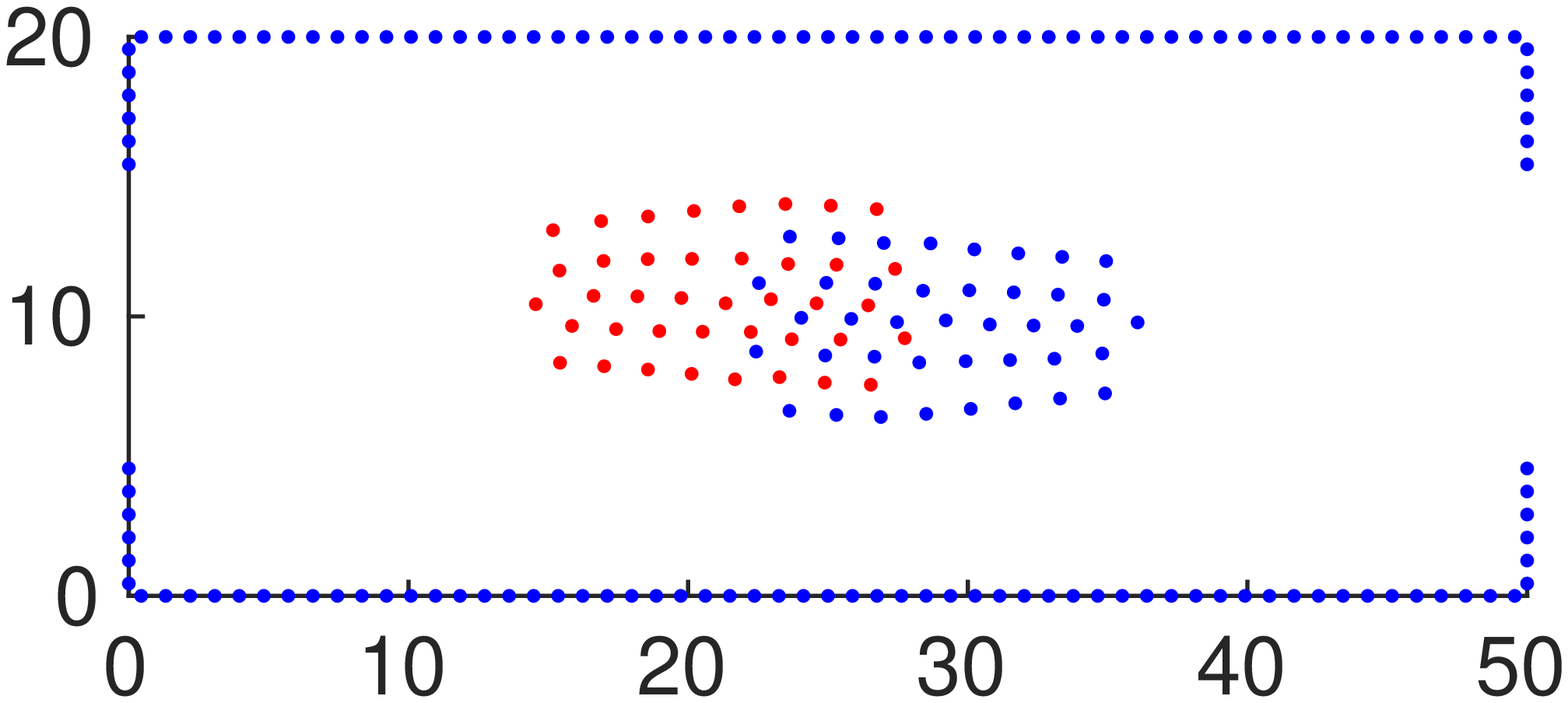}\hfill
	\includegraphics[width=.33\linewidth]{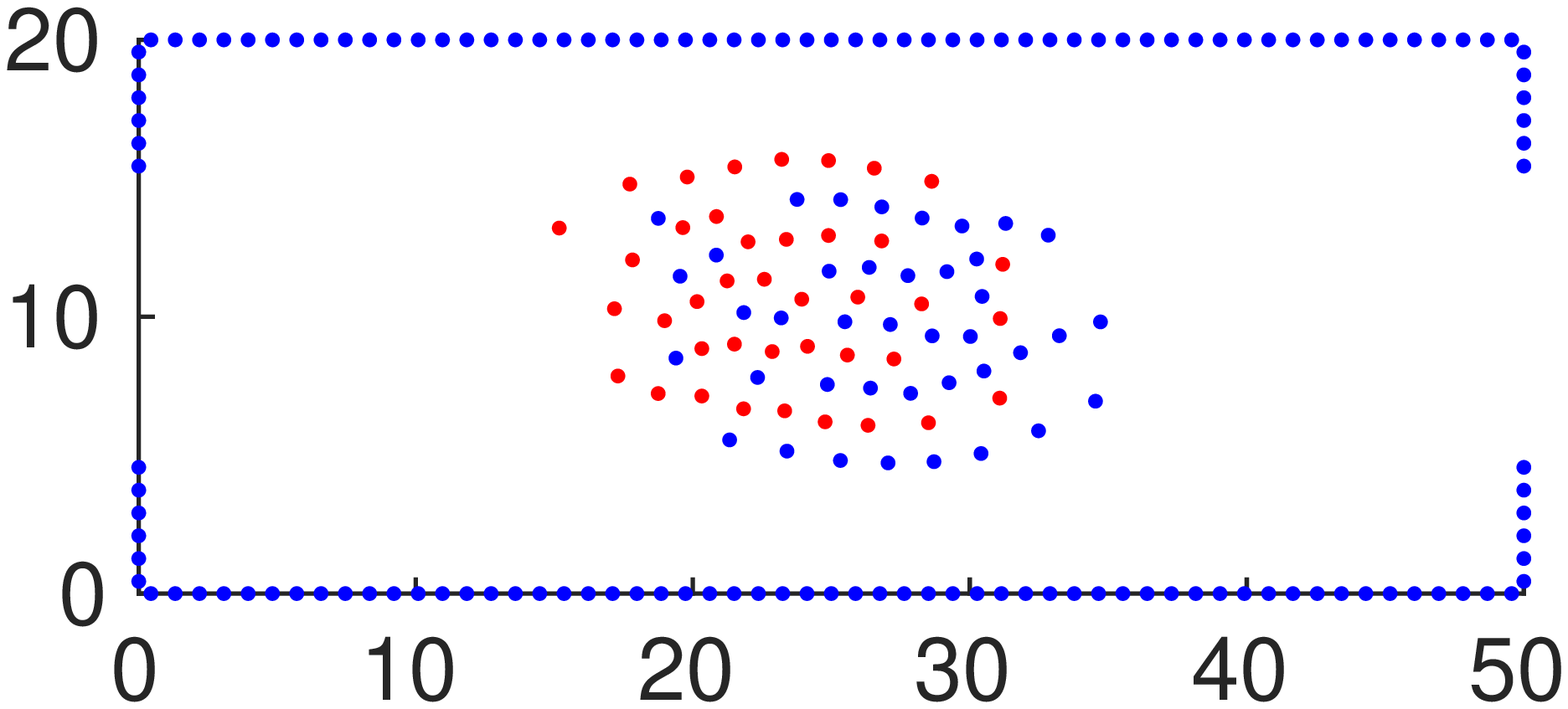}\\
	\includegraphics[width=.33\linewidth]{sys_with_dir_t_20.eps}\hfill
	\includegraphics[width=.33\linewidth]{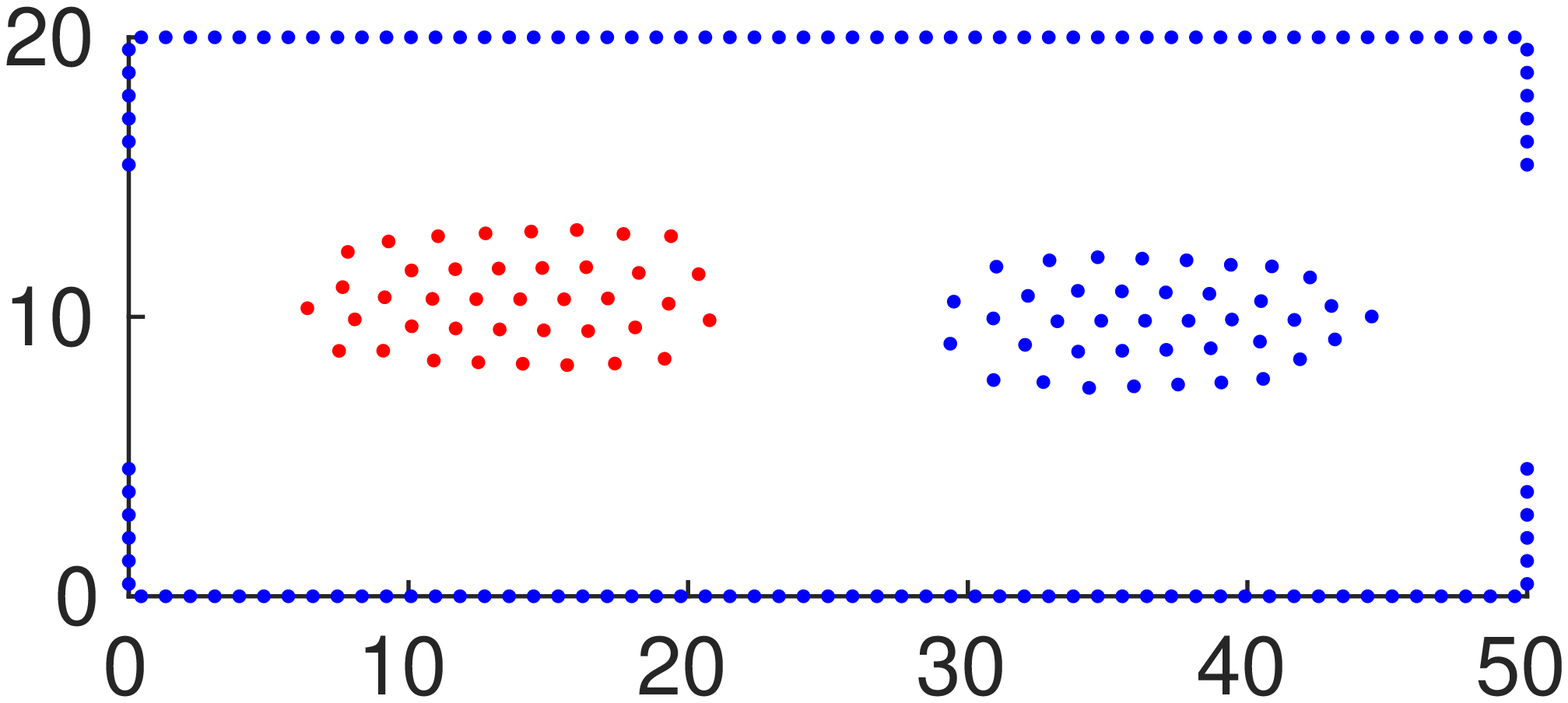}\hfill
	\includegraphics[width=.33\linewidth]{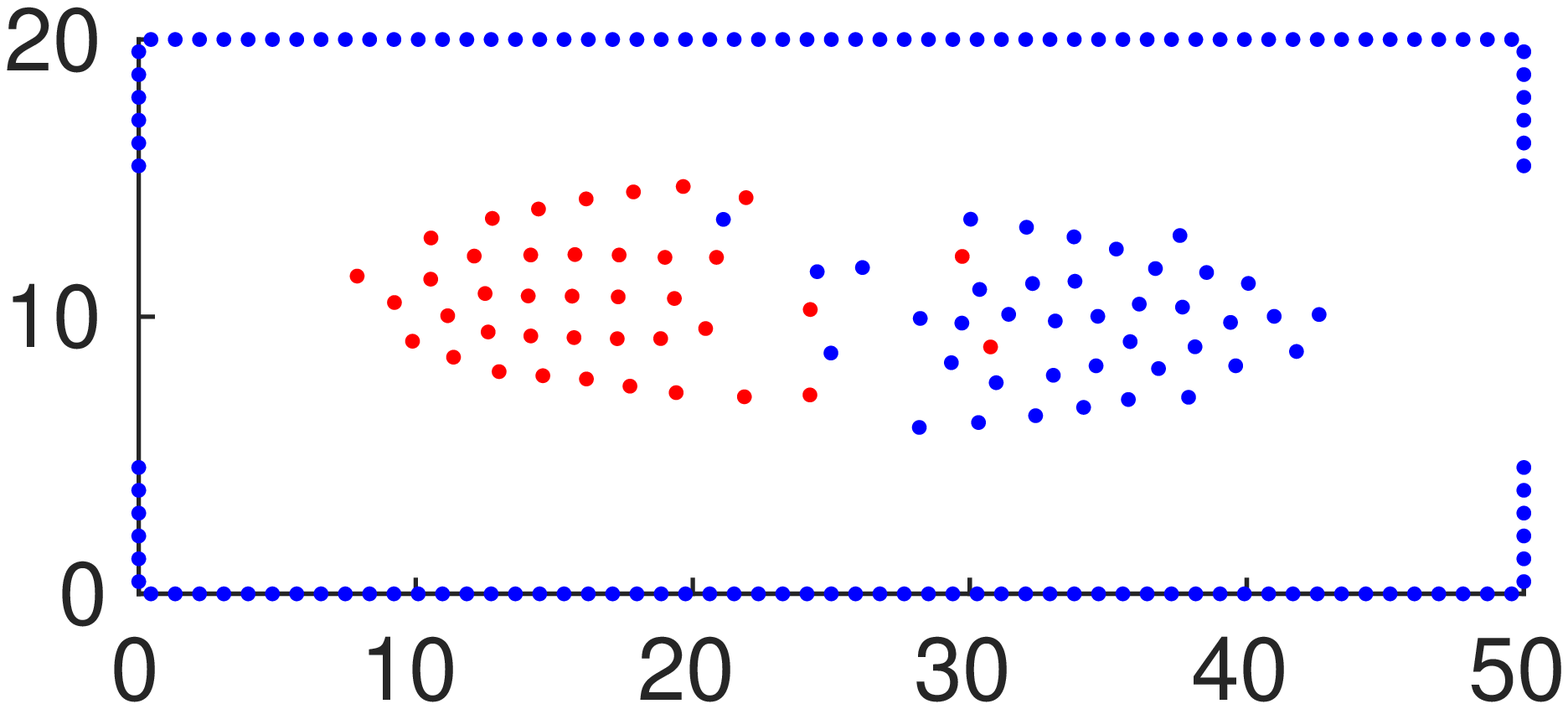}\\		
	\includegraphics[width=.33\linewidth]{sys_with_dir_t_25.eps}\hfill
	\includegraphics[width=.33\linewidth]{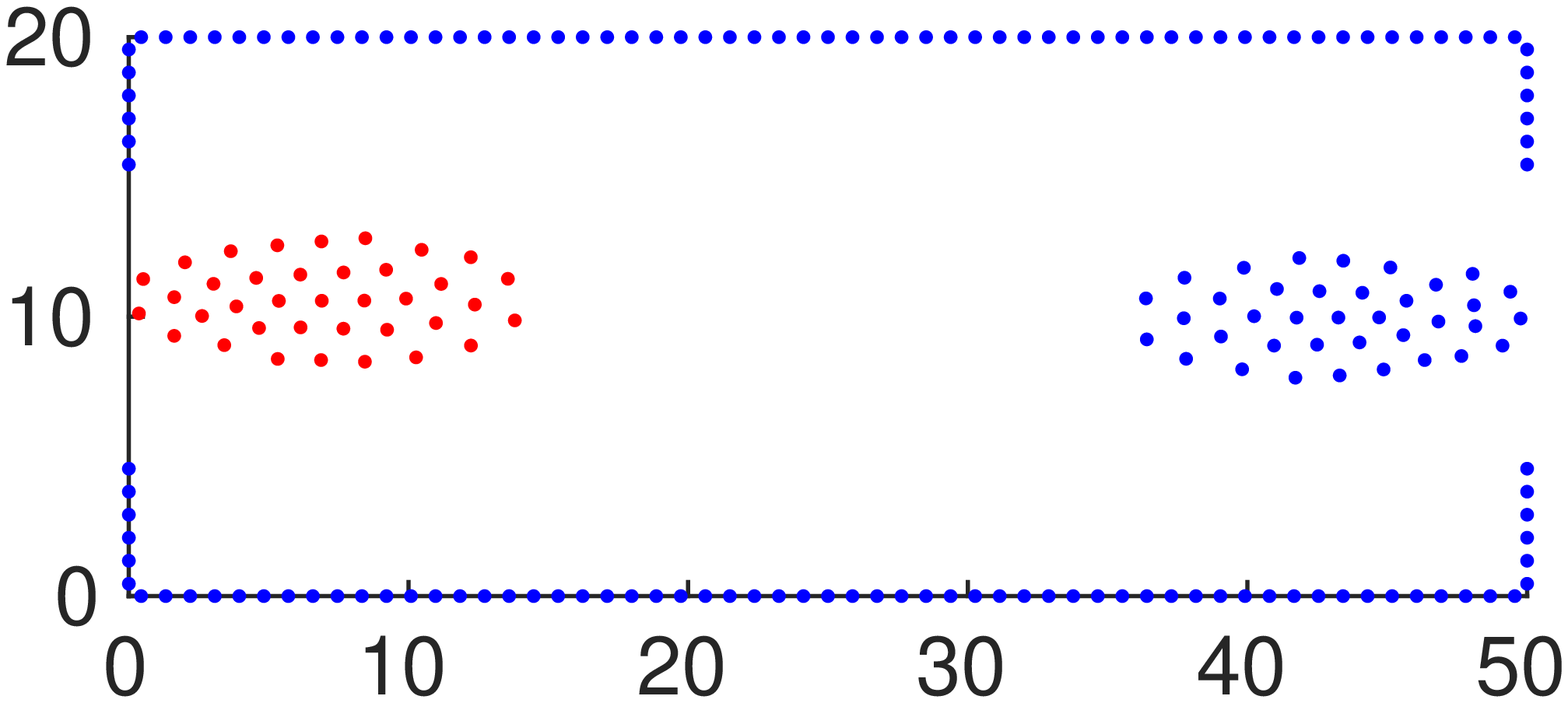}\hfill
	\includegraphics[width=.33\linewidth]{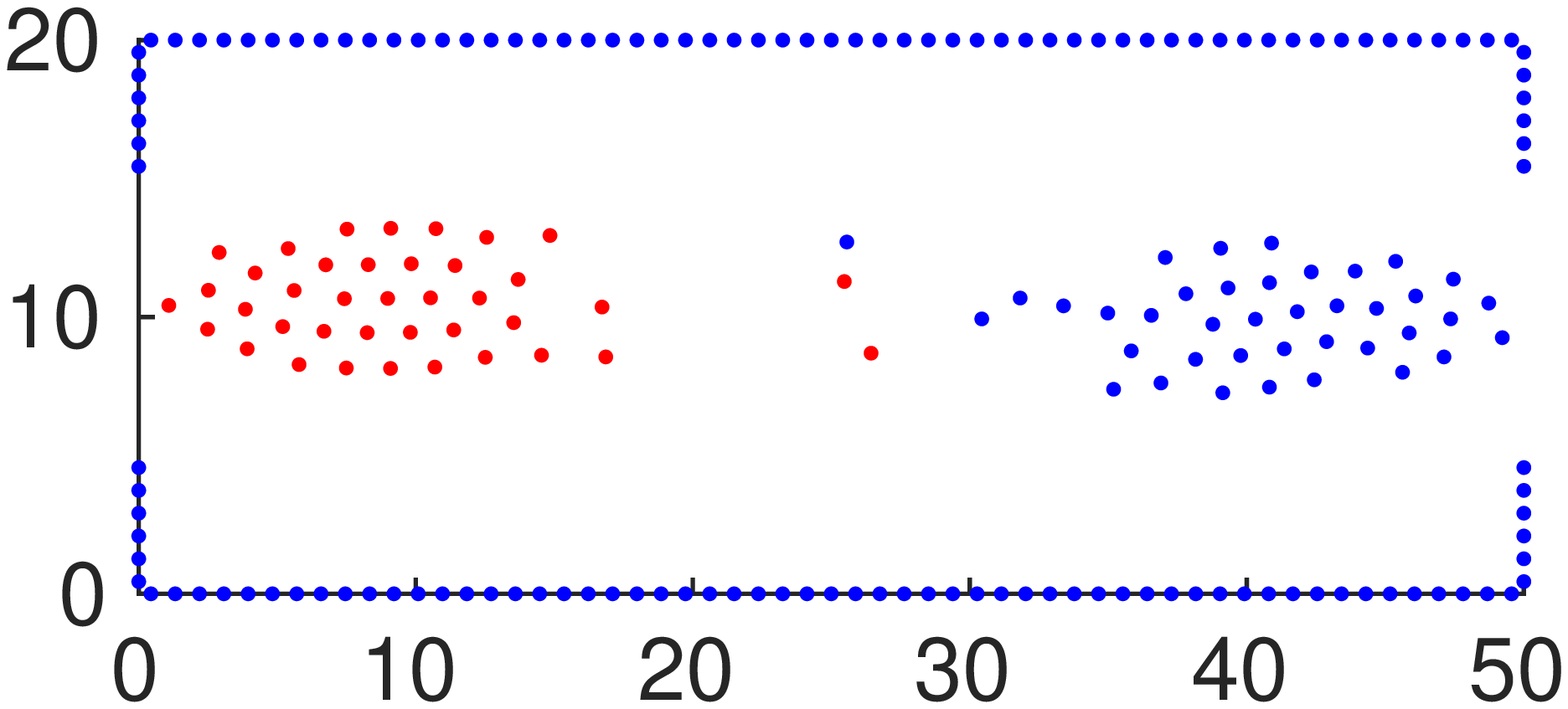}
	
	\caption{Distribution of grid particles in the vision-based model (first column), social force model with $K_n = 100$ (second column) and social force model with $K_n = 1000$ (third column) at time $t = 8, 10, 15, 20, 25$ (top to bottom).}
	\label{fig:11}
\end{figure}

\begin{figure}
	\includegraphics[scale=0.5]{density_with_dir_t_10.eps}\\
	\includegraphics[scale=0.5]{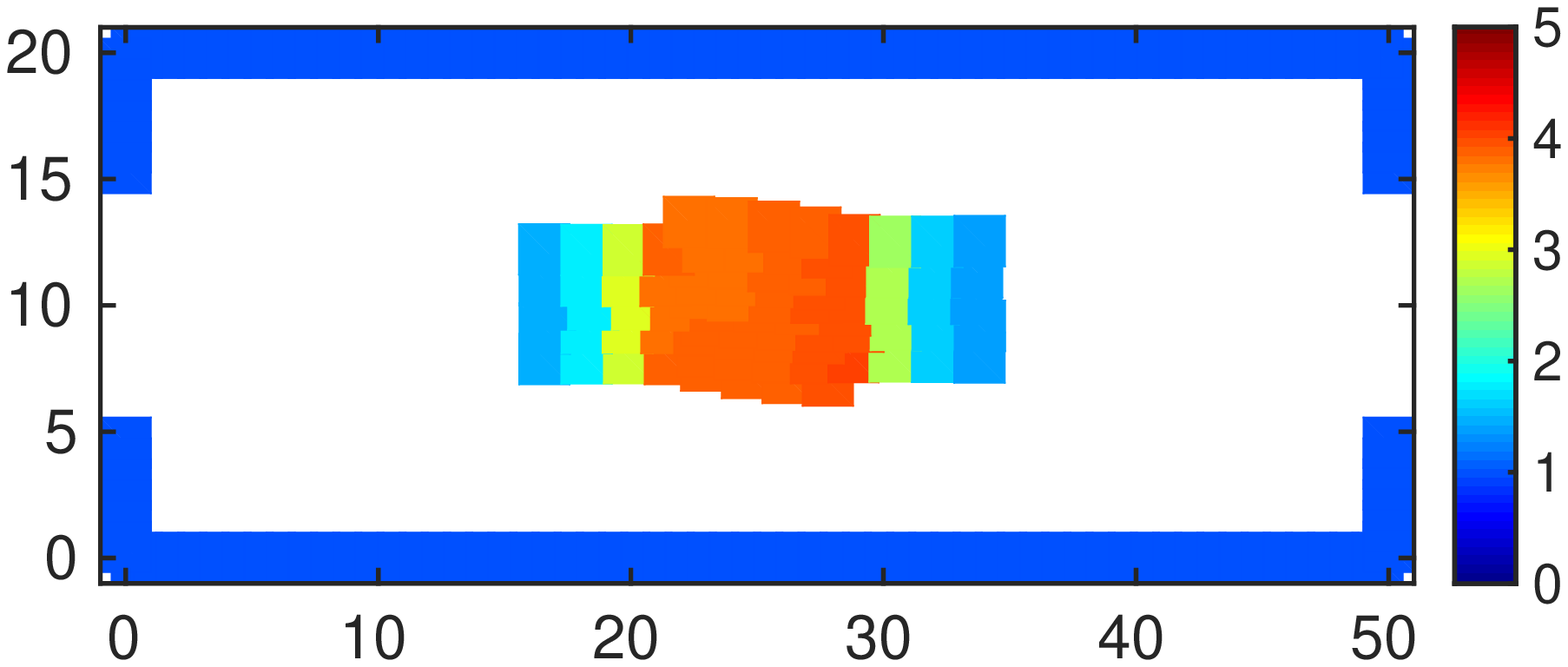}\\
	\includegraphics[scale=0.5]{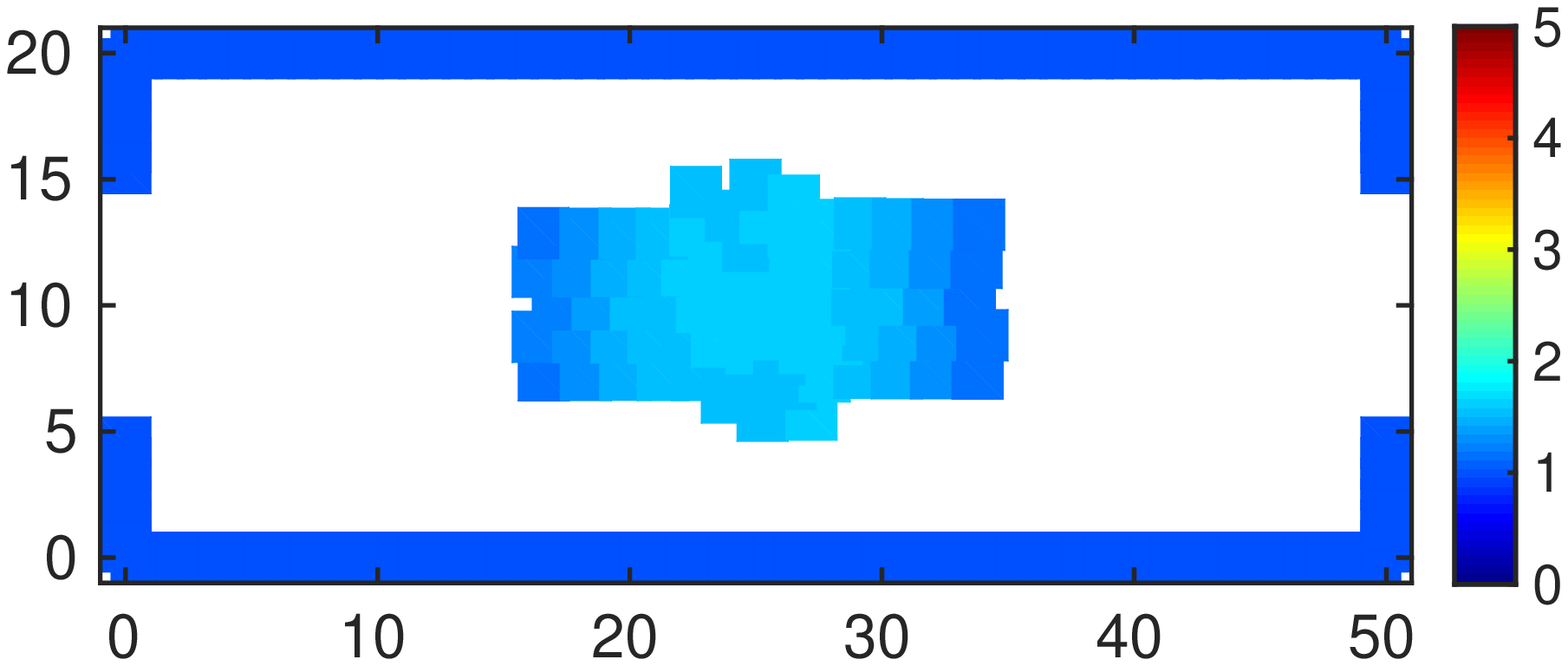}
	\caption{Density of pedestrians in the vision-based model (top), social force model with $K_n = 100$ (middle) and social force model with $K_n = 1000$ (bottom) at time $t = 10$.}
	\label{fig:12} 
\end{figure}

Figure \ref{fig:11} shows the time evolution of the grid particles in vision-based model and social force pedestrian models for different $k_n$ at time $t = 8$, $ t = 10$, $t = 15$, $t = 20$, and $t= 25$. One can observe that collision avoidance between pedestrians of the vision-based model is almost similar with a social force pedestrian model for $k_n = 100$. For bigger values of $k_n$ as in the third column of figure \ref{fig:11} one observes larger differences. Figure \ref{fig:12} shows the corresponding density plots at time $t = 10$ in vision-based model and social force pedestrian models for different $k_n$.

\begin{figure}
	\includegraphics[width=.5\linewidth]{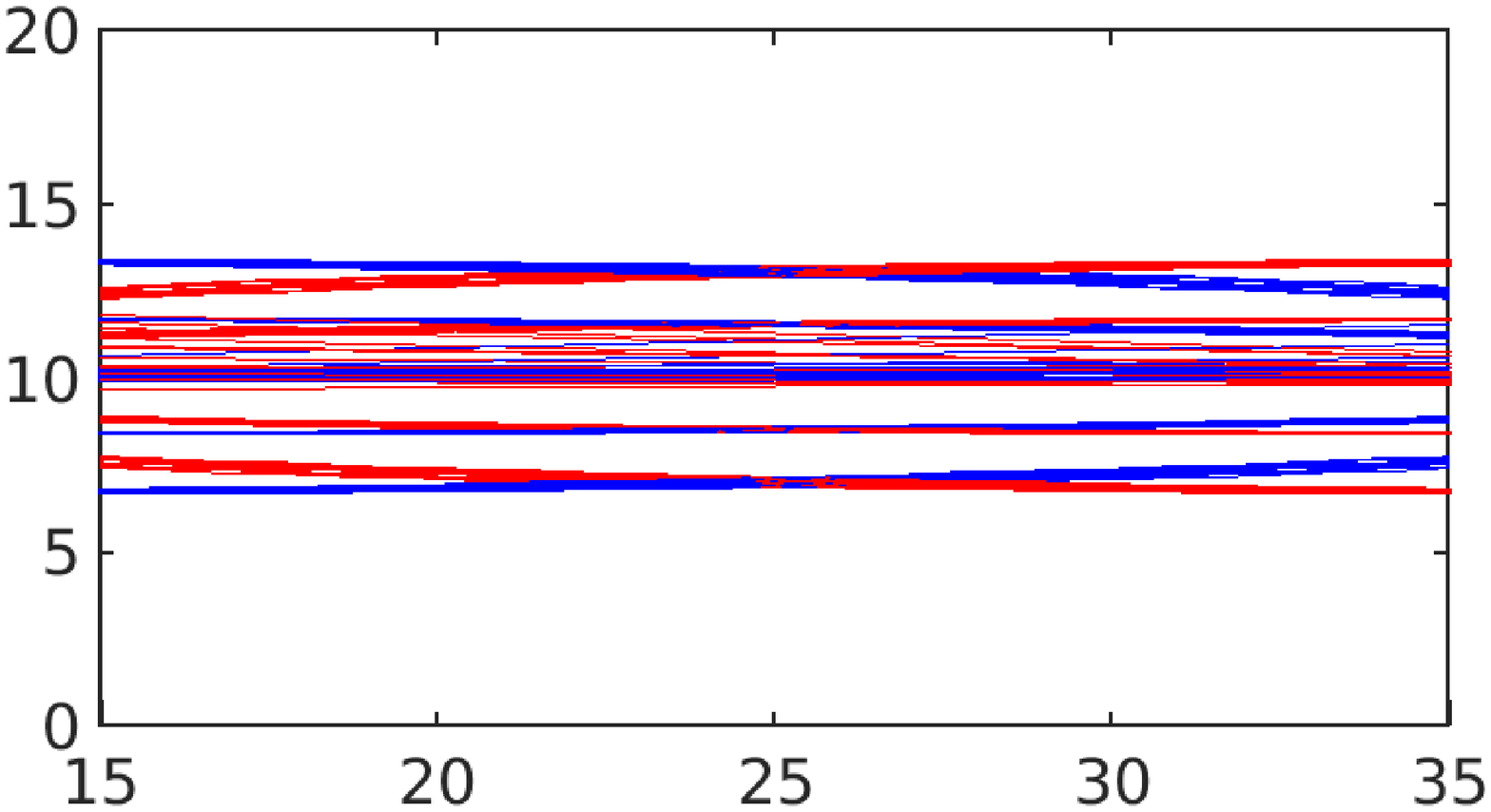}
	\includegraphics[width=.5\linewidth]{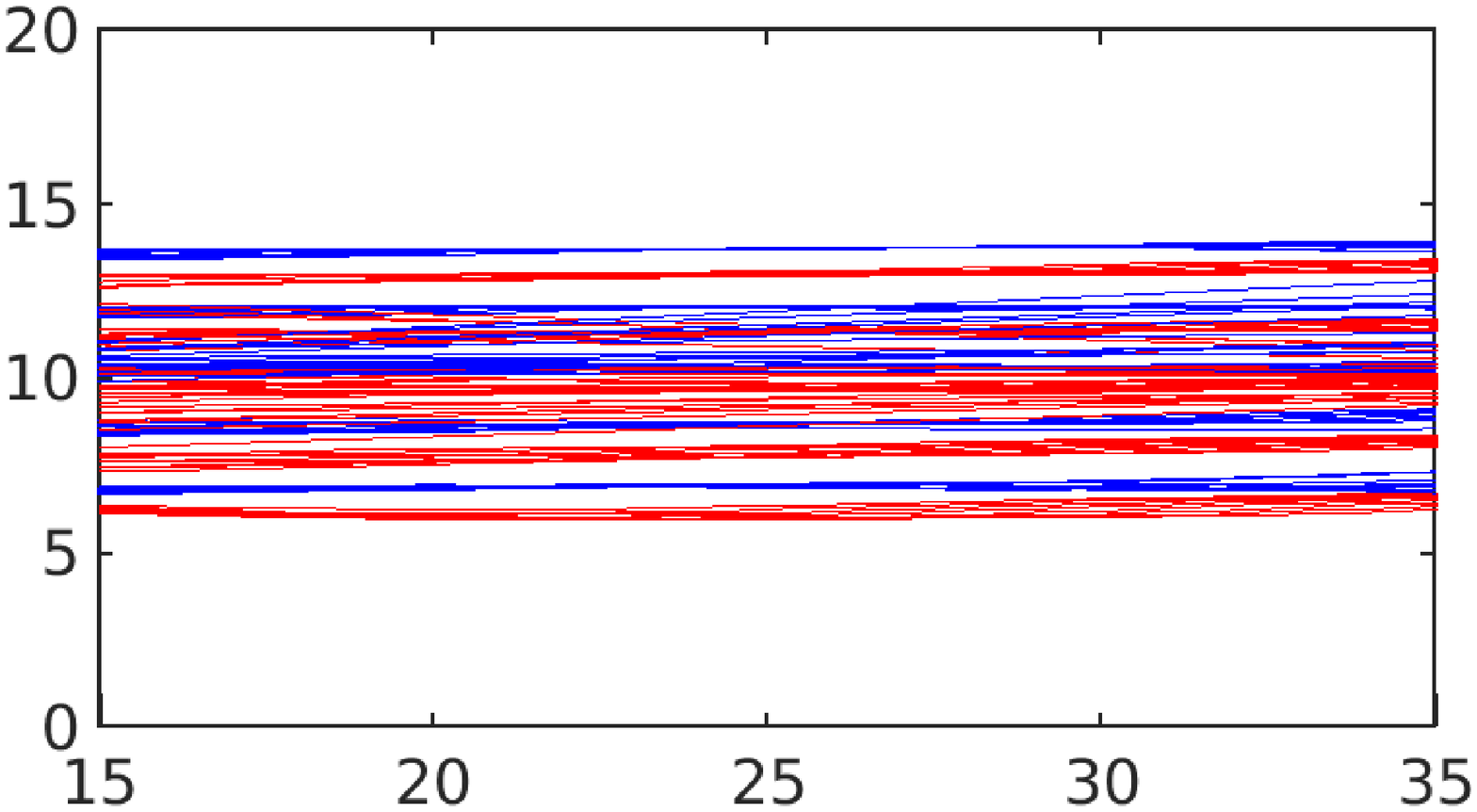}\\
	\includegraphics[width=.5\linewidth]{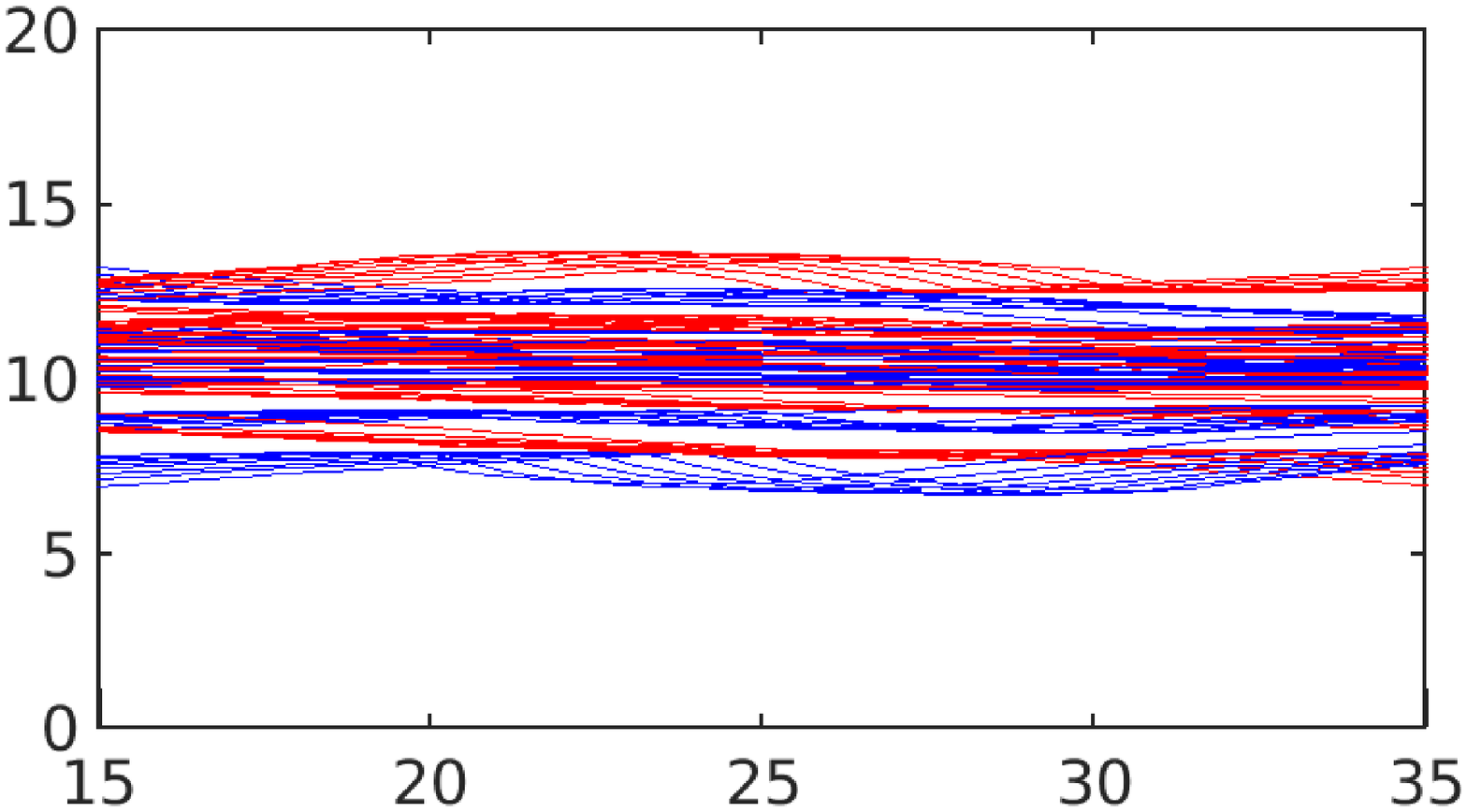}
	\includegraphics[width=.5\linewidth]{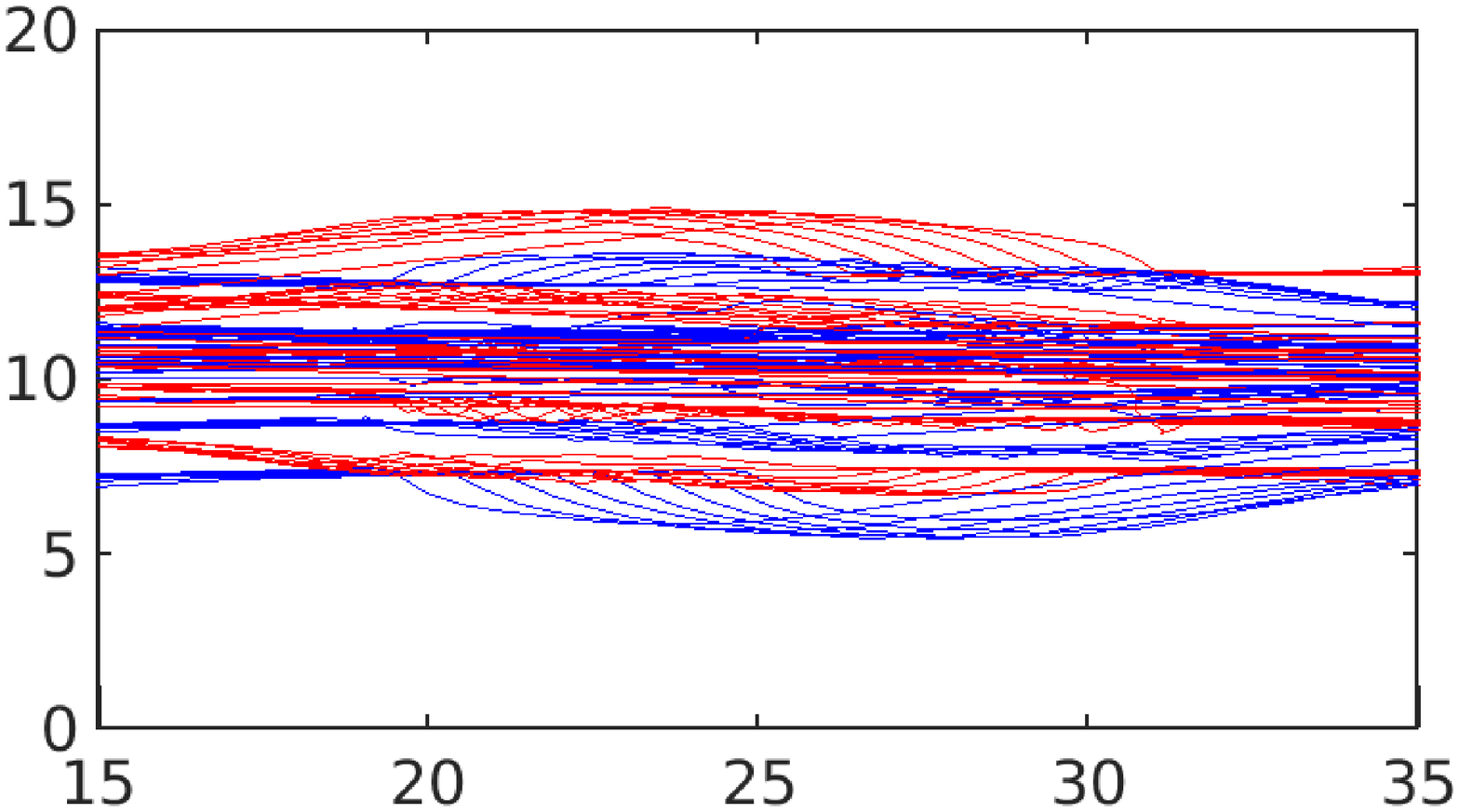}\\
	\includegraphics[width=.5\linewidth]{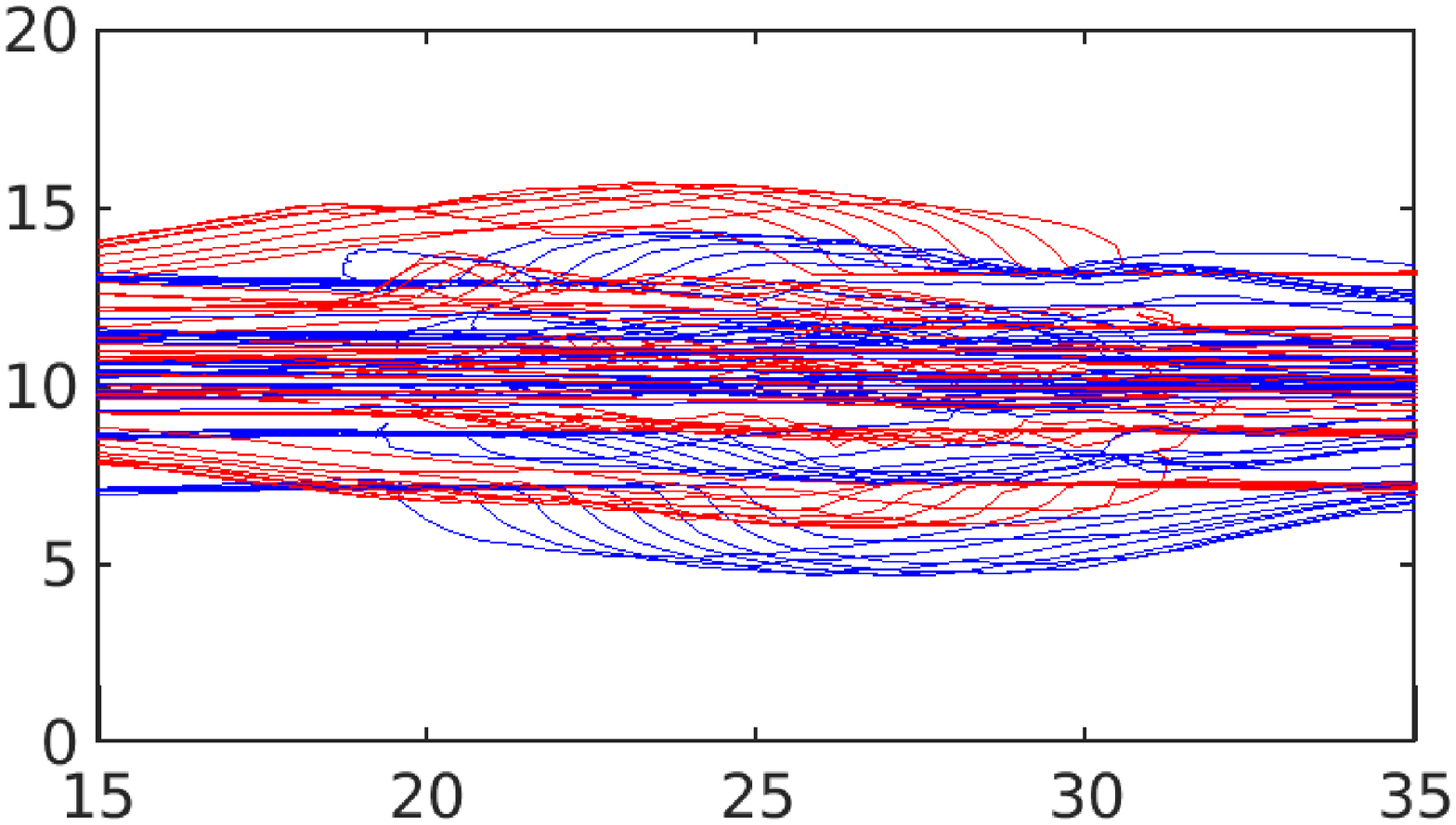}
	\includegraphics[width=.5\linewidth]{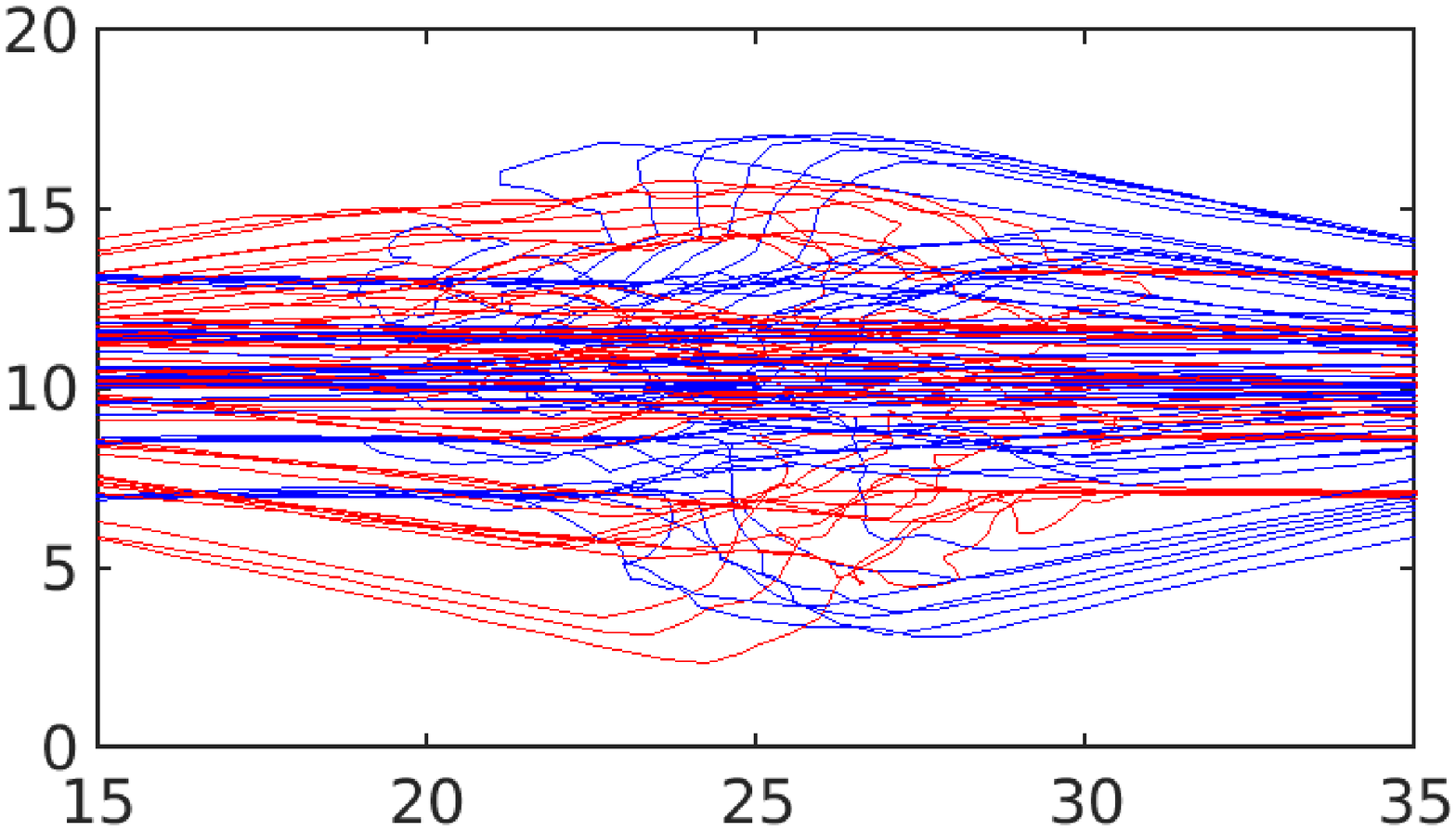}
	
	\caption{Trajectory of the grid particles in the models such as vision-based model without control on direction (top left), vision-based model with control on direction (top right), social force model for $k_n$ being 100 (middle left), 500 (middle right), 1000 (bottom left), 2000 (bottom right) respectively.}
	\label{fig:13}
\end{figure}

Figure \ref{fig:13} shows the trajectory of the grid particles in the vision-based model and social force model. For larger values of $k_n$ in the social force model one observes trajectory which are similar to those obtained in \cite{ondrej}.

\begin{figure}[htbp]
	\centering
	\includegraphics[width=0.8\linewidth]{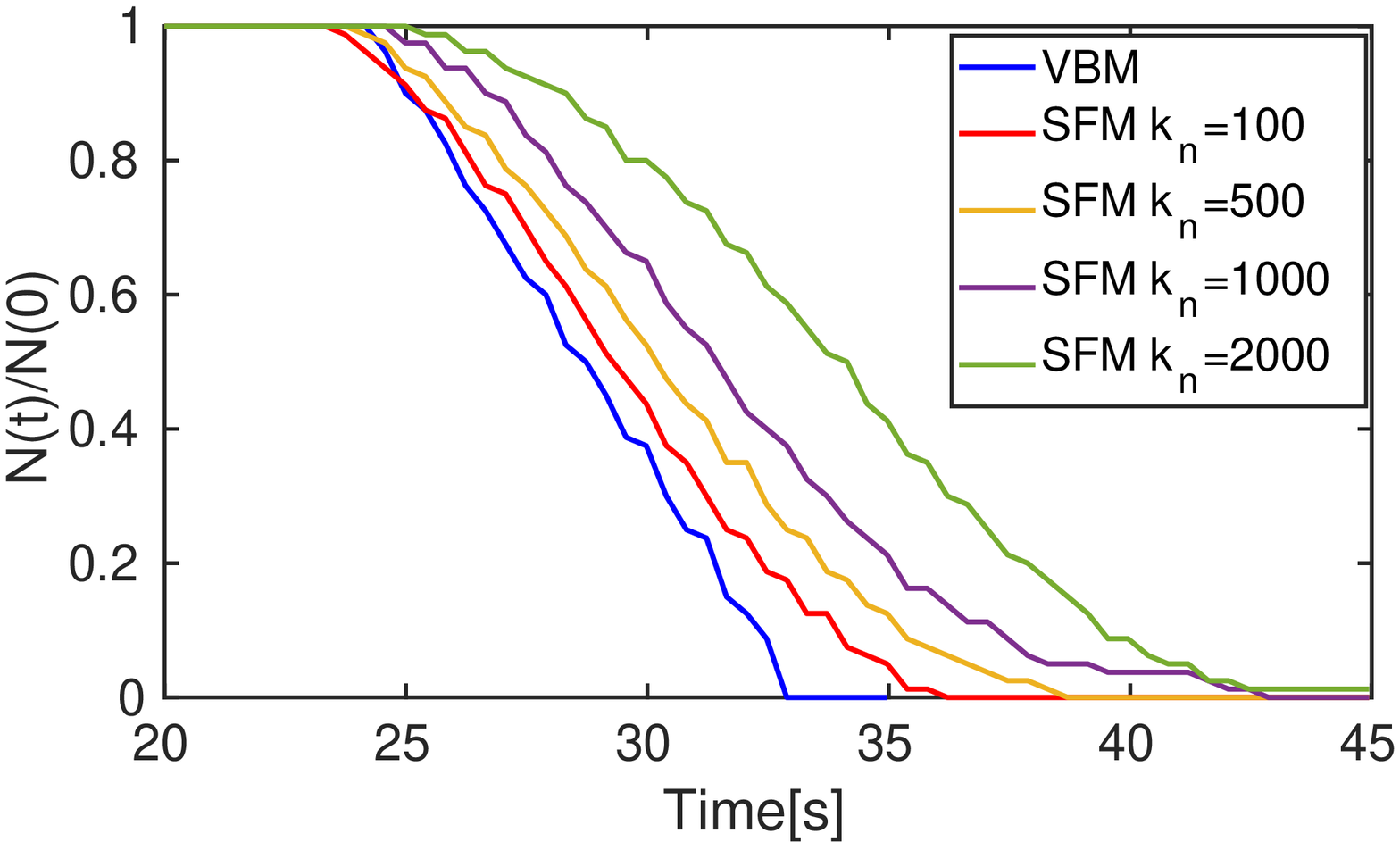}
	\caption{Ratio of initial and actual grid particles with respect to time in vision-based model and social force model with different $k_n$.}
	\label{fig:14}
\end{figure}

Figure \ref{fig:14} shows the percentage of grid particles being in the computational domain for the vision-based model and social force model for the different $k_n$ with respect to time. One observes that the evacuation time is larger in the social force model and also, increasing for  larger values of $k_n$.


\subsection{Computational time}\label{sec:5.5}
The computational time for a simulation of the macroscopic equation up to $t = 45s$ are given for the vision-based model and the social force based model in  table \ref{table:1}. 	One  observes that the social force model has almost double the computational costs than the vision-based model. This is mainly due to the solution of the   Eikonal equation with source depending on  $\rho (t,x)$, which  has to be solved at every time step.

\begin{table}
	\centering
	\begin{tabular}{|c|c|c|c|c|}
		\hline
		Models & RUN TIME\\
		\hline
		VBM without control on direction & 02:00\\	    	  
		\hline
		VBM with control on direction & 02:01\\
		\hline
		Social force model with $k_n = 100$ & 04:23 \\
		\hline
		Social force model with $k_n = 500$ & 04:29 \\
		\hline
		Social force model with $k_n = 1000$ & 04:31 \\
		\hline
		Social force model with $k_n = 2000$ & 04:35 \\
		\hline
	\end{tabular}
	\caption{Computing time of the vision-based model and social force model.}
	\label{table:1}
\end{table}


\section{Conclusion}\label{sec:6}

We have considered the  macroscopic vision-based pedestrian model presented in \cite{degond}.  A mesh-free particle method to solve the governing equations is presented and used for the computation of several numerical examples analysing  non-local and  local vision-based model with and without additional  social force terms. Also, we have presented a comparison between macroscopic models based on vision-based model and on social force pedestrian model  coupled with the eikonal equation for a simple bi-directional situation.  Results obtained here with a mesh free particle method for macroscopic equations are in accordance with the ones obtained for the microscopic vision based model in \cite{ondrej}.
Future research topics are in particular the consideration of more complex situations, for example, interaction with moving objects and comparison with experimental data.


\section*{Acknowledgment}
	This work is supported by the German research foundation, DFG grant KL 1105/27-1 and  by the DAAD PhD program MIC "Mathematics in Industry and Commerce".




\end{document}